\documentclass[12pt]{article}
\usepackage{color, amsmath, amssymb, amsthm, graphicx, hyperref, bbm, tikz, natbib, algorithm, algpseudocode,verbatim}
\usepackage[top=1in, bottom=1in, left=1in, right=1in]{geometry}
\allowdisplaybreaks

\usepackage[toc, page]{appendix}
\usepackage{afterpage}
\usepackage{placeins}

\usepackage{setspace}

\newtheorem{theorem}{Theorem}

\newtheorem{proposition}{Proposition}
\newtheorem{corollary}{Corollary}
\newtheorem{lemma}{Lemma}
\newtheorem{remark}{Remark}
\newtheorem{assumption}{Assumption}

\title{Bootstrapping $\ell_p$-Statistics in High Dimensions}
\author{Alexander Giessing\footnote{Department of ORFE, Princeton University, Princeton, NJ 08544, USA. E-mail: giessing@princeton.edu.} \:\:\: and\:\:\: Jianqing Fan\footnote{Department of ORFE, Princeton University, Princeton, NJ 08544, USA. E-mail: jqfan@princeton.edu.} \thanks{The project is supported by DMS-1662139 and DMS-1712591, NIH grant 2R01-GM072611-13, and ONR grant N00014-19-1-2120.}}
\begin{document}
\maketitle

\begin{abstract}
	This paper considers a new bootstrap procedure to estimate the distribution of high-dimensional $\ell_p$-statistics, i.e. the $\ell_p$-norms of the sum of $n$ independent $d$-dimensional random vectors with $d \gg n$ and $p \in [1, \infty]$. We provide a non-asymptotic characterization of the sampling distribution of $\ell_p$-statistics based on Gaussian approximation and show that the bootstrap procedure is consistent in the Kolmogorov-Smirnov distance under mild conditions on the covariance structure of the data. As an application of the general theory we propose a bootstrap hypothesis test for simultaneous inference on high-dimensional mean vectors. We establish its asymptotic correctness and consistency under high-dimensional alternatives, and discuss the power of the test as well as the size of associated confidence sets. We illustrate the bootstrap and testing procedure numerically on simulated data.
\end{abstract}

\noindent
\hspace{10pt} \quad{} {\small \textbf{Keywords:} Bootstrap; high-dimensional inference; Berry-Esseen bound; anti-concen-
	
\hspace{10pt} tration; Gaussian approximation; Gaussian comparison inequality.}

\section{Introduction}
Let $X=\{X_i\}_{i=1}^n$ be a random sample of independent and centered random vectors in $\mathbb{R}^d$, where dimension $d = d_n$ may grow with sample size $n$. Consider the re-scaled sum
\begin{align*}
S_n^X = \big(S_{n,1}^X, \ldots, S_{n,d}^X\big)' : = \frac{1}{\sqrt{n}}\sum_{i=1}^n X_i,
\end{align*}
and define the \emph{$\ell_p$-statistic} $T_{n,p}$ by
\begin{align}\label{eq:Introduction-1}
T_{n,p} := \|S_n^X\|_p =\begin{cases}
\left(\sum_{k=1}^d |S_{n,k}^X|^p\right)^{1/p}, & p \in [1, \infty)\\
\: \max_{1 \leq k \leq d} |S_{n,k}^X|, & p = \infty.
\end{cases}
\end{align}
This paper is concerned with developing a bootstrap procedure to estimate the distribution of $\ell_p$-statistics when the dimension $d$ exceeds the sample size $n$. This distribution is of interest in many statistical applications. In particular, $\ell_2$-statistic $T_{n,2}$ and maximum statistic $T_{n,\infty}$ are frequently applied to a broad spectrum of statistical problems such as testing of multiple means, construction of simultaneous confidence regions, and model selection \citep{bai1996effect,chen2010two,fan2015power}.

In low dimensions, when the dimension $d$ is fixed, the asymptotic properties of $\ell_p$-statistics are well-understood: If the data are i.i.d. with finite second moments, the central limit theorem (CLT) applied to the re-scaled sum $S_n^X$ and the continuous mapping theorem guarantee that $T_{n,p}\overset{d}{\rightarrow} \|Z\|_p$, where $Z \sim N\left(0, \mathrm{E}\left[X_1X_1'\right]\right)$. Thus, the limiting distribution of $T_{n,p}$ depends on the data only through the first two moments. Closed-form expressions of the limiting distribution of $\ell_p$-statistics remain somewhat elusive, but for $T_{n,2}$ and $T_{n,\infty}$ tractable characterizations exist under additional assumptions on the covariance structure.

The situation is very different in high dimensions. If the dimension $d$ grows faster than $\sqrt{n}$, the classical CLT does no longer apply to the re-scaled sum $S_n^X$. So, to approximate the distribution of $T_{n,p}$ one has to target directly the scalar random variable $\|S_n^X\|_p$. Since $\|S_n^X\|_p$ is a highly non-linear function of the random sample $X$, this calls for a non-parametric approach. In this direction,~\cite{chernozhukov2013GaussianApproxVec, chernozhukov2015ComparisonAnti, chernozhukov2017CLTHighDim} have made important progress by developing a non-parametric multiplier bootstrap procedure to approximate the distribution of the maximum statistic $T_{n,\infty}$. In this paper, we further develop this line of research. While there exist specialized results for high-dimensional sum-of-squares type  $T_{n,2}^2$-statistics~\citep{bai1996effect,bentkus2003DependenceBerryEsseen, chen2010two, fan2015power, pouzo2015bootstrap, xu2019L2Asymptotics}, a unified investigation on the weak convergence of general $\ell_p$-statistics remains highly challenging due to the lack of smoothness of the $\ell_p$-norm. In this sense, our work solves a long-standing open problem initiated by the aforementioned pioneering work.

The primary methodological contribution of this paper is a bootstrap procedure for $\ell_p$-statistics with $p \in [1, \infty]$. Our bootstrap procedure draws inspiration from above observation that in low dimensions the limiting distribution of $\ell_p$-statistics depends only on the first two moments of the data. Specifically, the proposed algorithm involves sampling bootstrap data from a Gaussian distribution that is parameterized by an estimate of the covariance matrix. The algorithm works with any estimate of the covariance matrix; it is easy to implement and very versatile. In particular, it can be combined with estimates of the covariance matrix that leverage special structures such as low rank, (approximate) sparsity or bandedness. The algorithm is best understood as a hybridization of non-parametric and parametric bootstrap, and we call it the \emph{Gaussian parametric bootstrap}.

A secondary methodological contribution is a bootstrap hypothesis test for testing many linear restrictions on high-dimensional mean vectors. This hypothesis test is based on the Gaussian parametric bootstrap for $\ell_p$-statistics and it is asymptotically correct and consistent under certain high-dimensional alternatives. We give precise recommendations on how to choose the exponent $p \in [1, \infty]$ based on characteristics of the random sample (tails and covariance structure) and to maximize the power for given alternative hypotheses. For small exponents $p$, the test is useful when the goal is to identify significant subsets from a large collection of means, e.g. sets of genes in micro-array and genetic sequence studies. Whereas for large exponents $p$, the test is powerful when the purpose is to detect significant singletons, e.g. anomaly detection in materials science and medical imaging.

The two main theoretical contributions of this paper are a non-asymptotic characterization of the sampling distribution of $\ell_p$-statistics $T_{n,p}$ in high dimensions and the consistency of the Gaussian parametric bootstrap. The non-asymptotic characterization is based on a Gaussian approximation, i.e. a proxy statistic constructed from Gaussian random vectors. The quality of the Gaussian and bootstrap approximation improves as the sample size $n$ increases and shows a subtle interplay between dimension $d$, exponent $p$, and the tail distribution of the data. Among other things, we demonstrate that if the data has light tails the approximation errors vanish for $ \log d = o(n)$ and all $p \in [1, \infty]$; whereas if the data is heavy-tailed with at most $s \geq 4$ finite moments the approximation errors are negligible for $d \log d= o(n^{s/4})$ and all $p \in [1,  s]$. These theoretical results provide a comprehensive view on the asymptotic distribution theory of $\ell_p$-statistics and are relevant in guiding practitioners in choosing between different $\ell_p$-statistics given the properties of the random sample at hand. Qualitatively, our numerical experiments lend further support to these theoretical findings.

Establishing the Gaussian approximation and consistency of the bootstrap is non-trivial and we develop a significant amount of new technical tools. The following three technical results are of interest beyond the scope of this paper: First, we derive an abstract Berry-Esseen-type CLT for $\ell_p$-statistics in high dimensions, which extends and improves the known Berry-Esseen-type CLTs for $p=2$~\citep{bentkus2003DependenceBerryEsseen} and $p = \infty$~\citep{chernozhukov2017CLTHighDim}. Second, we establish an anti-concentration inequality for $\ell_p$-norms of random vectors with log-concave probability measure. For $p \in \{2, \infty\}$ this inequality is sharper than related inequalities by~\cite{goetze2019LargeBall} and~\cite{chernozhukov2017Nazarov}. Third, we develop a Gaussian comparison inequality to compare the distributions of $\ell_p$-norms of different Gaussian random vectors in Kolmogorov-Smirnov distance. For $p = \infty$ this inequality improves the corresponding result in~\cite{chernozhukov2015ComparisonAnti}.

\emph{Organization.} The paper is organized as follows. We introduce the Gaussian parametric bootstrap in Section~\ref{sec:Methodology} and present our main theoretical results on the Gaussian approximation of $\ell_p$-statistics and the consistency of the Gaussian parametric bootstrap in Section~\ref{sec:MainResults}. We develop applications to testing high-dimensional mean vectors in Section~\ref{sec:SimultaneousTesting} and report results from several numerical experiments in Section~\ref{sec:NumericalExperiments}. In Appendix~\ref{apx:MainTheory} we discuss technical results, including the abstract Berry-Esseen-type CLT, the anti-concentration inequalities for $\ell_p$-statistics, and the new Gaussian comparison theorems. Appendix~\ref{apx:Proofs} contains proofs to all our results.

\emph{Notation.} For non-negative real-valued sequences $\{a_n\}_{n\geq 1}$ and $\{b_n\}_{n \geq 1}$, the relation $a_n \lesssim b_n$ means that there exists an absolute constant $c > 0$ independent of $n, d, p$ and an integer $n_0 \in \mathbb{N}$ such that $a_n \leq c b_n$ for all $n \geq n_0$. We write $a_n \asymp b_n$ if  $a_n \lesssim b_n$ and $b_n \lesssim a_n$. We define $a_n \vee b_n = \max\{a_n, b_n\}$  and $a_n \wedge b_n = \min\{a_n, b_n\}$. For a vector $a \in \mathbb{R}^d$ and $p \in [1, \infty)$ we write $\|a\|_p = (\sum_{k=1}^d |a_k|^p)^{1/p}$. Also, we write $\|a\|_\infty = \max_{1 \leq k \leq d} |a_k|$. For a scalar random variable $\xi$ and $\alpha \in (0, 2]$ we define the $\psi_\alpha$-Orlicz norm by $\|\xi \|_{\psi_\alpha} = \inf\{t > 0 : \mathrm{E}[\exp(|\xi|^\alpha/t^\alpha)] \leq 2\}$. For a sequence of scalar random variables $\{\xi_n\}_{n \geq 1}$ we write $\xi_n = O_p(a_n)$ if $\xi_n/a_n$ is stochastically bounded. For any symmetric real-valued matrix $M \in \mathbb{R}^{d \times d}$ we denote its largest and smallest eigenvalue by $\lambda_{\max}(M)$ and $\lambda_{\min}(M)$, respectively. We denote its operator norm by $\|M\|_{op}$ (its largest singular value) and $\|M\|_{2 \rightarrow p} = \sup_{\|u\|_2 \leq 1} \|Mu\|_p$. We write $M \succeq 0$ to indicate that $M$ is positive semi-definite. For any convex body $K \subset \mathbb{R}^d$ we write $\mathrm{Vol}(K) = \int_K d \lambda^d$, where $\lambda^d$ is the Lebesgue measure in $d$ dimensions.

\section{Methodology}\label{sec:Methodology}
We introduce the new Gaussian parametric bootstrap for $\ell_p$-statistics and discuss its relation to the non-parametric Gaussian multiplier bootstrap.

\subsection{Gaussian parametric bootstrap}\label{subsec:GPB}
Let $X=\{X_i\}_{i=1}^n$ be a random sample of independent and centered random vectors. The Gaussian parametric bootstrap algorithm requires as input a consistent and positive semi-definite estimate $\widehat{\Sigma}_n$ of the (averaged) population covariance matrix,
\begin{align*}
\Sigma_n : = \mathrm{E}\left[\frac{1}{n} \sum_{i=1}^n X_iX_i'\right].
\end{align*}
We will discuss candidates for $\widehat{\Sigma}_n$ in subsequent sections. Let $V^X \mid X \sim N(0, \widehat{\Sigma}_n)$ and define the Gaussian parametric bootstrap estimate of the $\ell_p$-statistic $T_{n,p}$ by
\begin{align}\label{eq:sec:Methodology-1}
T_{n,p}^* : = \left\| V^X\right\|_p = \begin{cases}
\left(\sum_{k=1}^d \left|V^X_k\right|^p\right)^{1/p}, &  p \in [1, \infty)\\
\: \max_{1 \leq k \leq d} |V_k^X|, & p = \infty.
\end{cases}
\end{align}

The rationale for this bootstrap statistic is easiest to understand in low dimensions: If dimension $d$ is fixed and the data $X= \{X_i\}_{i=1}^n$ is i.i.d. with finite second moments, the CLT and the continuous mapping theorem imply that $T_{n,p}\overset{d}{\rightarrow} \|Z\|_p$, where $Z \sim N\left(0, \mathrm{E}\left[X_1X_1'\right]\right)$. Hence, in this scenario, the bootstrap statistic $T_{n,p}^*$ is just the parametric bootstrap estimate of the limiting random variable $\|Z\|_p$. Of course, if $d \geq \sqrt{n}$ and the data is non-identically distributed, the CLT does not apply and the limiting random variable $Z$ needs not to exist. The gist of the theoretical results in Sections~\ref{subsec:GaussianApprox} and~\ref{subsec:BootstrapConsistency} is that we do not need the CLT to hold for the distributions of $T_{n,p}^*$ and $T_{n,p}$ to be close. For $\ell_p$-statistics this result is new, but it is in line with similar results on linear regression functions, empirical processes in infinite-dimensional Banach spaces, as well as maximum and spectral statistics in high dimensions~\cite[e.g.][]{bickel1983bootstrapping, radulovic1998can, chernozhukov2013GaussianApproxVec, roellin2013stein, lopes2019BootstrappingSpectral}.

\subsection{Relation to the Gaussian multiplier bootstrap}\label{subsec:GMB}
The Gaussian multiplier bootstrap was first proposed by~\cite{chernozhukov2013GaussianApproxVec} in the context of the maximum statistic $T_{n,\infty}$. It is a special case of the wild bootstrap method~\citep[][]{wu1986jacknife, liu1988bootstrap, mammen1993bootstrap} and its adaptation to general $\ell_p$-statistics $T_{n,p}$ is straightforward:

Let $g= \{g_i\}_{i=1}^n$ be a sequence of i.i.d. standard normal random variables independent of the random sample $X = \{X_i\}_{i=1}^n$. The Gaussian multiplier bootstrap algorithm builds on the centered random sample $X_1 - \bar{X}_n, \ldots, X_n - \bar{X}_n$, where $\bar{X}_n := n^{-1}\sum_{i=1}^n X_i$. We set
\begin{align*}
S^{gX}_n  :=  \left(S^{gX}_{n1}, \ldots, S^{gX}_{nd}\right)' := \frac{1}{\sqrt{n}}\sum_{i=1}^n g_i (X_i - \bar{X}_n),
\end{align*}
and define the Gaussian multiplier bootstrap estimate of the $\ell_p$-statistic $T_{n,p}$ by
\begin{align}\label{eq:subsec:BootstrapProcedures-1}
T_{n,p}^{g} : = \left\|S^{gX}_n\right\|_p
\end{align}

Since $S_n^{gX} \mid X \sim N(0, \widehat{\Sigma}_{\mathrm{naive}})$ with $\widehat{\Sigma}_{\mathrm{naive}} = n^{-1}\sum_{i=1}^n(X_i - \bar{X}_n)(X_i - \bar{X}_n)'$, the Gaussian multiplier bootstrap statistic is in fact equivalent to a Gaussian parametric bootstrap statistic based on the sample covariance matrix $\widehat{\Sigma}_{\mathrm{naive}}$. The key advantage of the Gaussian parametric over the Gaussian multiplier bootstrap is that it allows for more refined estimates of the population covariance matrix $\Sigma_n$ that leverage additional structure such as low-rank, (approximate) sparsity, and bandedness. This is particularly important in high dimensions where the sample covariance matrix $\widehat{\Sigma}_{\mathrm{naive}}$ is a poor estimate of the population covariance matrix.


\section{Theoretical analysis}\label{sec:MainResults}
We present a non-asymptotic characterization of $\ell_p$-statistics via Gaussian approximation and establish the consistency of the Gaussian parametric bootstrap procedure.

\subsection{Assumptions}\label{subsec:Assumptions}
Unless otherwise stated, $X=\{X_i\}_{i=1}^n$ denotes a random sample of independent and centered random vectors in dimension $d$, where $d = d_n$ grows with the sample size $n$. We analyze the theoretical properties of $\ell_p$-statistics and the Gaussian parametric bootstrap under the following three different assumptions on the tails of random vectors.

\begin{assumption}[Sub-Gaussian]\label{assumption:SubGaussian}
	Let $X = \{X_i\}_{i=1}^n$ be a sequence of independent and centered random vectors in $\mathbb{R}^d$ such that for all $1 \leq i \leq n$,
	\begin{align*}
	\forall u \in \mathbb{R}^d :\:\: \left\| u'X_i \right\|_{\psi_2} \lesssim \mathrm{E}\left[(u'X_i)^2\right]^{1/2}.
	\end{align*}
\end{assumption}

\begin{assumption}[Sub-Exponential]\label{assumption:SubExponential}
	Let $X = \{X_i\}_{i=1}^n$ be a sequence of independent and centered random vectors in $\mathbb{R}^d$ such that for all $1 \leq i \leq n$,
	\begin{align*}
	\forall u \in \mathbb{R}^d :\:\: \left\| u'X_i \right\|_{\psi_1} \lesssim \mathrm{E}\left[(u'X_i)^2\right]^{1/2}.
	\end{align*}
\end{assumption}

\begin{assumption}[Finite $s$th moments]\label{assumption:FiniteMoments}
	Let $X = \{X_i\}_{i=1}^n$ be a sequence of independent and centered random vectors in $\mathbb{R}^d$ such that for some $s \geq 3$ and all $1 \leq i \leq n$,
	\begin{align*}
	\forall u \in \mathbb{R}^d :\:\: \mathrm{E}\left[|u'X_i|^s\right]^{1/s} \lesssim K_s \mathrm{E}\left[(u'X_i)^2\right]^{1/2}.
	\end{align*}
\end{assumption}

Assumption~\ref{assumption:SubGaussian} is one of the many equivalent definitions of sub-Gaussian random vectors~\citep[e.g][]{antonini1997subgaussian, vershynin2018HighDimProb}. This specific formulation is useful for applications in high-dimensional statistics because $\mathrm{E}\left[(u'X_i)^2\right]^{1/2} \leq \|\mathrm{E}[X_iX_i']\|_{op}\|u\|_2$. Hence, we can easily incorporate characteristics of the covariance matrix such as sparsity, bandedness, low-rank, etc. Assumptions~\ref{assumption:SubExponential} and~\ref{assumption:FiniteMoments} relax and generalize Assumption~\ref{assumption:SubGaussian} in an obvious way. Most importantly, if $X$ satisfy Assumption~\ref{assumption:FiniteMoments} for all $s \geq 1$ and with $K_s = \sqrt{s}$ ($K_s = s$) then $X$ is sub-Gaussian (sub-Exponential) and also satisfy Assumption~\ref{assumption:SubGaussian} (Assumption~\ref{assumption:SubExponential}).

\subsection{Gaussian approximation}\label{subsec:GaussianApprox}
In this section we show that the distribution of the $\ell_p$-statistic $T_{n,p}$ can be approximated by the distribution of a proxy statistic based on Gaussian random vectors. This result rationalizes the Gaussian parametric bootstrap procedure in high dimensions. It is also relevant for establishing bootstrap consistency in the next section.

Let $Z =\{Z_i\}_{i=1}^n$ be a sequence of independent multivariate Gaussian random vectors $Z_i \sim N(0,\mathrm{E}[X_iX_i'])$ which are independent of $X = \{X_i\}_{i=1}^n$. We define the \emph{Gaussian proxy statistic} of the $\ell_p$-statistic $T_{n,p}$ as
\begin{align}\label{eq:subsec:GaussianApprox-1}
\widetilde{T}_{n,p} = \|S_n^Z\|_p, \hspace{30pt} \mathrm{where} \hspace{30pt} S_n^Z : = \frac{1}{\sqrt{n}}\sum_{i=1}^nZ_i.
\end{align}

To state the Gaussian approximation result we need to define the following additional quantities: the rank of the (averaged) covariance matrices of the $X_i$'s,
\begin{align}\label{eq:subsec:GaussianApprox-2}
r_n := \mathrm{rank}\left(\mathrm{E}\left[\frac{1}{n}\sum_{i=1}^n X_iX_i'\right]\right),
\end{align}
the smallest and largest variances of the $X_i$'s,
\begin{align}\label{eq:subsec:GaussianApprox-3}
\sigma_{n,\min}^2 :=   \min_{1\leq k \leq d} \min_{1 \leq i \leq n}\mathrm{E}\left[ X_{ik}^2\right]\hspace{30pt} \mathrm{and} \hspace{30pt} \sigma_{n,\max}^2 :=   \max_{1 \leq k \leq d} \max_{1 \leq i \leq n}\mathrm{E}[X_{ik}^2],
\end{align}
and the largest ratio of the variances of the $X_i$'s,
\begin{align}\label{eq:subsec:GaussianApprox-4}
\kappa_n^2 := \max_{1 \leq k \leq d} \left(\max_{1 \leq i \leq n} \mathrm{E}[X_{ik}^2] \Big/  \min_{1 \leq i \leq n} \mathrm{E}[X_{ik}^2] \right).
\end{align}

Our first theorem shows that the distribution of $\widetilde{T}_{n,p}$ can approximate the distribution of $T_{n,p}$ in Kolmogorov-Smirnov distance uniformly over all $p \in [1, \infty]$.

\begin{theorem}[Gaussian approximation]\label{theorem:GaussianApprox}
	\noindent	
	\begin{itemize}
		\item[(i)] For all $p \in [1, \infty)$ and $X$ satisfying Assumption~\ref{assumption:SubGaussian},
		\begin{align}\label{eq:theorem:GaussianApprox-1}
		\sup_{t \geq 0} \left|\mathrm{P}(T_{n,p} \leq t) - \mathrm{P}(\widetilde{T}_{n,p} \leq t)\right|\lesssim \sqrt{\frac{p^3 (\log d) r_n^{1/p}}{n^{1/3}}\frac{ \sigma_{n,\max}^2}{\sigma_{n,\min}^2}}.
		\end{align}
		\item[(ii)] For all $p \in [1, \infty)$ and $X$ satisfying Assumption~\ref{assumption:SubExponential},
		\begin{align}\label{eq:theorem:GaussianApprox-2}
		\sup_{t \geq 0} \left|\mathrm{P}(T_{n,p} \leq t) - \mathrm{P}(\widetilde{T}_{n,p} \leq t)\right|\lesssim \sqrt{\frac{p^3 (\log d)^2r_n^{1/p}}{n^{1/3}}\frac{\sigma_{n,\max}^2}{\sigma_{n,\min}^2}}.
		\end{align}
		\item[(iii)] For all $p \in [\log d, \infty]$ and $X$ satisfying either Assumption~\ref{assumption:SubGaussian} or~\ref{assumption:SubExponential},
		\begin{align}\label{eq:theorem:GaussianApprox-3}
		\sup_{t \geq 0} \left|\mathrm{P}(T_{n,p} \leq t) - \mathrm{P}(\widetilde{T}_{n,p} \leq t)\right|\lesssim \left(\frac{\kappa_n^2 \log^7 d}{n}\right)^{1/6}.
		\end{align}
		\item[(iv)] For $X$ satisfying Assumption~\ref{assumption:FiniteMoments} with $s \geq 4$ and all $p \in [1, s]$,
		\begin{align}\label{eq:theorem:GaussianApprox-4}
		\sup_{t \geq 0} \left|\mathrm{P}(T_{n,p} \leq t) - \mathrm{P}(\widetilde{T}_{n,p} \leq t)\right| \lesssim (K_s \vee \sqrt{s})\sqrt{\frac{p^3 d^{4/(3s)}r_n^{1/p}}{n^{1/3}}\frac{\sigma_{n,\max}^2}{\sigma_{n,\min}^2} }.
		\end{align}
	\end{itemize}
\end{theorem}

\begin{remark}
	This result is a special case of an abstract Berry-Esseen-type CLT for $\ell_p$-norms of sums of high-dimensional random vectors. We present this more general result together with a discussion of the related literature in Appendix~\ref{subsec:BerryEsseen}.
\end{remark}
\begin{remark}	
	The dependence of the bounds on $\sigma_{n,\max}^2$, $\sigma_{n,\min}^2$, and $\kappa_n^2$ is not necessarily optimal; e.g., if the $X_i$'s exhibit variance decay in the sense of~\cite{lopes2020bootstrapping}, directly applying the abstract Berry-Esseen-type CLT in Appendix~\ref{subsec:BerryEsseen} can yield better bounds. Moreover, we can always replace $\sigma_{n,\min}^2$ by the larger quantity $\min_{1 \leq k \leq d} \sqrt{n^{-1} \sum_{i=1}^n \mathrm{E}[X_{ik}^2]}$. 
\end{remark}

The theorem reveals that even in high dimensions the distribution of $T_{n,p}$ depends on the data mostly through the first and second moments, i.e. mean zero and covariance matrix $\Sigma_n$. This insight significantly simplifies the task of estimating the distribution of $T_{n,p}$ and is the rationale for the Gaussian parametric bootstrap procedure.

Another striking aspect of this result is the dependence on exponent $p\in [1, \infty]$. Namely, as the exponent $p$ crosses the threshold $\log d$, the upper bounds in $(i)-(iv)$ undergo a phase transition from polynomial in $r_n$ to logarithmic in $d$. This phase transition is directly related to similar behavior of the variance of $\ell_p$-norms of Gaussian random vectors~\cite[][]{paouris2018Dvoretzky}. We discuss this technical aspect in greater detail in Appendix~\ref{subsec:AntiConcentration}.

Since this Gaussian approximation result is non-asymptotic we can take limits (with respect to $n,d, p$) in any order. Given the scope of the paper, we are most interested in the high-dimensional setting with $n,d \rightarrow \infty$ and $p \in [1, \infty)$ fixed. For this asymptotic regime we note the following: The bounds in cases $(i)$, $(ii)$, and $(iv)$ imply that the larger the exponent $p$ and the stronger the moment conditions on the $X_i$'s, the faster $d$ can grow (relative to $n$) while still guaranteeing that the distributions of $T_{n,p}$ and $\widetilde{T}_{n,p}$ are close. Case $(iii)$ (with $p = \infty$) covers the case of the max-statistic considered in~\cite{chernozhukov2013GaussianApproxVec, chernozhukov2015ComparisonAnti, chernozhukov2017CLTHighDim} and improves their bound by removing the dependence on the inverse of $\sigma_{n,\min}^2$. If the $X_i$'s are identically distributed then $\kappa_n^2 = 1$ and the bound is independent of any characteristic of the covariance matrix of the data (rank, eigenvalues, or diagonal values). 

Since the $T_{n,2}$ statistic is of particular interest in many statistical applications, we provide the following easy corollary with a short discussion.

\begin{corollary}[Gaussian approximation of $T_{n,2}$]\label{corollary:GaussianApprox-L2}
	\noindent
	\begin{itemize}
		\item[(i)] If $X$ is sub-Gaussian (satisfies Assumption~\ref{assumption:SubGaussian}), then
		\begin{align}\label{eq:corollary:GaussianApprox-L2-1}
		\sup_{t \geq 0} \left|\mathrm{P}(T_{n,2} \leq t) - \mathrm{P}(\widetilde{T}_{n,2} \leq t)\right|\lesssim \sqrt{\frac{(\log d) r_n^{1/2}}{n^{1/3}}\frac{ \sigma_{n,\max}^2}{\sigma_{n,\min}^2}}.
		\end{align}
		\item[(ii)] If $X$ has finite $s \geq 4$ moments (satisfies Assumption~\ref{assumption:FiniteMoments} with $s \geq 4$), then
		\begin{align}\label{eq:corollary:GaussianApprox-L2-2}
		\sup_{t \geq 0} \left|\mathrm{P}(T_{n,2} \leq t) - \mathrm{P}(\widetilde{T}_{n,2} \leq t)\right|\lesssim \left(K_s \vee \sqrt{s}\right) \sqrt{\frac{ d^{4/(3s)} r_n^{1/2}}{n^{1/3}}\frac{ \sigma_{n,\max}^2}{\sigma_{n,\min}^2}}.
		\end{align}
	\end{itemize}
\end{corollary}
\begin{remark}
	A similar result holds for sub-Exponential random variables satisfying Assumption~\ref{assumption:SubExponential}.
\end{remark}
The result that is most related to Corollary~\ref{corollary:GaussianApprox-L2} is the dimension-dependent Berry-Essen bound by~\cite{bentkus2003DependenceBerryEsseen}.~\cite{bentkus2003DependenceBerryEsseen} addresses a slightly more general problem than we do: He derives a Berry-Esseen-type CLT for $S_n^X$ that holds uniformly over the class of Euclidean balls with arbitrary radii and arbitrary centers. In contrast,  our Corollary~\ref{corollary:GaussianApprox-L2} corresponds to a Berry-Esseen-type CLT for $S_n^X$ that holds uniformly over the class of localized Euclidean balls with arbitrary radii but center fixed to the origin. The upper bound in Theorem 1.1 in~\cite{bentkus2003DependenceBerryEsseen} is at least of order $d^{3/2} n^{-1/2}$. It appears that part of the reason why we obtain a better dependence on dimension $d$ (relative to $n$) is that we consider only localized Euclidean balls.

There is a rich literature on the closely related problem of Gaussian approximations of quadratic forms~\citep[e.g.][]{bentkus1997uniform, goetze2014ExplicitRates,  pouzo2015bootstrap, spokoiny2015bootstrap, goetze2019LargeBall, xu2019L2Asymptotics}. The Berry-Esseen-type bounds in this literature often feature a better dependence on the sample size $n$, but either have a worse dependence on dimension $d$ relative to $n$, leave the dependence on $d$ wholly unaddressed, or do not apply to degenerate distributions (i.e. low-rank covariance matrix). In general, the existing bounds appear to be less useful for applications to high-dimensional statistics than our results in this section.

\subsection{Bootstrap consistency}\label{subsec:BootstrapConsistency}
In this section we provide non-asymptotic bounds on the Kolmogorov-Smirnov distance between the distributions of the $\ell_p$-statistic $T_{n,p}$ and the Gaussian parametric bootstrap statistic $T_{n,p}^*$. As corollary we also show the consistency of the Gaussian parametric bootstrap.

Recall from Section~\ref{subsec:GPB} that the Gaussian parametric bootstrap requires a positive semi-definite estimate $\widehat{\Sigma}_n$ of the (averaged) population covariance matrix $\Sigma_n$. The non-asymptotic bounds in this section depend on the following quantities
\begin{align}\label{eq:subsec:BootstrapConsistency-1}
\widehat{\Delta}_{op} := \|\widehat{\Sigma}_n - \Sigma_n\|_{op} \hspace{30pt} \mathrm{and} \hspace{30pt} \widehat{\Delta}_p := \|\mathrm{vec}(\widehat{\Sigma}_n - \Sigma_n)\|_p, \hspace{15pt} p\in [1, \infty].
\end{align}
Note that $\widehat{\Delta}_p$ corresponds to the entry-wise $\ell_p$-norm of $\widehat{\Sigma}_n - \Sigma_n$ with $\widehat{\Delta}_2$ being the Frobenius norm. To establish the bootstrap consistency, we use
\begin{align*}
&\sup_{t \geq 0} \left|\mathrm{P}(T_{n,p} \leq t) - \mathrm{P}(T^*_{n,p} \leq t \mid X)\right|\\
&\quad{}\leq \sup_{t \geq 0} \left|\mathrm{P}(T_{n,p} \leq t) - \mathrm{P}(\widetilde{T}_{n,p} \leq t )\right| + \sup_{t \geq 0} \left|\mathrm{P}(\widetilde{T}_{n,p} \leq t ) - \mathrm{P}(T^*_{n,p} \leq t  \mid X )\right|.
\end{align*}
The first term on the right hand side in above display is deterministic and can be bounded by using  Theorem~\ref{theorem:GaussianApprox}. The second term is stochastic and can be handled by the Gaussian comparison inequality in Appendix~\ref{subsec:GaussianComparison}.

The following theorem shows that the distributions of $T_{n,p}$ and $T_{n,p}^*$ are close in Kolmogorov-Smirnov distance uniformly over all $p \in [1, \infty]$ and for generic estimates $\widehat{\Sigma}_n$.

\begin{theorem}[Consistency of the Gaussian parametric bootstrap]\label{theorem:ConsistencyGPB}
	
	\noindent	
	\begin{itemize}
	\item[(i)] For all $p \in [1, \infty)$ and $X$ satisfying Assumption~\ref{assumption:SubGaussian},
	\begin{align}\label{eq:theorem:ConsistencyGPB-1}
	\sup_{t \geq 0} \left|\mathrm{P}(T_{n,p} \leq t) - \mathrm{P}(T^*_{n,p} \leq t \mid X)\right|\lesssim \sqrt{\frac{p^3 (\log d)r_n^{1/p} }{n^{1/3}}\frac{\sigma_{n, \max}^2}{\sigma_{n, \min}^2}} + \sqrt{\frac{p^2r_n^{1/p}}{d^{1/p}}\frac{\widehat{\Delta}_p}{\sigma_{n, \min}^2}}.
	\end{align}
	\item[(ii)] For all $p \in [\log d, \infty]$ and $X$ satisfying Assumption~\ref{assumption:SubGaussian},
	\begin{align}\label{eq:theorem:ConsistencyGPB-2}
	\sup_{t \geq 0} \left|\mathrm{P}(T_{n,p} \leq t) - \mathrm{P}(T^*_{n,p} \leq t \mid X)\right|\lesssim \left(\frac{\kappa_n^2 \log^7 d}{n}\right)^{1/6} + \kappa_n(\log d) \sqrt{\frac{ \widehat{\Delta}_{op} \wedge \widehat{\Delta}_\infty }{\sigma_{n,\max}^2}}.
	\end{align}
	\item[(iii)] For $X$ satisfying Assumption~\ref{assumption:FiniteMoments} with $s \geq 4$ and all $p \in [1, s]$,
	\begin{align}\label{eq:theorem:ConsistencyGPB-3}
	\begin{split}
	&\sup_{t \geq 0} \left|\mathrm{P}(T_{n,p} \leq t) - \mathrm{P}(T^*_{n,p} \leq t \mid X)\right|\\
	&\quad{}\quad{}\lesssim (K_s \vee \sqrt{s})\sqrt{\frac{p^3 d^{4/(3s)}r_n^{1/p}}{n^{1/3}}\frac{\sigma_{n,\max}^2}{\sigma_{n,\min}^2} } + \sqrt{\frac{p^2r_n^{1/p}}{d^{1/p}}\frac{\widehat{\Delta}_p}{\sigma_{n, \min}^2}} .
	\end{split}
	\end{align}
\end{itemize}	
\end{theorem}
\begin{remark}
The second term on the right hand side of inequalities~\ref{eq:theorem:ConsistencyGPB-1}--\ref{eq:theorem:ConsistencyGPB-3} reflects the difference between $\widetilde{T}_{n,p}$ and $T_{n,p}^*$ in the Kolmogrov-Smirnov distance. Theorem~\ref{theorem:ConsistencyGPB} $(i)$ and $(ii)$ hold also for sub-Exponential random variables with the obvious modifications.
\end{remark}

Theorem~\ref{theorem:ConsistencyGPB} is only practically relevant in combination with estimates $\widehat{\Sigma}_n$ for which the stochastic estimation errors $\widehat{\Delta}_p$ and $\widehat{\Delta}_{op} \wedge \widehat{\Delta}_\infty$ are small. In Appendix~\ref{subsec:AuxResults} we provide bounds on these quantities for several different estimates $\widehat{\Sigma}_n$. For the remainder of this section we consider the special case $\widehat{\Sigma}_n = \widehat{\Sigma}_{\mathrm{naive}} := n^{-1}\sum_{i=1}^n(X_i - \bar{X}_n)(X_i - \bar{X}_n)'$. We define the \emph{naive Gaussian parametric bootstrap} estimate based on the sample covariance matrix $\widehat{\Sigma}_{\mathrm{naive}}$ by
\begin{align*}
T_{n,p, \mathrm{naive}}^*  : = \|V^\mathrm{naive}\|_p, \hspace{30pt} V^\mathrm{naive} \mid X \sim N(0,\widehat{\Sigma}_{\mathrm{naive}}).
\end{align*}
Since $T_{n,p, \mathrm{naive}}^*$ is equivalent to the Gaussian multiplier statistic $T_{n,p}^g$, the following result is also a statement about the Gaussian multiplier bootstrap.

\begin{corollary}[Consistency of the naive Gaussian parametric bootstrap]\label{corollary:ConsistencyGPB-SampleCovariance}
	Suppose that $X$ satisfies Assumption~\ref{assumption:SubGaussian}. Let $\zeta \in (0,1)$ arbitrary and set $\lambda_n \asymp \sqrt{\frac{\log d + \log (2/\zeta)}{n}} \bigvee \frac{\log d + \log (2/\zeta)}{n}$.
	\begin{itemize}
	\item[(i)] For all $p \in [1, \infty)$ with probability at least $1 - \zeta$,
	\begin{align}\label{eq:corollary:ConsistencyGPB-SampleCovariance-1}
	\begin{split}
	&\sup_{t \geq 0} \left|\mathrm{P}(T_{n,p} \leq t) - \mathrm{P}(T^*_{n,p, \mathrm{naive}} \leq t \mid X)\right|\\
	&\quad{}\lesssim \sqrt{\frac{p^3 (\log d)r_n^{1/p} }{n^{1/3}}\frac{\sigma_{n, \max}^2}{\sigma_{n, \min}^2}} +  \sqrt{p^2 \lambda_n d^{1/p} r_n^{1/p}\frac{\sigma_{n, \max}^2}{\sigma_{n, \min}^2}}.
	\end{split}
	\end{align}
	\item[(ii)] For all $p \in [\log d, \infty]$ with probability at least $1- \zeta$,
	\begin{align}\label{eq:corollary:ConsistencyGPB-SampleCovariance-2}
	\sup_{t \geq 0} \left|\mathrm{P}(T_{n,p} \leq t) - \mathrm{P}(T^*_{n,p,\mathrm{naive}} \leq t \mid X)\right| \lesssim \left(\frac{\kappa_n^2 \log^7 d}{n}\right)^{1/6} + \sqrt{\lambda_n \kappa_n^2 \log^2 d}.
	\end{align}
\end{itemize}
\end{corollary}

\begin{remark}\label{remark:corollary:ConsistencyGPB-SampleCovariance}
	The bound in case $(ii)$ depends on the covariance matrix only through the ratio $\kappa_n^2 \geq 1$. If the $X_i$'s are identically distributed then $\kappa_n^2 = 1$. For $p= \infty$ this is a useful improvement over the bounds in Theorem 4.1 and Proposition 4.1 in~\cite{chernozhukov2017CLTHighDim}.
\end{remark}

The main message of this corollary is that in high dimensions the naive Gaussian parametric and the Gaussian multiplier bootstrap can be consistent for large exponents $p \geq \log d$ but may fail to be consistent for small exponents $p \in [1, \log d)$.
More precisely, cases $(i)$ and $(ii)$ imply that the naive Gaussian parametric and the Gaussian multiplier bootstrap are consistent in probability for small $p \in [1, \log d)$ if $d^{2/p} \log d = o(n)$ and for large $p \in [\log d, \infty]$ if $\log^7 d = o(n)$. Using the Borel-Cantelli lemma, one can easily turn this into sufficient conditions for ``almost sure'' bootstrap consistency.

\subsection{Bootstrap consistency under structured covariance matrices}\label{subsec:BootstrapConsistencyStructure}
We establish two refined consistency results for the Gaussian parametric bootstrap in high dimensions. In particular, we significantly improve the rates of bootstrap consistency for small exponents $p \in [1, \log d)$ (cf. Corollary~\ref{corollary:ConsistencyGPB-SampleCovariance} $(i)$) by exploiting certain sparsity and bandedness properties of the covariance matrix. We do not present results for large exponents $p \in [\log d, \infty]$ because in this regime sparsity and bandedness properties cannot be leveraged (and are also not needed) to further improve the rates given in Corollary~\ref{corollary:ConsistencyGPB-SampleCovariance} $(ii)$.

To keep the discussion simple, we now assume that $X = \{X_i\}_{i=1}^n$ is a random sample of i.i.d. random vectors in $\mathbb{R}^d$ with mean zero and covariance matrix $\Sigma = (\sigma_{jk})_{j,k=1}^d$. We will drop the subscript $n$ in $r$, $\sigma_{\min}^2$, and $\sigma_{\max}^2$.

\begin{assumption}[Approximately sparse covariance matrix]\label{assumption:ApproxSparseCovariance}
	There exist constants $\gamma \in [0,1)$, $\theta \in [1, \infty]$ and $R_{\gamma, \theta} > 0$ such that
	\begin{align}\label{eq:assumption:ApproxSparseCovariance-1}
	\max_{1 \leq j\leq d} \left(\sum_{k=1}^d |\sigma_{jk}|^{\gamma \theta}\right)^{1/\theta} \leq R_{\gamma, \theta}.
	\end{align}
\end{assumption}
For $\gamma=0$ this assumption is most restrictive and implies that the covariance matrix is sparse with at most $R_{0,\theta}^\theta$ non-zero entries in each row. The covariance matrix of an $\mathrm{AR}$-process is a prominent example  satisfying this assumption for some positive $\gamma > 0$.

\begin{assumption}[Approximately bandable covariance matrix]\label{assumption:ApproxBandableCovariance}
	There exist constants $\alpha \in (0, \infty]$ and $\theta \in [1, \infty]$ such that for all $1 \leq \ell \leq d-1$,
	\begin{align}\label{eq:assumption:ApproxBandableCovariance-1}
	\max_{1 \leq k\leq d} \left(\sum_{j=1}^d \big\{ |\sigma_{jk}|^\theta: |j-k| > \ell  \big\}\right)^{1/\theta} \leq B_\theta \ell^{-\alpha},
	\end{align}
	for some  $B_\theta > 0$.
\end{assumption}
The larger $\alpha > 0$, the more the covariance matrix $\Sigma$ resembles a diagonal matrix. Covariance matrices of $\mathrm{MA}$-processes  satisfies this assumption for some finite $\alpha > 0$.

For $\theta=1$ Assumptions~\ref{assumption:ApproxSparseCovariance} and~\ref{assumption:ApproxBandableCovariance} reduce to two frequently adopted assumptions in the literature on high-dimensional covariance estimation~\citep[e.g.][and references therein]{bickel2008covariance, bickel2008regularization, cai2010optimal, cai2011adaptive, acella-medina2018robust}. The larger $\theta$, the milder are the restrictions imposed on the covariance matrix. 

Under Assumption~\ref{assumption:ApproxSparseCovariance} it is natural to estimate the covariance matrix via thresholding of the naive sample covariance~\citep[e.g.][]{bickel2008covariance, Lam2009sparsistency}. For simplicity, here we only consider the hard-thresholding operator; Appendix~\ref{subsec:AuxResults} contains results for more general thresholding operators. For a matrix $M = (m_{jk})_{j,k=1}^d$ and $\lambda > 0$, we define the hard-thresholding operator by
\begin{align}\label{eq:subsec:BootstrapConsistencyStructure-1}
T_\lambda(M) := \big(m_{jk}\mathbf{1}\{|m_{jk}| > \lambda\}\big)_{j,k=1}^d.
\end{align}

Under Assumption~\ref{assumption:ApproxBandableCovariance} it is common to estimate the covariance matrix via banding of the naive sample covariance~\citep[e.g.][]{bickel2008regularization}: For a given $\ell > 0$, define
\begin{align}\label{eq:subsec:BootstrapConsistencyStructure-2}
B_\ell(M) := \big(m_{jk}\mathbf{1}\{|j-k| \leq \ell\}\big)_{j,k=1}^d.
\end{align}

Recall that the Gaussian parametric bootstrap procedure requires a positive semi-definite estimate of the covariance matrix. If $\lambda_{\min}(\Sigma)$ and sample size $n$ are sufficiently large,
\cite{bickel2008covariance} and \cite{bickel2008regularization} show that $T_{\lambda}(\widehat{\Sigma}_{\mathrm{naive}})$ and $B_{\ell}(\widehat{\Sigma}_{\mathrm{naive}})$ are positive definite with probability one. If the sample size is small we suggest projecting these estimates onto the cone of positive semi-definite matrices. Since the resulting positive semi-definite projections $T_{\lambda}^+(\widehat{\Sigma}_{\mathrm{naive}})$ and $B_{\ell}^+(\widehat{\Sigma}_{\mathrm{naive}})$ maintain the same order of $\ell_p$-error as the original estimates, this projection step does not add any additional theoretical challenges. Indeed, define
\begin{align}
T_{\lambda}^+(\widehat{\Sigma}_{\mathrm{naive}}) := \arg \min_{S \succeq 0}  \| \mathrm{vec}( T_{\lambda}(\widehat{\Sigma}_{\mathrm{naive}}) - S)\|_p,
\end{align}
and observe that by triangular inequality and contraction property of projections,
\begin{align}
\begin{split}
 \| \mathrm{vec}( T_{\lambda}^+(\widehat{\Sigma}_{\mathrm{naive}}) - \Sigma)\|_p \leq 2  \| \mathrm{vec}( T_{\lambda}(\widehat{\Sigma}_{\mathrm{naive}}) - \Sigma)\|_p.
\end{split}
\end{align}

The same reasoning applies to $B_{\ell}^+(\widehat{\Sigma}_{\mathrm{naive}})$. In the following, we therefore tacitly assume that this projection step has been applied and drop the superscript ``+''.

We define the Gaussian parametric bootstrap statistics based on $ T_{\lambda}(\widehat{\Sigma}_{\mathrm{naive}})$ and $B_{\ell}(\widehat{\Sigma}_{\mathrm{naive}})$, respectively, by
\begin{align}\label{eq:subsec:BootstrapConsistencyStructure-3}
&T_{n,p, \lambda}^* := \|V^\lambda\|_p, \hspace{30pt}  V^\lambda \mid X \sim N\big(0, T_{\lambda}(\widehat{\Sigma}_{\mathrm{naive}})\big),
\end{align}
and
\begin{align}\label{eq:subsec:BootstrapConsistencyStructure-4}
&T_{n,p, \ell}^*  := \|V^\ell\|_p, \hspace{30pt}  V^\ell \mid X \sim N\big(0, B_{\ell}(\widehat{\Sigma}_{\mathrm{naive}})\big),
\end{align}
where thresholding level $\lambda > 0$ and banding parameter $\ell > 0$ will be specified below. The next two corollaries refine Corollary~\ref{corollary:ConsistencyGPB-SampleCovariance}.

\begin{corollary}[Consistency of the Gaussian parametric bootstrap under approximate sparsity]\label{corollary:ConsistencyGPB-ApproxSparsity}
	Let $X = \{X_i\}_{i=1}^n$ be a random sample of i.i.d. random vectors in $\mathbb{R}^d$ with mean zero and covariance matrix $\Sigma$. Suppose that $\Sigma$ satisfies Assumption~\ref{assumption:ApproxSparseCovariance}.
	\begin{itemize}
		\item[(i)] Set $\lambda_n \asymp \sqrt{\frac{\log d + \log (2/\zeta)}{n}} \bigvee \frac{\log d + \log (2/\zeta)}{n}$ with $\zeta \in (0,1)$ arbitrary. If in addition Assumption~\ref{assumption:SubGaussian} holds, then for all $p \in [\theta, \infty)$ with probability at least $1 - \zeta$,
		\begin{align}\label{eq:corollary:ConsistencyGPB-ApproxSparsity-1}
		\begin{split}
		&\sup_{t \geq 0} \left|\mathrm{P}(T_{n,p} \leq t) - \mathrm{P}(T^*_{n,p, \lambda_n} \leq t \mid X)\right|\\
		&\quad{}\lesssim \sqrt{\frac{p^3 (\log d)r^{1/p} }{n^{1/3}}\frac{\sigma_{\max}^2}{\sigma_{\min}^2}} +\sqrt{\frac{p^2\lambda_n r^{1/p}}{\lambda_n^\gamma}\frac{R_{\gamma, p}}{\sigma_{\max}^{2\gamma}}\frac{\sigma_{\max}^2}{\sigma_{\min}^2}}.
		\end{split}
		\end{align}
		\item[(ii)] Set $\lambda_n \asymp \sqrt{\frac{s \wedge \log d}{n}}$. If in addition Assumption~\ref{assumption:FiniteMoments} holds with $s \geq 4 \vee \theta$, then for all $p \in [2 \vee \theta, s]$,
		\begin{align}\label{eq:corollary:ConsistencyGPB-ApproxSparsity-2}
		\begin{split}
		&\sup_{t \geq 0} \left|\mathrm{P}(T_{n,p} \leq t) - \mathrm{P}(T^*_{n,p,\lambda_n} \leq t \mid X)\right|\\
		&\quad{}\lesssim (K_s \vee \sqrt{s})\sqrt{\frac{p^3 d^{4/(3s)}r^{1/p}}{n^{1/3}}\frac{\sigma_{\max}^2}{\sigma_{\min}^2} } + O_p\left(K_s^{1-\gamma} \sqrt{\frac{p^2\lambda_n d^{2/s} r^{1/p}}{(\lambda_n d^{2/s})^\gamma}\frac{R_{\gamma, p}}{\sigma_{\max}^{2\gamma}}\frac{\sigma_{\max}^2}{\sigma_{\min}^2}}\right).
		\end{split}
		\end{align}
	\end{itemize}
\end{corollary}

\begin{corollary}[Consistency of the Gaussian parametric bootstrap under approximate bandedness]\label{corollary:ConsistencyGPB-ApproxBandable}
	Let $X = \{X_i\}_{i=1}^n$ be a random sample of i.i.d. random vectors in $\mathbb{R}^d$ with mean zero and covariance matrix $\Sigma$. Suppose that $\Sigma$ satisfies Assumption~\ref{assumption:ApproxBandableCovariance}.
	\begin{itemize}
		\item[(i)] Set $\ell_n= B_p^{p/(1 + p\alpha)} \sigma_{\max}^{-2p/(1 + p\alpha)}\lambda_n^{-p/(1 + p\alpha)}$, where $\lambda_n \asymp \sqrt{\frac{\log d + \log (2/\zeta)}{n}} \bigvee \frac{\log d + \log (2/\zeta)}{n}$ and $\zeta \in (0,1)$ arbitrary. If in addition Assumption~\ref{assumption:SubGaussian} holds, then for all $p \in [\theta, \infty)$ with probability at least $1 - \zeta$,
		\begin{align}\label{eq:corollary:ConsistencyGPB-ApproxBandable-1}
		\begin{split}
		&\sup_{t \geq 0} \left|\mathrm{P}(T_{n,p} \leq t) - \mathrm{P}(T^*_{n,p, \ell_n} \leq t \mid X)\right|\\
		&\quad{}\lesssim \sqrt{\frac{p^3 (\log d)r^{1/p} }{n^{1/3}}\frac{\sigma_{\max}^2}{\sigma_{\min}^2}} +\sqrt{ \frac{p^2 \lambda_n r^{1/p}}{\lambda_n^{1/(1 + p\alpha)}} \frac{B_p^{1/(1 + p\alpha)}}{\sigma_{\max}^{2/(1 + p\alpha)}}\frac{\sigma_{\max}^2}{\sigma_{\min}^2}}.
		\end{split}
		\end{align}
		\item[(ii)] Set $\ell_n= B_p^{p/(1 + p\alpha)} \sigma_{\max}^{-2p/(1 + p\alpha)}\lambda_n^{-p/(1 + p\alpha)}$, where $\lambda_n \asymp \sqrt{\frac{s \wedge \log d}{n}}$. If in addition Assumption~\ref{assumption:FiniteMoments} holds with $s \geq 4 \vee \theta$, then for all $p \in [2 \vee \theta, s]$,
		\begin{align}\label{eq:corollary:ConsistencyGPB-ApproxBandable-2}
		\begin{split}
		&\sup_{t \geq 0} \left|\mathrm{P}(T_{n,p} \leq t) - \mathrm{P}(T^*_{n,p,\ell_n} \leq t \mid X)\right|\\
		&\quad{}\lesssim (K_s \vee \sqrt{s})\sqrt{\frac{p^3 d^{4/(3s)}r^{1/p}}{n^{1/3}}\frac{\sigma_{\max}^2}{\sigma_{\min}^2} } + O_p\left(K_s^{\frac{p\alpha}{1 + p\alpha}} \sqrt{\frac{p^2\lambda_n d^{2/s} r^{1/p}}{(\lambda_n d^{2/s})^{\frac{1}{1+p\alpha}}}\frac{B_p^{\frac{1}{1+p\alpha}}}{\sigma_{\max}^{\frac{2}{1+p\alpha}}}\frac{\sigma_{\max}^2}{\sigma_{\min}^2}}\right).
		\end{split}
		\end{align}
	\end{itemize}
\end{corollary}
\begin{remark}
	For large exponents $p \in [\log d, \infty]$ the bootstrap statistics $T_{n,p, \lambda_n}^*$ and $T_{n,p, \ell_n}^*$ satisfy the upper bounds in Corollary~\ref{corollary:ConsistencyGPB-SampleCovariance} $(ii)$.
\end{remark}
The main takeaway from these two corollaries is that under reasonable assumptions on the covariance structure and the tails of the data there exist Gaussian parametric bootstrap statistics $T_{n,p}^*$ that are consistent in high dimensions for any fixed $p \in [1, \infty)$.

In particular, inequality~\eqref{eq:corollary:ConsistencyGPB-ApproxSparsity-1} (inequality~\eqref{eq:corollary:ConsistencyGPB-ApproxBandable-1}) implies that if the data is sub-Gaussian and the population covariance matrix is approximately sparse (approximately bandable) the Gaussian parametric bootstrap based on the thresholded covariance matrix (the banded covariance matrix) is consistent in probability for all $p \in [1, \infty)$ provided that $\log d = o(n^3)$.

\section{Application: Testing high-dimensional mean vectors}\label{sec:SimultaneousTesting}
As an application of the Gaussian parametric bootstrap we now present a bootstrap hypothesis test based on $\ell_p$-statistics for testing linear restrictions on high-dimensional mean vectors. We show that this test is asymptotic correct and consistent. Moreover, we discuss the effect of the exponent $p$ on the size of simultaneous confidence sets and the power of the test. Lastly, we discuss an extension of the generic testing framework to simultaneous inference on high-dimensional linear models.

\subsection{Setup and test statistic}\label{subsec:SetupTesting}
Given a random sample $X = \{X_i\}_{i=1}^n$ of i.i.d. random vectors in $\mathbb{R}^d$ with unknown mean $\mu$ and unknown covariance matrix $\Sigma$ we are interested in testing the high-dimensional linear restrictions
\begin{align}\label{eq:sec:SimultaneousTesting-1}
H_0: \:\: M\mu = m_0 \hspace{20pt} \mathrm{vs.} \hspace{20pt} H_1: \:\: M \mu \neq m_0,
\end{align}
for some $M \in \mathbb{R}^{d' \times d}$ and $m_0 \in \mathbb{R}^{d'}$ when dimension $d$ and number of restrictions $d'$ may exceed the sample size $n$.

We propose to test hypothesis~\eqref{eq:sec:SimultaneousTesting-1} on the basis of the $\ell_p$-statistic
\begin{align}\label{eq:sec:SimultaneousTesting-2}
S_{n,p} := \left\|\frac{1}{\sqrt{n}} \sum_{i=1}^n (MX_i - m_0)\right\|_p, \hspace{30pt} p \geq 1,
\end{align}
and, given a nominal level $\alpha \in (0,1)$, reject the null hypothesis if and only if
\begin{align}\label{eq:sec:SimultaneousTesting-3}
S_{n,p} \geq c^*_{n,p}(1- \alpha),
\end{align}
where $c^*_{n,p}(\alpha)$
is the $\alpha$-quantile of the Gaussian parametric bootstrap estimate
\begin{align}\label{eq:sec:SimultaneousTesting-4}
S_{n,p}^* := \|V^0\|_p, \hspace{30pt} V^0 \mid X \sim N(0, \widehat{\Omega}_n),
\end{align}
and $\widehat{\Omega}_n$ is a positive semi-definite estimate of $\Omega = M \Sigma M'$.

A distinguishing feature of this bootstrap hypothesis test is the exponent $p \in [1, \infty]$ and we show that the exponent $p$ has significant impact on the asymptotic correctness and the power of the test. In practice, tests based on $\ell_p$-statistics $S_{n,p}$ with exponents $1$, $2$, and $\infty$ are of particular interest. For one, the $\ell_1$-statistic $S_{n,1}$ and the maximum statistic $S_{n,\infty}$ lie at opposite ends of the spectrum of possible exponents $p$ and therefore have power functions that are complementary in a sense to be made precise below. For another, the maximum statistic $S_{n,\infty}$ can also be applied to the problem of multiple hypothesis testing. Since the bootstrap test based on $S_{n,\infty}$ accounts for the dependence between the multiple tests, it is (asymptotically) less conservative than the Bonferroni adjustment.
Lastly, the sum-of-squares type statistic $S_{n,2}$ is essentially a feasible version of Hotelling's $T^2$-statistic in high dimensions and as such interesting in its own right~\citep{fan2015power}.

Let $\mathcal{H}_0 = \{\mu \in \mathbb{R}^d : M\mu = m_0\}$  and $\mathcal{H}_1 = \mathcal{H}_0^c$.  Write $\Omega =(\omega_{jk})_{j,k=1}^{d'}$, and
\begin{align}\label{eq:sec:SimultaneousTesting-5}
r_\omega : = \mathrm{rank}(\Omega), \hspace{20pt} \omega^2 : = (\omega_{kk})_{k=1}^{d'}, \hspace{20pt} \omega_{\min}^2 := \min_{1 \leq k \leq d'} \omega_k^2, \hspace{20pt} \omega_{\max}^2 := \max_{1 \leq k \leq d}\omega_k^2.
\end{align}
Let $\widehat{\Omega}_n$ be a positive semi-definite estimate of $\Omega$ and define
\begin{align}\label{eq:sec:SimultaneousTesting-6}
\widehat{\Gamma}_{op} := \|\widehat{\Omega}_n - \Omega\|_{op} \hspace{30pt} \mathrm{and} \hspace{30pt} \widehat{\Gamma}_p := \|\mathrm{vec}(\widehat{\Omega}_n - \Omega)\|_p, \hspace{15pt} p\in [1, \infty].
\end{align}

We also introduce the following high-level assumption.

\begin{assumption}[Asymptotic sufficient conditions]\label{assumption:SufficientConditions} At least one of the following statements holds true.
	\begin{itemize}
		\item[(i)] Assumption~\ref{assumption:SubGaussian} holds, $p \in [1, \infty)$,
		\begin{align*}
		(\log^3 d')r_\omega^{3/p} \omega_{\max}^6\omega_{\min}^{-6} = o(n), \hspace{20pt} and \hspace{20pt} \widehat{\Gamma}_p = o_p\left(r_\omega^{-1/p}{d'}^{1/p} \omega_{\min}^2\right).
		\end{align*}
		\item[(ii)] Assumption~\ref{assumption:SubGaussian} holds, $p \in [\log d', \infty]$,
		\begin{align*}
		\log^7 d' = o(n), \hspace{20pt} and \hspace{20pt} \widehat{\Gamma}_{op}\wedge\widehat{\Gamma}_\infty = o_p\left( (\log d')^{-2}\omega_{\max}^2\right).
		\end{align*}
		\item[(iii)] Assumption~\ref{assumption:FiniteMoments} holds with $s \geq 4$, $p \in [1, s]$,
		\begin{align*}
		(K_s^2 \vee s)^3 {d'}^{4/s}r_\omega^{3/p} \omega_{\max}^6\omega_{\min}^{-6}= o( n), \hspace{20pt} and \hspace{20pt} \widehat{\Gamma}_p = o_p\left(r_\omega^{-1/p}{d'}^{1/p} \omega_{\min}^2\right).
		\end{align*}
	\end{itemize}
\end{assumption}
We emphasize that under rather mild conditions there exist estimates $\widehat{\Omega}_n$ such that $\widehat{\Gamma}_n$ and $\widehat{\Gamma}_{op} \wedge\widehat{\Gamma}_\infty$ satisfy the conditions in Assumption~\ref{assumption:SufficientConditions}; see Appendix~\ref{subsec:AuxResults} for details.

\subsection{Asymptotic correctness}\label{subsec:AsympCorrectness}
In this section we show that the bootstrap hypothesis test has asymptotic correct size. We state the theorem in a non-asymptotic fashion to match the results from previous sections.

\begin{theorem}[Asymptotic size $\alpha$ test]\label{theorem:SizeAlphaTest}
	Let $\xi$ be an arbitrary real-valued random variable, whose role will be discussed afterwards.
	\begin{itemize}
		\item[(i)] For all $p \in [1, \infty)$ and $X$ satisfying Assumption~\ref{assumption:SubGaussian},
		\begin{align}\label{eq:theorem:SizeAlphaTest-1}
		\begin{split}
		&\sup_{\alpha \in (0,1)} \sup_{\mu \in \mathcal{H}_0} \left|\mathrm{P}_\mu\big(S_{n,p} + \xi \leq c^*_{n,p}(\alpha)\big) - \alpha \right|\\
		&\quad{}\quad{}\lesssim \sqrt{\frac{p^3 (\log d')r_\omega^{1/p} }{n^{1/3}}\frac{\omega_{\max}^2}{\omega_{\min}^2}} + \inf_{\delta > 0}\left\{ \sqrt{\frac{p^2r_\omega^{1/p}}{{d'}^{1/p}}\frac{\delta}{\omega_{\min}^2}} + \mathrm{P}\left(\widehat{\Gamma}_p > \delta\right) \right\} \\
		&\quad{}\quad{} \quad{}+ \inf_{\eta > 0} \left\{ \sqrt{\frac{p r_\omega^{1/p}}{{d'}^{2/p}}\frac{\eta^2}{\omega_{\min}^2}} + \mathrm{P} \left(|\xi| > \eta \right)\right\}.
		\end{split}
		\end{align}
		\item[(ii)] For all $p \in [\log d, \infty]$ and $X$ satisfying Assumption~\ref{assumption:SubGaussian},
		\begin{align}\label{eq:theorem:SizeAlphaTest-2}
		\begin{split}
		&\sup_{\alpha \in (0,1)} \sup_{\mu \in \mathcal{H}_0} \left|\mathrm{P}_\mu\big(S_{n,p} + \xi \leq c^*_{n,p}(\alpha)\big) - \alpha \right|\\
		&\quad{}\quad{}\lesssim \left(\frac{\log^7 d'}{n}\right)^{1/6} +  \inf_{\delta > 0}\left\{(\log d') \sqrt{\frac{ \delta }{\omega_{\max}^2}} + \mathrm{P}\left(\widehat{\Gamma}_{op} \wedge\widehat{\Gamma}_\infty  > \delta\right) \right\}\\
		&\quad{}\quad{} \quad{} +\inf_{\eta > 0} \left\{ (\log d') \sqrt{\frac{ \eta^2 }{\omega_{\max}^2}} +\mathrm{P} \left(|\xi| > \eta \right)\right\}.
		\end{split}
		\end{align}
		\item[(iii)] For $X$ satisfying Assumption~\ref{assumption:FiniteMoments} with $s \geq 4$ and all $p \in [1, s]$,
		\begin{align}\label{eq:theorem:SizeAlphaTest-3}
		\begin{split}
		&\sup_{\alpha \in (0,1)} \sup_{\mu \in \mathcal{H}_0} \left|\mathrm{P}_\mu\big(S_{n,p} + \xi \leq c^*_{n,p}(\alpha)\big) - \alpha \right|\\
		&\quad{}\quad{}\lesssim (K_s \vee \sqrt{s})\sqrt{\frac{p^3 {d'}^{4/(3s)}r_\omega^{1/p}}{n^{1/3}}\frac{\omega_{\max}^2}{\omega_{\min}^2} } +  \inf_{\delta > 0}\left\{ \sqrt{\frac{p^2r_\omega^{1/p}}{{d'}^{1/p}}\frac{\delta}{\omega_{\min}^2}} + \mathrm{P}\left(\widehat{\Gamma}_p > \delta\right)\right\}\\
		&\quad{}\quad{} \quad{}+ \inf_{\eta > 0} \left\{ \sqrt{\frac{p r_\omega^{1/p}}{{d'}^{2/p}}\frac{\eta^2}{\omega_{\min}^2}} + \mathrm{P} \left(|\xi| > \eta \right)\right\}.
		\end{split}
		\end{align}
	\end{itemize}	
\end{theorem}
\begin{remark}
	For $\widehat{\Omega}_n = \widehat{\Omega}_{\mathrm{naive} }:= n^{-1}\sum_{i=1}^nM(X_i - \bar{X}_n)(X_i - \bar{X}_n)'M'$ these bounds also hold for quantiles $c^{*g}_{n,p}(\alpha)$ obtained via the Gaussian multiplier bootstrap procedure.
\end{remark}

A special feature of this result is the real-valued random variable $\xi$. For now, assume that $\xi \equiv 0$ and let $\eta\downarrow 0 $ arbitrarily fast. In this case, Theorem~\ref{theorem:SizeAlphaTest} provides non-asymptotic error bounds on the type I error of the bootstrap hypothesis test based on $\ell_p$-statistic $S_{n,p}$.

Next, consider the case in which $\xi$ is not identical to zero. Then, Theorem~\ref{theorem:SizeAlphaTest} is a statement about the test statistic $R_{n,p} : = S_{n,p} + \xi$, where $\xi$ may be interpreted as approximation error. This is particularly useful if we want to test hypotheses about a parameter $\beta_0 \in \mathbb{R}^d$ for which there exists an estimator $\hat{\beta}$ that admits the expansion
\begin{align}
\sqrt{n}M(\hat{\beta} - \beta_0) = \frac{1}{\sqrt{n}} \sum_{i=1}^n (MX_i - m_0) + r_n.
\end{align}
In this case, the triangle inequality yields $|\xi| \leq \|r_n\|_p$. The primary example that we have in mind is the de-biased lasso estimator for linear models~\citep[e.g.][]{vandegeer2014on, zhang2014confidence}. We elaborate on this idea in detail in Section~\ref{subsec:DebiasedLasso}.


\subsection{Confidence sets for high-dimensional parameters}\label{subsec:ConfidenceSets}
We can use Theorem~\ref{theorem:SizeAlphaTest} to construct consistent confidence sets $\mathcal{C}_{n,p} \subset \mathbb{R}^d$ for a high-dimensional parameter $\mu_0 \in \mathbb{R}^d$. To this end, set $M = I_d$, $m_0 = \mu_0$, and define
\begin{align}
\mathcal{C}_{n,p} : = \left\{ \mu \in \mathbb{R}^d : \left\|\frac{1}{n}\sum_{i=1}^n X_i - \mu \right\|_p \leq \frac{c_{n,p}^*(1-\alpha)}{\sqrt{n}}\right\},
\end{align}
for a given nominal level $\alpha \in (0,1)$.
Then, under Assumption~\ref{assumption:SufficientConditions}, Theorem~\ref{theorem:SizeAlphaTest} guarantees that
\begin{align}\label{eq:sec:SimultaneousTesting-7-1}
\mathrm{P}_{\mu_0}\left(\mu_0 \in \mathcal{C}_{n,p}\right) \rightarrow 1 - \alpha.
\end{align}

Given the collection $\{\mathcal{C}_{n,p}\}_{p \geq 1}$ a practitioner will be most interested in knowing which of these confidence sets is ``smallest''. To answer this question, we study how $p \in [1, \infty]$ affects the volume of $\mathcal{C}_{n,p}$ as $d, n \rightarrow \infty$. 
To simplify matters, we only consider $\Sigma = \sigma^2 I_d$.

Obviously, the confidence sets $\mathcal{C}_{n,p}$ are just $\ell_p$-norm balls with center $\bar{X}_n$ and radii $c^*_{n,p}(1- \alpha)/\sqrt{n}$. Recall that the volume of centered $d$-dimensional $\ell_p$-balls with radius $r > 0$, say $\mathcal{B}_p^d(r)$, is given by
\begin{align}\label{eq:sec:SimultaneousTesting-7}
	\mathrm{Vol}\left(\mathcal{B}_p^d(r)\right) = \frac{(2r)^d \Gamma\left(1 + 1/p\right)^d}{
	\Gamma\left(1 + d/p\right)}.
\end{align}

Also, by Lemma~\ref{lemma:UpperBoundQuantiles}, Remark~\ref{remark:lemma:UpperBoundQuantiles}, and Lemma 2 in~\cite{schechtman1990volume}, with probability approaching one, for all $\alpha \in (0, 1/2)$,
\begin{align}\label{eq:sec:SimultaneousTesting-8}
c^*_{n,p}(1-\alpha)/\sqrt{n} \asymp \begin{cases}
\sigma d^{1/p} \sqrt{p/n},  & p < \log d\\
\sigma \sqrt{(\log d)/n}, & p \geq \log d.
\end{cases}
\end{align}

Whence, by~\eqref{eq:sec:SimultaneousTesting-7},~\eqref{eq:sec:SimultaneousTesting-8}, and Sterling's formula we have 
\begin{align}\label{eq:sec:SimultaneousTesting-9}
\mathrm{Vol}\left(\mathcal{C}_{n,p}\right) \asymp \begin{cases}
\displaystyle \left(\frac{ep}{c_p}\right)^{d/p} \left(\frac{p}{d}\right)^{1/2} \left(\frac{4 \sigma^2 p}{n}\right)^{d/2},  & p < \log d\\[15pt]
\displaystyle \left(\frac{ep}{c_p d} \wedge \frac{e}{c_p}\right)^{d/p}  \left(\frac{p}{d} \wedge 1\right)^{1/2} \left(\frac{4\sigma^2 \log d}{n}\right)^{d/2}, & p \geq \log d,
\end{cases}
\end{align}
where $c_{p}^{1/p} \in (0.8856, 1]$. It is now easy to check that (asymptotically) the volume of $\mathcal{C}_{n,p}$ is a monotonically increasing function of the exponent $p$. In other words, confidence sets based on $\ell_p$-statistics $S_{n,p}$ with small exponents are less conservative than confidence sets based on, say, the maximum statistic $S_{n,\infty}$. Asymptotically, $\mathcal{C}_{n,1}$ is the smallest confidence set.


\subsection{Consistency under high-dimensional alternatives}\label{subsec:ConsistencyTesting}
We now analyze the consistency of the bootstrap hypothesis test under high-dimensional alternatives. Let $Z \sim N(0, I_{d'})$ and define
\begin{align}
\mathcal{A}_p := \left\{ \left(\mu_n\right)_{n \in \mathbb{N}}, \mu_n \in \mathbb{R}^{d_n} :  \frac{ \mathrm{E}\|\Omega^{1/2}Z\|_p \vee \sqrt{\mathrm{Var}\|\Omega^{1/2}Z\|_p}}{ \sqrt{n}\|M\mu_n - m_0\|_p} = o(1)\right\},
\end{align}
and its ``complement''
\begin{align}
\mathcal{Z}_p := \left\{\left(\mu_n\right)_{n \in \mathbb{N}}, \mu_n \in \mathbb{R}^{d_n} : \frac{ \sqrt{n}\|M\mu_n - m_0\|_p}{\mathrm{E}\|\Omega^{1/2}Z\|_p \vee \sqrt{\mathrm{Var}\|\Omega^{1/2}Z\|_p}} = o(1)\right\}.
\end{align}

In words, $\mathcal{A}_p$ contains alternatives $(\mu_n)_{n \in \mathbb{N}}$ whose signals $\sqrt{n}\|M\mu_n - m_0\|_p$ asymptotically dominate the mean and standard deviation of the Gaussian proxy statistic $\|\Omega^{1/2}Z\|_p$; whereas $\mathcal{Z}_p$ consists of alternatives whose signals are asymptotically negligible compared to mean and standard deviation of $\|\Omega^{1/2}Z\|_p$.

The following result shows that the bootstrap hypothesis test is consistent for all $(\mu_n)_{n \in \mathbb{N}} \in \mathcal{A}_p$ and inconsistent for all $(\mu_n)_{n \in \mathbb{N}} \in \mathcal{Z}_p$.

\begin{theorem}[Consistency under high-dimensional alternatives]\label{theorem:ConsistencyLocalAlternatives}
	Suppose that Assumption~\ref{assumption:SufficientConditions} holds and $\sqrt{\mathrm{Var}\|\Omega^{1/2}Z\|_p} = o\left(\mathrm{E}\|\Omega^{1/2}Z\|_p\right)$.
	\begin{itemize}
		\item[(i)] 	For $\alpha \in (0,1)$ and all $\left(\mu_n\right)_{n \in \mathbb{N}} \in \mathcal{A}_p$,
		\begin{align*}
		\lim_{n \rightarrow \infty}\mathrm{P}_{\mu_n}\Big(S_{n,p} > c_{n,p}^*(1- \alpha)\Big) = 1.
		\end{align*}
		\item[(ii)]	For $\alpha \in (0,1/2)$ and all $\left(\mu_n\right)_{n \in \mathbb{N}} \in \mathcal{Z}_p$,
		\begin{align*}
		\lim_{n \rightarrow \infty} \mathrm{P}_{\mu_n}\Big(S_{n,p} > c_{n,p}^*(1- \alpha)\Big) < 1.
		\end{align*}
	\end{itemize}
\end{theorem}

\begin{remark}\label{remark:theorem:ConsistencyLocalAlternatives}
	Under mild moment conditions the ``relative standard deviation'' $\sqrt{\mathrm{Var}\|\Omega^{1/2}Z\|_p}/ \\ \mathrm{E}\|\Omega^{1/2}Z\|_p$ tends to zero as the dimension $d \rightarrow \infty$ grows~\citep{boucheron2013concentration, biau2015HighDimPNorms}. In particular, by the Gaussian Poincar{\'e} inequality, $\sqrt{\mathrm{Var}\|\Omega^{1/2}Z\|_p} \leq \|\Omega^{1/2}\|_{2 \rightarrow p} \leq \|\Omega^{1/2}\|_{op}$, where the first inequality holds for all $p \in [1, \infty]$ and the second for at least all $p \geq 2$.
\end{remark}


\subsection{Power and the role of the exponent $p$}\label{subsec:Power}
It is part of statistical folklore that sum-of-squares type statistics have good power against ``dense'' alternatives, i.e alternatives whose signals in $M\mu$ are spread out over a large number of coordinates, whereas maximum type statistics are more powerful against ``sparse'' alternatives, i.e. alternatives with only a few strong signals in $M\mu$ \cite{fan2015power}. Theorem~\ref{theorem:ConsistencyLocalAlternatives} allows us to verify this statement more formally. Let $M = I_d$, $m_0 = 0$, $\Sigma = \sigma^2I_d$, and define the set of alternatives
\begin{align}\label{eq:sec:SimultaneousTesting-12}
\mathcal{D}_{\delta,s} := \left\{ \mu \in \mathbb{R}^s \times \{0\}^{d-s}: \delta/c \leq \mu_k/\sigma  \leq \delta c,\: 1 \leq k \leq s \right\},
\end{align}
where $c \geq 1 $ is an absolute constant, $\delta > 0$ regulates the signal strength, and $s \in \{1, \ldots, d\}$ controls the sparsity.

Given this setup, we ask the following question: What is the minimum signal strength $\delta \equiv \delta(n, d, s, p)$ needed for the bootstrap test based on $S_{n,p}$ to reject the null hypothesis $H_0: \mu = 0$ at significance level $\alpha \in (0,1/2)$ when $\mu \in \mathcal{D}_{\delta, s}$?

By Remark~\ref{remark:theorem:ConsistencyLocalAlternatives} $\sqrt{\mathrm{Var}\|\Omega^{1/2}Z\|_p} \leq \sigma$ and by Lemma 2 in~\cite{schechtman1990volume} $\mathrm{E}\|\Omega^{1/2}Z\|_p \asymp \sqrt{p} d^{1/p}$ for $p < \log d$ and $\mathrm{E}\|\Omega^{1/2}Z\|_p \asymp \sigma \sqrt{\log d}$ for $p \geq \log d$. Thus, by Theorem~\ref{theorem:ConsistencyLocalAlternatives} $(ii)$, a necessary condition for correctly rejecting the null hypothesis (with probability approaching one) is
\begin{align}\label{eq:sec:SimultaneousTesting-13}
\sqrt{n}\|\mu\|_p \gtrsim  \begin{cases}
\sigma \sqrt{p} d^{1/p},  & p < \log d\\
\sigma \sqrt{\log d}, & p \geq \log d.
\end{cases}
\end{align}

Now, suppose that $s \asymp d$, i.e. $\mathcal{D}_{\delta, s}$ contains only dense alternatives. Then, for $p \in [1, \log d)$ \eqref{eq:sec:SimultaneousTesting-13} holds if $\delta \gtrsim \sqrt{p/n}$, whereas for $p \in [\log d, \infty]$ \eqref{eq:sec:SimultaneousTesting-13} holds only if $\delta \gtrsim \sqrt{(\log d)/n}$. Thus, bootstrap tests based on $\ell_p$-statistics with small exponents are more powerful in detecting dense alternatives than those based on $\ell_p$-statistics with large exponents.

Next, assume that $s \ll d$, i.e. $\mathcal{D}_{\delta, s}$ contains only sparse alternatives. Then, for $p \in [1, \log d)$ \eqref{eq:sec:SimultaneousTesting-13} holds if $\delta \gtrsim \sqrt{p(d/s)^{2/p}/n}$, whereas for $p \in [\log d, \infty]$ \eqref{eq:sec:SimultaneousTesting-13} holds already if $\delta \gtrsim \sqrt{(\log d)/n}$. Therefore, tests based on $\ell_p$-statistics with large exponents are more responsive to sparse alternatives than those based on $\ell_p$-statistics with small exponents.

\subsection{Simultaneous inference on high-dimensional linear models}\label{subsec:DebiasedLasso}
The bootstrap hypothesis test based on the $\ell_p$-statistic $S_{n,p}$ can be combined with the de-biased Lasso estimator~\citep[][]{vandegeer2014on, zhang2014confidence} to conduct simultaneous inference on high-dimensional linear models. This approach extends the one by~\cite{zhang2017SimultaneousInf}, who propose a bootstrap test for the de-biased lasso estimator based on the Gaussian multiplier bootstrap for the maximum statistic $S_{n, \infty}$.

Consider the high-dimensional sparse model
\begin{align}
Y_i = X_i'\beta_0 + \varepsilon_i, \hspace{20pt} i=1, \ldots, n,
\end{align}
with response $Y_i \in \mathbb{R}$, i.i.d. predictors $X_i \in \mathbb{R}^d$ with mean $\mu$ and covariance matrix $\Sigma$, i.i.d. errors $\varepsilon_i$ (independent of $X_i$) with mean 0 and variance $\sigma_\varepsilon^2$, and sparse regression vector $\beta_0$. We are interested in testing the linear hypothesis
\begin{align}
H_0: \:\: M\beta_0 = m_0 \hspace{20pt} \mathrm{vs.} \hspace{20pt}  H_1:\:\: M\beta_0 \neq m_0.
\end{align}
Write $Y = (Y_1, \ldots,Y_n) \in \mathbb{R}^n$, $\varepsilon = (\varepsilon_1, \ldots, \varepsilon_n) \in \mathbb{R}^n$, and $\mathbf{X} = [X_1, \ldots, X_n]' \in \mathbb{R}^{n \times d}$. For $\lambda > 0$ define the ordinary lasso estimate by
\begin{align}
\hat{\beta}_\lambda : = \arg\min_{\beta \in \mathbb{R}^d} \|Y - \mathbf{X} \beta \|_2^2/n + 2\lambda \|\beta\|_1,
\end{align}
and the de-biased lasso estimate by
\begin{align}
\breve{\beta} : = \hat{\beta}_\lambda + \widehat{\Theta}\mathbf{X}'(Y - \mathbf{X}\hat{\beta}_\lambda)/n,
\end{align}
where $\widehat{\Theta}$ is a suitable approximation of the inverse of the Gram matrix $\widehat{\Sigma} = \mathbf{X}'\mathbf{X}/n$. Define the $\ell_p$-statistic
\begin{align}
R_{n,p} : = \sqrt{n}\|M\breve{\beta} - m_0\|_p,
\end{align}
and observe that
\begin{align}\label{eq:subsec:DebiasedLasso-1}
\begin{split}
\sqrt{n}M(\breve{\beta} - \beta_0) &= \frac{1}{\sqrt{n}}\sum_{i=1}^n M\Sigma^{-1} X_i \varepsilon_i \\
&\quad{}+ \underbrace{M(\widehat{\Theta} - \Sigma^{-1})\mathbf{X}' \varepsilon/\sqrt{n}}_{=:r_1}  - \underbrace{\sqrt{n}M(\widehat{\Theta} \widehat{\Sigma} - I_d)(\hat{\beta}_\lambda - \beta_0)}_{=:r_2}.
\end{split}
\end{align}

Further, note that the first term on the right hand side in above display is the re-scaled sum of $n$ i.i.d. random vectors with mean zero and covariance matrix $\sigma_\varepsilon^2 M\Sigma^{-1}M'$. Hence, $R_{n,p} = S_{n,p} + \xi$, where $S_{n,p} = \|n^{-1/2}\sum_{i=}^n M\Sigma^{-1} X_i \varepsilon_i\|_p$ and $|\xi| \leq \|r_1\|_p + \|r_2\|_p$. Under mild assumptions, $\|r_1\|_p \leq \|M(\widehat{\Theta} - \Sigma^{-1})\|_{q \rightarrow p} \|\mathbf{X}'\varepsilon\|_q/\sqrt{n}$ and $\|r_2\|_p \leq \|M(\widehat{\Theta} \widehat{\Sigma} - I_d)\|_{q \rightarrow p}\|\hat{\beta}_\lambda - \beta_0\|_q$, $q \geq 1$, are negligible~\cite[][]{vandegeer2014on}. Thus, based on the expansion~\eqref{eq:subsec:DebiasedLasso-1} and the discussion in Section~\ref{subsec:AsympCorrectness} we can approximate the distribution of $R_{n,p}$ under the null hypothesis by the distribution of the Gaussian parametric bootstrap estimate
\begin{align}
S^*_{n,p} : = \left\|V^{\mathrm{debias}}\right\|_p, \hspace{20pt}  \mathrm{where} \hspace{20pt} V^{\mathrm{debias}}  \mid \{Y, \mathbf{X}\} \sim N(0, \hat{\sigma}_\varepsilon^2 M\widehat{\Theta}M'),
\end{align}
and $\hat{\sigma}_\varepsilon^2$ is a consistent estimate of the error variance $\sigma_{\varepsilon}^2$~\citep[][]{fan2012variance}.

We can now use the quantiles of $S^*_{n,p}$ to compute (bootstrap) critical values for the $\ell_p$-statistic $S_{n,p}$ and to construct confidence sets for $M\beta_0$.

\section{Numerical experiments}\label{sec:NumericalExperiments}
The purpose of the numerical experiments is in this section is threefold. First, they show that for small exponents $p \in [1, \log d)$ the Gaussian parametric bootstrap outperforms the Gaussian multiplier bootstrap, while for large exponents $p \in [\log d, \infty)$ both bootstrap procedures perform similarly. Second, they confirm the theoretical claims from Section~\ref{sec:MainResults} that for heavy-tailed data the accuracy of the Gaussian parametric and multiplier bootstrap suffers as the exponent $p$ increases. Third, they show that the exponent $p$ affects the power of the bootstrap hypothesis test as described in Section~\ref{sec:SimultaneousTesting}.

\subsection{Data generation}\label{subsec:Setup}
We generate vectors $X_1, \ldots, X_n \in \mathbb{R}^d$ via a Gaussian copula model
\begin{align}
X_{ij} = F^{-1}\left(\Phi(Y_{ij}) \right), \hspace{20pt} 1 \leq i \leq n, \hspace{20pt} 1 \leq j \leq d,
\end{align}
where the random vectors $Y_1, \ldots, Y_n \in \mathbb{R}^d$ are sampled independently and identically from a centered Gaussian distribution with sparse covariance matrix $\Sigma$, $\Phi$ is the cdf of the $N(0,1)$ distribution, and $F$ is the distribution function of either the uniform distribution on $[-1,1]$ (``light-tailed'') or Student's $t$-distribution with $4$ degrees of freedom (``heavy-tailed''). We create the sparse and low-rank covariance matrix $\Sigma$ in two steps: First, define the block diagonal matrix $\widetilde{\Sigma} = \mathrm{diag}(\Lambda, \ldots, \Lambda) \in \mathbb{R}^{d \times d}$, where $\Lambda = (\Lambda_{jk})_{j,k=1}^{d/100}$ with $\Lambda_{jk} = 0.8^{j+k -2}$ for all $1 \leq j,k \leq d/100$ is a rank-one matrix. Then, (randomly) generate a permutation matrix $P$ and set $\Sigma = P\widetilde{\Sigma} P'$. The matrix $\Sigma$ is positive semi-definite, sparse with $d/100$ non-zero elements in each row, and has rank $100$. The permutation matrix $P$ is generated only once and is the same throughout all Monte Carlo simulations.

\subsection{Specific implementation of a hard-thresholding estimator}
The Gaussian parametric bootstrap procedure requires as input a positive semi-definite estimate of the population covariance matrix $\Sigma$. To exploit the sparsity of $\Sigma$ while also ensuring positive semi-definiteness of the estimate, we propose the following two-step procedure:

First, compute a pilot estimate via correlation thresholding~\citep{fan2011high} of the sample covariance matrix $\widehat{\Sigma}_{\mathrm{naive}} = \left(\hat{\sigma}_{jk}\right)_{j,k =1}^d$, i.e. for $\lambda > 0$ compute
\begin{align}\label{eq:subsec:DataGeneration-1}
\widehat{\Sigma}_n(\lambda) := T_\lambda^{\mathrm{cor}}(\widehat{\Sigma}_{\mathrm{naive}}) = \left(\hat{\sigma}_{jk}\mathbf{1}\left\{\frac{|\hat{\sigma}_{jk}|}
{\sqrt{\hat{\sigma}_{jj}\hat{\sigma}_{kk}}} \geq  \lambda \right\}\right)_{j,k=1}^d.
\end{align}
Then, project the pilot estimate $\widehat{\Sigma}_n(\lambda)$ onto the cone of positive semi-definite matrices by setting all negative eigenvalues equal to 0. Denote the resulting estimate by $\widehat{\Sigma}_n^+(\lambda)$.

It remains to choose the thresholding level $\lambda > 0$. We proceed as in~\cite{bickel2008covariance, bickel2008regularization} and select $\lambda$ by cross-validation: At each fold $\nu \in \{1, \ldots, N\}$, randomly split the sample $X= \{X_i\}_{i=1}^n$ into two sub-samples $X^1$ and $X^2$ of sizes $n_1 = \lceil n/3 \rceil$ and $n_2 = n - n_1$, respectively. Denote by $\widehat{\Sigma}_{1, \nu}$ and $\widehat{\Sigma}_{2,\nu}$ the sample covariance matrices of the $\nu$th split based on $X^1$ and $X^2$. Let $\widehat{\Sigma}_{1,\nu}^{+}(\lambda)$ be the correlation-thresholded and projected estimate based on $\widehat{\Sigma}_{1, \nu}$.  Define the cross-validated risk at level $\lambda > 0$ by
\begin{align}
\widehat{R}(\lambda) : = \frac{1}{N}\sum_{\nu=1}^{N}\left\|\mathrm{vec}\left(\widehat{\Sigma}_{1, \nu}^+(\lambda)- \widehat{\Sigma}_{2, \nu} \right)\right\|_2,
\end{align}
and select the ``optimal'' thresholding level as
\begin{align}
\hat{\lambda} : = \arg\min_{\lambda\in [0,1]} \widehat{R}(\lambda).
\end{align}

In practice, we set $N = 10$ and minimize the risk $\widehat{R}(\lambda)$ over a grid $\mathcal{G} \subset [0, 1]$ with $|\mathcal{G}| =40$ equally spaced points. The algorithm is sensitive to the number of folds and grid points; increasing $N$ and $|\mathcal{G}|$ beyond 10 and 40, respectively, can improve the accuracy of the bootstrap approximation (at the cost of additional computational complexity).

\subsection{Performance of Gaussian parametric and multiplier bootstrap}\label{subsec:PerformanceBootstrap}
To assess the performance of the Gaussian parametric and the multiplier bootstrap in finite samples, we provide two types of plots:
\begin{itemize}
	\item \emph{Kolmogorov-Smirnov distance.} We plot side-by-side boxplots of the Kolmogorov-Smirnov distances between the estimated distributions of the $\ell_p$-statistic $T_{n,p}$ and (a) the Gaussian proxy statistic, $\widetilde{T}_{n,p}$, (b) the Gaussian parametric bootstrap statistic based on the naive sample covariance, $T^*_{n,p, \mathrm{naive}}$, (c) the Gaussian parametric bootstrap statistic based on the thresholding estimate $\widehat{\Sigma}_n^+(\hat{\lambda})$, $T^*_{n,p,\lambda}$, (d) the Gaussian multiplier bootstrap, $T^g_{n,p}$. These boxplots give insight into the overall quality of the bootstrap procedures. Note that (b) and (d) are the same, but are implemented by two different algorithms.
	
	Since the true distribution of the $\ell_p$-statistic $T_{n,p}$ is unknown, we evaluate it based on 5000 Monte Carlo samples. To estimate the distributions of $\widetilde{T}_{n,p}$, $T^*_{n,p, \mathrm{naive}}$, $T^*_{n,p,\lambda}$, and $T^g_{n,p}$ we generate 1000 Monte Carlo samples of $X = \{X_i \in \mathbb{R}^d, 1\leq i \leq n\}$ and 1000 bootstrap samples for each Monte Carlo sample $X = \{X_i \in \mathbb{R}^d, 1\leq i \leq n\}$. We report results for sample size $n =  200$, dimension $d= 1000$, and exponents $p \in\{1, 2, \log d, \infty\}$.
	
	\item \emph{Lower tail probabilities.} We plot point estimates of $\mathrm{P}\left(T_{n,p} \leq q_{0.95}\right)$, where $q_{0.95}$ is the 95\% quantile of the distribution of $\widetilde{T}_{n,p}$, $T^*_{n,p, \mathrm{naive}}$, $T^*_{n,p,\lambda}$, and $T^g_{n,p}$, respectively. These point estimates clarify the pointwise accuracy of the bootstrap procedures. They can also be interpreted as the relative frequencies of the coverage of 95\% simultaneous confidence sets for the parameter $\mu = \mathrm{E}[X_1]$ under $H_0: \mu = 0$.
	
	We estimate these probabilities as follows. First, we draw 1000 Monte Carlo samples $X^{(1)}, \ldots, X^{(1000)}$, where $X^{(m)}= \{X_i^{(m)} \in \mathbb{R}^d, 1 \leq i \leq n\}$, and compute the associated $\ell_p$-statistics $T_{n,p}^{(m)}$, $1 \leq m \leq 1000$. For each Monte Carlo sample $X^{(m)}$, we generate 1000 bootstrap samples and construct bootstrap estimates $\hat{q}_{0.95}^{(m)}$ of the 95\% quantile of the distributions of $\widetilde{T}_{n,p}$, $T^*_{n,p, \mathrm{naive}}$, $T^*_{n,p,\lambda}$, and $T^g_{n,p}$, respectively. Then, we estimate $\mathrm{P}\left(T_{n,p} \leq q_{0.95}\right)$ as $1000^{-1} \sum_{m=1}^{1000} \mathbf{1}\{T_{n,p}^{(m)} \leq \hat{q}_{0.95}^{(m)}\}$. Again we report results for sample size $n =  200$, dimension $d= 1000$, and exponents $p \in\{1, 2, \log d, \infty\}$.
\end{itemize}

\begin{figure}[h!]
	\centering
	
	\input{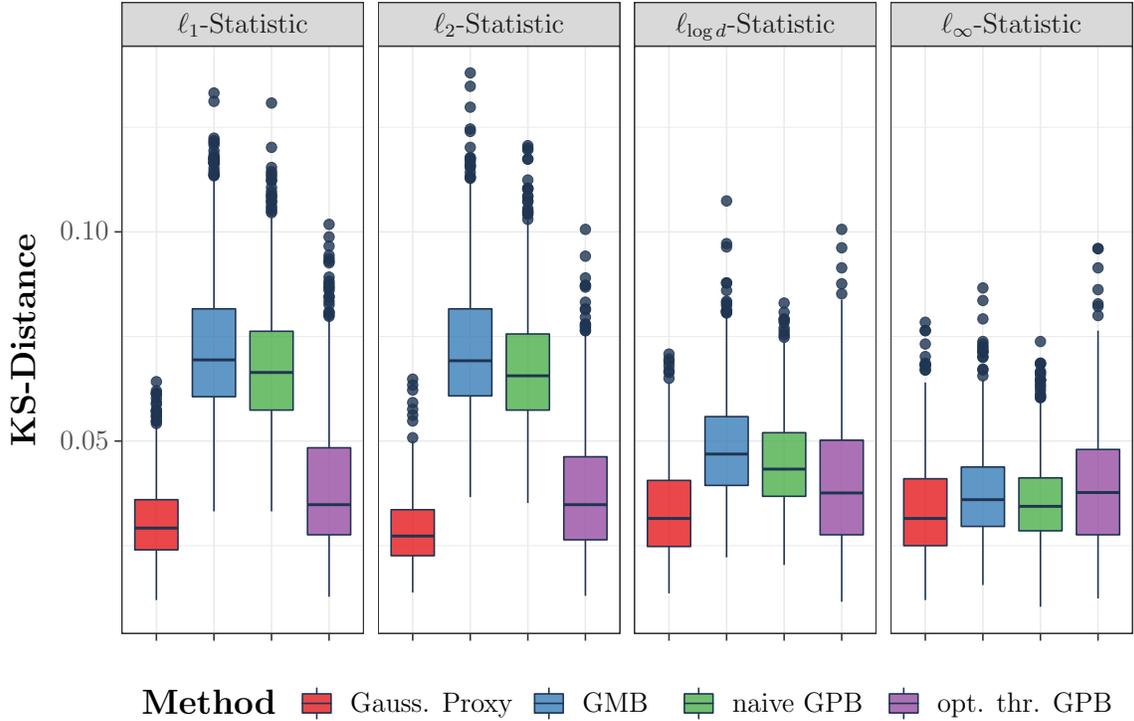}
	\caption{Boxplots of 1000 Kolmogorov-Smirnov distances between the distribution of the $\ell_p$-statistic and its bootstrap estimates based on Gaussian Proxy (Gauss. Proxy), Gaussian Multiplier Bootstrap (GMB), Naive Gaussian Parametric Bootstrap (naive GPB), and Gaussian Parametric Bootstrap based on $\widehat{\Sigma}_n^+(\hat{\lambda})$ (opt. thr. GPB). Sample size $n=200$, dimension $d= 1000$, $F$ cdf of Uniform$(-1,1)$. }
	\label{fig:1}
\end{figure}

\begin{figure}[h!]
	\centering
	\input{est_distr_KS_boxplots-n=200-d=1000-B=1000-nsim=1000-June-SLRHT.tex}
	\caption{Boxplots of 1000 Kolmogorov-Smirnov distances between the distribution of the $\ell_p$-statistic and its bootstrap estimates based on Gaussian Proxy (Gauss. Proxy), Gaussian Multiplier Bootstrap (GMB), Naive Gaussian Parametric Bootstrap (naive GPB), and Gaussian Parametric Bootstrap based on $\widehat{\Sigma}_n^+(\hat{\lambda})$ (opt. thr. GPB). Sample size $n=200$, dimension $d= 1000$, $F$ cdf of $t_4$. }
	\label{fig:2}
\end{figure}

\begin{figure}[h!]
	\centering
	\input{est_distr_coverage_boxplots-n=200-d=1000-B=1000-nsim=1000-June-SLRHU.tex}
	\caption{Relative frequencies of the simultaneous coverage of 1000 95\% confidence sets under $H_0: \mu = 0$ of the Gaussian Proxy (Gauss. Proxy), Gaussian Multiplier Bootstrap (GMB), Naive Gaussian Parametric Bootstrap (naive GPB), and Gaussian Parametric Bootstrap based on $\widehat{\Sigma}_n^+(\hat{\lambda})$ (opt. thr. GPB). The vertical bars indicate Monte Carlo standard errors. Sample size $n=200$, dimension $d= 1000$, $F$ cdf of Uniform$(-1,1)$.}
	\label{fig:3}
\end{figure}

\begin{figure}[h!]
	\centering
	\input{est_distr_coverage_boxplots-n=200-d=1000-B=1000-nsim=1000-June-SLRHT.tex}
	\caption{Relative frequencies of the simultaneous coverage of 1000 95\% confidence sets $H_0: \mu = 0$ of the Gaussian Proxy (Gauss. Proxy), Gaussian Multiplier Bootstrap (GMB), Naive Gaussian Parametric Bootstrap (naive GPB), and Gaussian Parametric Bootstrap based on $\widehat{\Sigma}_n^+(\hat{\lambda})$ (opt. thr. GPB). The vertical bars indicate Monte Carlo standard errors. Sample size $n=200$, dimension $d= 1000$, $F$ cdf of $t_4$. }
	\label{fig:4}
\end{figure}

In the following discussion the Gaussian proxy statistic $\widetilde{T}_{n,p}$ serves as an oracle estimator. It tells us how good the bootstrap procedures could be if we knew the true covariance matrix. Any difference between $\widetilde{T}_{n,p}$ and the other statistics solely arises from the different estimates of the covariance matrix.

Figure~\ref{fig:1} shows that if the data has light tails, the distribution of the Gaussian proxy statistic $\widetilde{T}_{n,p}$ provides an excellent approximation of the distribution of the $\ell_p$-statistic $T_{n,p}$ for all $p \in \{1, 2, \log d , \infty\}$. Moreover, the distribution of the Gaussian parametric bootstrap statistic $T^*_{n,p,\lambda}$ based on $\widehat{\Sigma}_n^+(\hat{\lambda})$ yields a comparably good approximation to the truth. In contrast, the distributions of the Gaussian multiplier statistic and the naive Gaussian parametric bootstrap are significantly poorer approximations to the truth. For large exponents $p \in \{\log d, \infty\}$ all four bootstrap approximations perform similarly. Thus, this plot fully supports every aspect of the theoretical results derived in Sections~\ref{subsec:BootstrapConsistency} and~\ref{subsec:BootstrapConsistencyStructure}.

Figure~\ref{fig:2} shows that if the data has heavy tails, Gaussian proxy statistic $\widetilde{T}_{n,p}$ and the bootstrap procedures yield poorer approximations to the truth. In particular, we see that the quality of the approximation worsens substantially as the exponent $p$ increases. This further corroborates the theoretical results derived in Sections~\ref{subsec:BootstrapConsistency} and~\ref{subsec:BootstrapConsistencyStructure}.

Figures~\ref{fig:3} and~\ref{fig:4} tell a similar, but more nuanced, story. From Figure~\ref{fig:3} we infer that if the data has light tails, the 95\% quantiles of the distributions of the Gaussian proxy statistic $\widetilde{T}_{n,p}$ and the Gaussian parametric bootstrap statistic $T^*_{n,p,\lambda}$ based on $\widehat{\Sigma}_n^+(\hat{\lambda})$ yield good approximations to the 95\% quantile of the true distribution for all exponents $p \in \{1, 2, \log d, \infty\}$. From Figure~\ref{fig:4} we learn that if the data has heavy tails, the approximations are fairly good for small exponents $p \in \{1, 2\}$, but fail spectacularly for large exponents $p \in \{\log d , \infty\}$. Moreover, Gaussian multiplier and naive Gaussian parametric bootstrap yield accurate estimates of the 95\% quantiles of the target distribution only for large exponents $p \in \{\log d, \infty\}$ and only when the data has light tails. This again supports the theoretical results from Sections~\ref{subsec:BootstrapConsistency} and~\ref{subsec:BootstrapConsistencyStructure}.

\subsection{Power of the bootstrap hypothesis test}

\begin{figure}[h!]
	\centering
	\input{power_plot_thr_GPB-n=200-d=400-B=1000-nsim=1000-SHU-L2.tex}
	\caption{Power functions under dense alternatives based on 1000 Monte Carlo samples. The gray bands indicate the Monte Carlo standard errors. Sample size $n=200$, dimension $d= 400$, $F$ cdf of Uniform$(-1,1)$.}
	\label{fig:5}
\end{figure}

\begin{figure}[h!]
	\centering
	\input{power_plot_thr_GPB-n=200-d=400-B=1000-nsim=1000-SHU-Lmax.tex}
	\caption{Power functions under sparse alternatives based on 1000 Monte Carlo samples. The gray bands indicate the Monte Carlo standard errors. Sample size $n=200$, dimension $d= 400$, $F$ cdf of Uniform$(-1,1)$.}
	\label{fig:6}
\end{figure}

To illustrate the effect of the exponent $p$ on size and power of the bootstrap hypothesis test we consider its power function in the following two high-dimensional testing scenarios:
\begin{itemize}
	\item \emph{Dense alternatives.} We test $H_0: \mu = 0$ vs. $H_1: \mu = \mu(\delta) \equiv \delta (1, \ldots, 1)' \in \mathbb{R}^d$ at a 5\% significance level. The signal strength $\delta$ is of order $O\left(1/\sqrt{nd}\right)$.
	\item \emph{Sparse alternatives.} We test $H_0: \mu = 0$ vs. $H_1: \mu = \mu(\delta) \equiv \delta (1, \ldots, 1,0, \ldots, 0)' \in \mathbb{R}^d$ at a 5\% significance level. The alternative has $2 \lceil \sqrt{\log d}/2\rceil$ non-zero entries and the signal strength $\delta$ is of order $O\left(\sqrt{ (\log d)/n}\right)$.
\end{itemize}
In Figures~\ref{fig:5} and~\ref{fig:6} we plot Monte Carlo estimates of the power function $\beta(\delta) = \mathrm{P}_{\mu(\delta)}\big(T_{n,p} > c^*_{n,p}(0.95)\big)$, where $c^*_{n,p}(0.95) = \inf\left\{t \in \mathbb{R} : \mathrm{P}(S_{n,p}^* \leq t \mid X)  \geq 0.05 \right\}$ and $S_{n,p}^*$ is the Gaussian parametric bootstrap test statistic based on $\widehat{\Sigma}_n^+(\hat{\lambda})$. The Monte Carlo estimate of $\beta(\delta)$ is based on 1000 Monte Carlo samples of $X = \{X_i \in \mathbb{R}^d, 1\leq i \leq n\}$ and 1000 bootstrap samples for each observed $X = \{X_i \in \mathbb{R}^d, 1\leq i \leq n\}$. The specific estimation procedure is identical to the one used to compute the lower tail probabilities in Section~\ref{subsec:PerformanceBootstrap}. We report results for sample size $n =  200$, dimension $d= 400$, exponents $p \in\{1, 2, \log d, \infty\}$, and light-tailed data.

Figure~\ref{fig:5} shows the power function for dense alternatives. We observe that tests based on $S_{n,1}$ and $S_{n,2}$ (they are nearly indistinguishable in the figure) are more powerful than those based on $S_{n,\log d}$ and $S_{n,\infty}$. This fully matches the theoretical predictions from Section~\ref{subsec:Power}. Figure~\ref{fig:6} displays the power function for sparse alternatives. In this case the bootstrap tests based on $S_{n,\log d}$ and $S_{n,\infty}$ are more powerful than those based on $S_{n,1}$ and $S_{n,2}$. The power functions associated with $S_{n,\log d}$ and $S_{n,\infty}$ are essentially the same with $S_{n,\log d}$ being slightly more powerful because of a larger constant (note $\mu$ has $2 \lceil \sqrt{\log d}/2\rceil = 4$ non-zero entries). Again, these findings fit well into the discussion in Section~\ref{subsec:Power}.

Lastly, at $\delta = 0$ the power functions of all four tests are about 0.05 in both Figures~\ref{fig:5} and~\ref{fig:6}. Thus, all four tests successfully control the type I error a 5\% significance level. For large values of $\delta$ all four tests unanimously reject the null hypothesis with probability (close to) one. This confirms the results from Sections~\ref{subsec:AsympCorrectness} and~\ref{subsec:ConsistencyTesting}.

\section{Conclusion}\label{sec:Conclusion}
In this paper we have introduced the Gaussian parametric bootstrap to estimate the distribution of $\ell_p$-statistics of high-dimensional random vectors. The procedure is versatile and user-friendly, since its implementation requires only a positive semi-definite estimate of the population covariance matrix. The main theoretical contributions state the consistency of the Gaussian parametric bootstrap under various conditions on the covariance structure of the data. To showcase the applicability of the Gaussian parametric bootstrap we propose a bootstrap hypothesis test for simultaneous inference on high-dimensional mean vectors. We discuss in detail asymptotic correctness, confidence sets, consistency under high-dimensional alternatives, and power of the test.

One of the current challenges in theoretical statistics is to understand when bootstrap procedures work in high-dimensional problems. At least for bootstrapping $\ell_p$-statistics of high-dimensional random vectors we can now give a definitive answer. The technical results in the appendix to this paper clarify that the success of bootstrapping $\ell_p$-statistics hinges on three factors: (a) $\ell_p$-norms of high-dimensional random vectors satisfy a Berry-Esseen-type central limit theorem under relatively mild moment conditions; (b) $\ell_p$-norms of Gaussian random vectors satisfy powerful anti-concentration inequalities; (c) the distributions of $\ell_p$-statistics under centered Gaussian distributions vary smoothly over their covariance matrices.

\expandafter\def\expandafter\appendixpagename%
\expandafter{\expandafter\Large\appendixpagename}
\begin{appendices}
\emph{Organization.} The appendices are divided into two parts. In Appendix~\ref{apx:MainTheory} we present additional results and technical lemmas. These include an abstract Berry-Esseen-type CLT for $\ell_p$-statistics (Theorem~\ref{theorem:BerryEsseen-Even-p}), Gaussian anti-concentration inequalities (Theorems~\ref{theorem:AntiConcentration-LpNorm-Even-p} and \ref{theorem:AntiConcentration-LpNorm}), Gaussian comparison inequalities (Theorem~\ref{theorem:Gaussian-Comparison-KolmogorovDistance}), smoothing inequalities (Section~\ref{subsec:SmoothingInequalities}), and auxiliary results for proving bootstrap consistency (Section~\ref{subsec:AuxResults}), for testing high-dimensional mean vectors (Section~\ref{subsec:AuxResultsApplication}), and concerning the partial derivatives of $\ell_p$-norms (Section~\ref{subsec:AuxResultsDerivatives}). In Appendix~\ref{apx:Proofs} we provide proofs to all results from the main text and the appendix.

\emph{Additional Notation.} We denote by $C^k(\mathbb{R}^d)$ the class of $k$ times continuously differentiable functions from $\mathbb{R}^d$ to $\mathbb{R}$, and by $C^k_b(\mathbb{R}^d)$ the class of all functions $f \in C^k(\mathbb{R}^d)$ with bounded support. For a real-valued matrix $M \in \mathbb{R}^{d \times d}$ we write $\|A\|_*$ to denote its nuclear norm (the $\ell_1$-norm of its singular values). For $\varepsilon > 0$ and a Borel set $A \in \mathcal{B}(\mathbb{R})$ define the $\varepsilon$-enlargement of $A$ as $A^\varepsilon := \{t \in \mathbb{R} : \inf_{s \in A} |t-s| \leq \varepsilon\}$. Moreover, we write $A^{-\varepsilon}$ to denote sets $B \in \mathcal{B}(\mathbb{R})$ for which  $B^{\varepsilon} = A$.

\section{Additional results and technical lemmas}\label{apx:MainTheory}
\subsection{Abstract Berry-Esseen-type CLT}\label{subsec:BerryEsseen}
Recall the setup from Section~\ref{sec:MainResults} in the main text. Consider a sequence $X=\{X_i\}_{i=1}^n$ of independent and centered random vectors in $\mathbb{R}^d$. Let $Z =\{Z_i\}_{i=1}^n$ be a sequence of independent multivariate Gaussian random vectors $Z_i \sim N(0, \mathrm{E}[X_iX_i'])$ which are independent of $X$. Define the scaled averages
\begin{align}\label{eq:subsec:BerryEsseen-1}
S_n^X : = \frac{1}{\sqrt{n}}\sum_{i=1}^nX_i \hspace{30pt} \mathrm{and} \hspace{30pt} S_n^Z : = \frac{1}{\sqrt{n}}\sum_{i=1}^nZ_i,
\end{align}
and the Kolmogorov-Smirnov distance
\begin{align}\label{eq:subsec:BerryEsseen-2}
\varrho_{n,p} := \sup_{t \geq 0} \left| \mathrm{P}\left(\|S_n^X\|_p \leq t\right) - \mathrm{P}\left(\|S_n^Z\|_p \leq t\right) \right|.
\end{align}

The upper bound in the original univariate Berry-Esseen inequality depends on the third moments of the $X_i$'s. In the multivariate case, the concept of third moments is less clear cut. For example, $\mathrm{E}[\|X\|^3]$ (for some norm $\| \cdot\|$) and $\sum_{k=1}^d\mathrm{E}[|X_k|^3]$ are both sensible generalizations of the univariate third moment. The bound in our Berry-Essen-type CLT depends on the following generalized third moments: For $a, b \geq 0$ arbitrary,
\begin{align}
M_{n, b}(a) & := \mathrm{E}\left[\frac{1}{n}\sum_{i=1}^n \Big( \|X_i\|_b^31\{\|X_i\|_b > a \} +  \|Z_i\|_b^31\{\|Z_i\|_b > a \}  \Big)\right], \label{eq:subsec:BerryEsseen-3}\\
L_{n,b} & := \mathrm{E}\left[\frac{1}{n}\sum_{i=1}^n \Big( \|X_i\|_{b}^3+  \|Z_i\|_{b}^3 \Big)\right].\label{eq:subsec:BerryEsseen-4}
\end{align}
We write $\overline{L}_{n,b}$ for an upper bound on $L_{n,b}$. Furthermore, we need two quantities based on the (average) covariance matrix of the $X_i$'s: the vector of its diagonal elements and its rank, i.e.
\begin{align}\label{eq:subsec:BerryEsseen-5}
\sigma_n^2 := \left(\mathrm{E}\left[\frac{1}{n}\sum_{i=1}^n X_{ik}^2\right]\right)_{k=1}^d \hspace{30pt} \mathrm{and} \hspace{30pt} r_n := \mathrm{rank}\left(\mathrm{E}\left[\frac{1}{n}\sum_{i=1}^n X_iX_i'\right]\right).
\end{align}

The following Berry-Esseen-type CLT for $\ell_p$-norms is our main theoretical contribution and central to all other results in this paper.

\begin{theorem}[Berry-Esseen-type CLT for $\ell_p$-norms]\label{theorem:BerryEsseen-Even-p}
	\noindent
		
	\begin{itemize}
		\item[(i)] For all $p \in [1, \infty)$ and all $\tau \in[1, \infty]$,
		\begin{align}\label{eq:theoremBerryEssen-Even-p-1}
		\varrho_{n,p}\lesssim \frac{M_{n,\tau p}\big( p^{1-1/(3\tau)}n^{1/3}\overline{L}_{n,\tau p}^{1/3}\big)}{p^{1 - 1/\tau} \overline{L}_{n,\tau p}} + \frac{(pd^{1/p})^{1-1/(3\tau)}\overline{L}_{n,\tau p}^{1/3}}{n^{1/6}}\frac{p^{1/2}r_n^{1/(2p)}}{\|\sigma_n\|_p}.
		\end{align}
		\item[(ii)] For all $p  \in [\log d, \infty]$,
		\begin{align}\label{eq:theoremBerryEssen-Even-p-2}
		\varrho_{n,p}\lesssim\frac{M_{n,\infty}\big(n^{1/3} (\log d)^{2/3}  \overline{L}_{n, \infty}^{1/3}\big)}{\overline{L}_{n, \infty}} + \frac{(\log d)^{7/6}}{n^{1/6}}\frac{\overline{L}_{n, \infty}^{1/3}}{\|\sigma_n\|_\infty}.
		\end{align}
	\end{itemize}
\end{theorem}
Observe that the upper bound on the Kolmogorov-Smirnov distance exhibits qualitatively different behavior depending on the magnitude of the exponent $p$ and the tails of the distribution of the $X_i$'s.

What is of interest here is that the upper bound undergoes a phase transition from polynomial dependence on $d$ (or $r_n$) to logarithmic dependence in $d$ as the exponent $p$ crosses the threshold $\log d$. It is easy to verify that for $p \asymp \log d$ and $\tau = 1$ the bounds in~\eqref{eq:theoremBerryEssen-Even-p-1} and~\eqref{eq:theoremBerryEssen-Even-p-2} are of the same order and that for $p \gtrsim \log d$ and $\tau \geq 1$ the bound in~\eqref{eq:theoremBerryEssen-Even-p-2} is tighter than the one in~\eqref{eq:theoremBerryEssen-Even-p-1}. Thus, the bound in~\eqref{eq:theoremBerryEssen-Even-p-1} is only relevant for exponents $p \lesssim \log d$. In this case, the tails of the distribution of the $X_i$'s come into play and the nuisance parameter $\tau \in [1, \infty]$ can be used to trade off moment conditions versus fractional powers of dimension $d$. To illustrate the basic idea of how to use $\tau$, let us consider the two boundary cases $\tau \in \{1, \infty\}$.  Denote by $\bar{\sigma}_{n,\min}^2 := \min_{1\leq k \leq d}\sigma_{n,k}^2$ the smallest diagonal element of the (averaged) covariance matrix of the $X_i$'s. For $\tau = 1$ eq.~\eqref{eq:theoremBerryEssen-Even-p-1} simplifies to
\begin{align*}
\varrho_{n,p}\lesssim \frac{M_{n,p}\big( p^{2/3}n^{1/3}L_{n,p}^{1/3}\big)}{\overline{L}_{n,p}} + \frac{p^{7/6} r_n^{1/(2p)} }{n^{1/6}\bar{\sigma}_{n, \min}} \left(\frac{ \overline{L}_{n,p} }{d^{1/p}}\right)^{1/3},
\end{align*}
while for $\tau = \infty$ eq.~\eqref{eq:theoremBerryEssen-Even-p-1} reduces to
\begin{align*}
\varrho_{n,p}\lesssim \frac{M_{n,\infty}\big( p n^{1/3}L_{n,\infty}^{1/3}\big)}{p\overline{L}_{n,\infty}} + \frac{p^{3/2} r_n^{1/(2p)} }{n^{1/6}\bar{\sigma}_{n, \min}}  \overline{L}_{n,\infty}^{1/3}.
\end{align*}
Typically, the first term on the right hand side of each of the two displays can be bounded independently of $d$ and is negligible as $n \rightarrow \infty$. Therefore, to decide which one of the two bounds is (asymptotically) tighter, we need to determined whether $d^{-1/p}\overline{L}_{n,p} \gtrless \overline{L}_{n,\infty}$. Clearly, this depends on the tails of the distribution of the $X_i$'s. For concreteness, if the $X_i$'s are sub-Gaussian, then $\overline{L}_{n,\infty} \asymp (\log d)^{3/2}$ while $d^{-1/p}\overline{L}_{n,p} \asymp d^{2/p}$. In the main text in Section~\ref{subsec:GaussianApprox}, we provide simplified and ready-to-use bounds under various conditions on the tails of the distribution of the $X_i$'s.

Theorem~\ref{theorem:BerryEsseen-Even-p} is non-asymptotic and holds for all $n, d, p$. However, it is only relevant in high-dimensional settings since in low-dimensional settings, in which the dimension $d$ is fixed or grows much slower than the square root of the sample size $n$, there exist sharper results~\citep[e.g.][and references therein]{bhattacharya1977RefinementsMultiCLT, goetze1991RateConvMultCLT, bentkus2003DependenceBerryEsseen, raic2019MultBerryEsseen}. Since $n^{-1/6}$ is the minimax optimal rate for CLTs in infinite dimensional Banach spaces, it is likely that in high-dimensional settings the bound in Theorem~\ref{theorem:BerryEsseen-Even-p} is nearly optimal in terms of dependence on $n$~\citep[see discussion in][]{chernozhukov2017CLTHighDim, bentkus1985lower}.

Theorem~\ref{theorem:BerryEsseen-Even-p} is related to Theorem 2.1 in~\cite{chernozhukov2017CLTHighDim}, which is a Berry-Esseeen-type CLT for hyper-rectangles. In fact, our proof builds on their idea of combining Stein's leave-one-out approach with Slepian's smart-path-interpolation and iterative arguments due by~\cite{bolthausen1984estimate}. All major technical differences between their and our proof can be traced back to the specific behavior of our new anti-concentration, Gaussian comparison, and smoothing inequalities in the regime $p \leq \log d$. If one is interested in results for large exponents $p \geq \log d$ only, one can simply combine the original proof from~\cite{chernozhukov2017CLTHighDim} with our new anti-concentration and smoothing inequalities. Without further modifications of their arguments one then obtains the following slight improvement of~\eqref{eq:theoremBerryEssen-Even-p-2}.

\begin{proposition}[Refined Berry-Esseen-type CLT for $\ell_p$-norms with large exponents]\label{proposition:Berry-Esseen-Large-p}
	For all $p \in [\log d , \infty]$,
	\begin{align*}
	\varrho_{n,p} \lesssim \frac{M_{n,\infty}\big(n^{1/3} (\log d)^{-1/3}  L_{n, \max}^{1/3}\big)}{L_{n, \max}} + \frac{(\log d)^{7/6}}{n^{1/6}}\frac{L_{n, \max}^{1/3}}{\|\sigma_n\|_\infty},
	\end{align*}
	where $L_{n, \max} \geq \max_{1 \leq k \leq d} \frac{1}{n}\sum_{i=1}^n \mathrm{E}\left[|X_{ik}|^3\right]$.
\end{proposition}
The second term on the right hand side in the bound of Proposition~\ref{proposition:Berry-Esseen-Large-p} is clearly smaller than the corresponding term in the bound~\eqref{eq:theoremBerryEssen-Even-p-2}. However, under the primitive conditions in Section~\ref{subsec:GaussianApprox} and for $p \geq \log d$, the terms $\overline{L}_{n, \max}$ and $\overline{L}_{n, p}$ will only differ by a factor of order $o(\log d)$.

\subsection{Anti-concentration inequalities}\label{subsec:AntiConcentration}
We begin with the following basic result for $\ell_p$-norms of random vectors with log-concave probability measure when $p \in 2 \mathbb{N}$ is an even integer. For a random variable $X \in \mathbb{R}^d$ with law $\nu$ and $A \in \mathcal{B}(\mathbb{R}^d)$ define $\mathrm{P}_\nu(X \in A) := \int_A d\nu $.

\begin{theorem}\label{theorem:AntiConcentration-LpNorm-Even-p}
	Let $X, X' \in \mathbb{R}^d$ be i.i.d. random vectors with law $\nu$. For $\varepsilon > 0$ arbitrary,
	\begin{align*}
	\sup_{\nu}\sup_{p \in 2\mathbb{N}} \sup_{ t \geq 0 }\mathrm{P}_\nu\left( t \leq \|X\|_p \leq t + \varepsilon \: \left\| \|X\|_p -\|X'\|_p \right\|_{\psi_1}\right) \lesssim \varepsilon,
	\end{align*}
	where the supremum in $\nu$ is taken over all log-concave probability measures on $\mathbb{R}^d$.
\end{theorem}
The main feature of this inequality is that it is dimension free in the sense that the left hand side depends on dimension $d$ only through the quantity $\left\| \|X\|_p - \|X'\|_p \right\|_{\psi_1}$.

Interestingly, the assumption that $\nu$ belongs to the class of log-concave probability measures is indeed necessary: For one, it is easy to see that if $d=1$ and the class of probability measures contains measures $\nu$ whose densities $d\nu$ have multiple modes (or a point mass), there exists an $\varepsilon > 0$ for which the inequality is violated. Thus, the densities $d\nu$ have to be continuous and unimodal. For another, the term $\left\| \|X\|_p - \|X'\|_p \right\|_{\psi_1}$ is finite for all $p \in 2 \mathbb{N}$ only if $\nu$ has at least sub-exponential tails. Together these two facts imply that the $\nu$ has to be log-concave.

The key idea behind the proof of Theorem~\ref{theorem:AntiConcentration-LpNorm-Even-p} is that for $p \in 2 \mathbb{N}$ we may interpret $\|X\|_p^p$ as a (multivariate) polynomial and invoke the distributional version of the Carbery-Wright inequality for random polynomials over convex bodies~\citep[][Theorem 8]{carbery2001Distributional}. We defer the detailed proof to Appendix~\ref{apx:Proofs}.

Specializing to a Gaussian random vector and the $\ell_2$-norm, Theorem~\ref{theorem:AntiConcentration-LpNorm-Even-p} yields the following important result:
\begin{corollary}\label{corollary:AntiConcentration-L2Norm}
	Let $X \in \mathbb{R}^d$ be a Gaussian random vector with mean $\mu \in \mathbb{R}^d$ and positive semi-definite covariance matrix $\Sigma$. For $\varepsilon > 0$ arbitrary,
	\begin{align*}
	\sup_{ t \geq 0 }\mathrm{P}\left( t \leq \|X\|_2 \leq t + \varepsilon\left(\mathrm{tr}(\Sigma^2) + \mu'\Sigma \mu\right)^{1/4}\right) \lesssim \varepsilon.
	\end{align*}
\end{corollary}
Note that this corollary holds for any fixed mean $\mu \in \mathbb{R}^d$ of $X$. In this sense it is an anti-concentration inequality about ellipsoids with arbitrary radii but fixed center $\mu\in \mathbb{R}^d$. In contrast, the well-known result by~\cite{nazarov2003MaximaPerimeter} is an anti-concentration inequality over ellipsoids with arbitrary radii and arbitrary centers. However, this stronger result requires the additional assumption that $\Sigma$ is positive definite and comes with a smaller standard deviation proxy, namely $\mathrm{tr}(\Sigma^{-2})^{-1/4}$. Our Corollary~\ref{corollary:AntiConcentration-L2Norm} sharpens the finite-dimensional analogue of Theorem 2.7 in~\cite{goetze2019LargeBall} by introducing the quadratic term $\mu'\Sigma\mu$ to the inequality.

Combining Theorem~\ref{theorem:AntiConcentration-LpNorm-Even-p} with an interpolation argument and fine properties of Gaussian measures (i.e. Plancherel's identity and careful truncation) yields the following theorem for general $\ell_p$-norms with exponent $p \in [1, \infty]$.

\begin{theorem}\label{theorem:AntiConcentration-LpNorm}
	Let $X \in \mathbb{R}^d$ be a centered Gaussian random vector with positive semi-definite covariance matrix $\Sigma = (\sigma_{kj})_{k,j = 1}^d$ of rank $r \geq 1$. Set $\sigma^2 = (\sigma_{kk})_{k=1}^d$. For $\varepsilon > 0$ arbitrary,
	\begin{align*}
	\sup_{p \in [1, \infty]} \sup_{t \geq 0} \mathrm{P}\left( t \leq  \|X\|_p \leq t + \varepsilon \: \frac{\|\sigma\|_p}{\omega_p(d, r)} \right) \lesssim \varepsilon,
	\end{align*}
	where
	\begin{align*}
	\omega_p(d, r) =
	\begin{cases}
	\sqrt{p r^{1/p}} &\mathrm{if} \:\: p \in [1, \infty),\\
	\sqrt{\log d} &\mathrm{if} \:\:  p \geq \log d.
	\end{cases}
	\end{align*}
\end{theorem}
\begin{remark}\label{remark:theorem:AntiConcentration-LpNorm}
	We will use the following refined inequality to prove Theorem~\ref{theorem:BerryEsseen-Even-p}. Let $p_+ = 2\lceil\frac{p}{2}\rceil$ be the smallest even integer larger than $p$. Then,
	\begin{align*}
	\sup_{p \in [1, \infty)} \sup_{t \geq 0} \sup_{q \in \{p, p_+\}} \mathrm{P}\left( t \leq  \|X\|_q \leq t + \varepsilon \: \frac{\|\sigma\|_p}{\sqrt{p r^{1/p}}} \right) \lesssim \varepsilon.
	\end{align*}
	We only need to show the validity of the inequality for $q = p_+$. To this end, invoke Theorem~\ref{theorem:AntiConcentration-LpNorm-Even-p} and, as in the proof of Theorem~\ref{theorem:AntiConcentration-LpNorm}, lower bound $\left\| \|X\|_{p_+} -\|X'\|_{p_+} \right\|_{\psi_1}$ by $c_0\|\sigma\|_p/ \sqrt{p r^{1/p}}$, where $c_0 > 0$ is an absolute constant.
\end{remark}

Several comments are in order: First, the term $\omega_p(d, r)$ undergoes a phase transition as $p$ crosses the threshold $\log d$. This phenomenon matches well-known phase transitions of the expected value and the variance of $\ell_p$-norms of isotropic Gaussian random vectors~\citep[e.g.][]{schechtman1990volume, paouris2018Dvoretzky}. Second, for $p = \infty$ our theorem improves Nazarov's inequality~\citep[i.e.][Theorem 1]{chernozhukov2017Nazarov} in a crucial detail: The term $\omega_p(d, r)$ depends on the inverse of the largest diagonal element of $\Sigma$, whereas the corresponding quantity in Nazarov's inequality depends on the inverse of the smallest diagonal element of $\Sigma$.  Third, if $X$ is isotropic then $\|\sigma\|_2r^{-1/4} \asymp \mathrm{tr}(\Sigma^2)^{1/4}$ and Corollary~\ref{corollary:AntiConcentration-L2Norm} and Theorem~\ref{theorem:AntiConcentration-LpNorm} are asymptotically equivalent. However, in general, the inequality in Corollary~\ref{corollary:AntiConcentration-L2Norm} is tighter.

In special cases, it is possible to obtain explicit constants for inequalities in this section. Notably, if the Gaussian random vector $X$ has a spherical distribution, the bounds on variance and moments of $\ell_p$-norms of Gaussian random vectors in Section 3 of~\cite{paouris2018Dvoretzky} are directly applicable.

\subsection{Gaussian comparison inequalities}\label{subsec:GaussianComparison}
The following result allows us to compare the distributions of $\ell_p$-norms of two centered Gaussian random vectors with (potentially) different covariance matrices.

\begin{theorem}\label{theorem:Gaussian-Comparison-KolmogorovDistance}
	Let $X$ and $Y$ be two independent Gaussian random vectors in $\mathbb{R}^d$ with mean zero and covariance matrices $\Sigma^X = (\Sigma^X_{jk})_{j,k=1}^d$ and $\Sigma^Y = (\Sigma^Y_{jk})_{j,k=1}^d$, respectively. Define $\sigma^2_X = (\Sigma^X_{kk})_{k=1}^d$, $\sigma^2_Y = (\Sigma^Y_{kk})_{k=1}^d$, $r_X = \mathrm{rank}(\Sigma^X)$, and $r_Y = \mathrm{rank}(\Sigma^Y)$. Set $\Delta_{op} = \big\|\Sigma^X - \Sigma^Y \big\|_{op}$ and $\Delta_p = \big\|\mathrm{vec}(\Sigma^X - \Sigma^Y) \big\|_p$.
	\begin{itemize}
		\item[(i)] For all $p \in [1, \infty)$,
		\begin{align}\label{eq:theorem:Gaussian-Comparison-KolmogorovDistance-1}
		\sup_{t \geq 0 } \Big|\mathrm{P}\left(\|X\|_p \leq t\right) - \mathrm{P}\left(\|Y\|_p \leq t\right) \Big| \lesssim \frac{\sqrt{p^2 d^{1/p} r_X^{1/p} \Delta_p} }{\|\sigma_X\|_p} \bigwedge \frac{ \sqrt{p^2 d^{1/p} r_Y^{1/p} \Delta_p} }{\|\sigma_Y\|_p}.
		\end{align}
		\item[(ii)] For all $p \in [\log d, \infty]$,
		\begin{align}\label{eq:theorem:Gaussian-Comparison-KolmogorovDistance-2}
		\sup_{t \geq 0 } \Big|\mathrm{P}\left(\|X\|_p \leq t\right) - \mathrm{P}\left(\|Y\|_p \leq t\right) \Big| \lesssim  \frac{ (\log d) \sqrt{\Delta_{op} \wedge \Delta_\infty} }{ \|\sigma_X\|_\infty \vee \|\sigma_Y\|_\infty}.
		\end{align}
	\end{itemize}
\end{theorem}
\begin{remark}
	In order to derive minimax lower bounds for the Gaussian parametric and Gaussian multiplier bootstrap it would be extremely useful to have complementary lower bounds on the Kolmogorov-Smirnov distance between the distributions of $\|X\|_p$ and $\|Y\|_p$.
\end{remark}
\begin{remark}
	Using Corollary~\ref{corollary:AntiConcentration-L2Norm} instead of Theorem~\ref{theorem:AntiConcentration-LpNorm} in the proof of above theorem we obtain the following alternative bound for $p =2$:
	\begin{align}\label{eq:remark:theorem:Gaussian-Comparison-KolmogorovDistance-1}
	\sup_{t \geq 0 } \Big|\mathrm{P}\left(\|X\|_2 \leq t\right) - \mathrm{P}\left(\|Y\|_2 \leq t\right) \Big| \lesssim \sqrt{\frac{d^{1/2} \Delta_2 }{\|\mathrm{vec}(\Sigma^X)\|_2 \vee \|\mathrm{vec}(\Sigma^Y)\|_2}}.
	\end{align}
\end{remark}
The most interesting aspect of this result is that the upper bound on the Kolmogorov-Smirnov distance shows qualitatively different behavior depending on the magnitude of the exponent $p$. Let $\sigma_{X,\min}^2: = \min_{1 \leq k \leq d}\Sigma_{kk}^X$ or $\sigma_{Y,\min}^2 := \min_{1 \leq k \leq d}\Sigma_{kk}^Y$. We can now further simplify~\eqref{eq:theorem:Gaussian-Comparison-KolmogorovDistance-1} to
\begin{align*}
\sup_{t \geq 0 } \Big|\mathrm{P}\left(\|X\|_p \leq t\right) - \mathrm{P}\left(\|Y\|_p \leq t\right) \Big| \lesssim \sqrt{\frac{p^2r_X^{1/p}}{d^{1/p}}\frac{\Delta_p}{\sigma_{X, \min}^2}} \bigwedge \sqrt{\frac{p^2r_Y^{1/p}}{d^{1/p}}\frac{\Delta_p}{\sigma_{Y, \min}^2}}.
\end{align*}
This bound is useful because it depends on the dimension $d$ only via the difference $\Delta_p$ (note that $r_X d^{-1}, r_Y d^{-1} \leq 1$). It is therefore a key ingredient for proving consistency of the naive Gaussian multiplier and the Gaussian parametric bootstrap statistics.

For $p=2$~\cite{goetze2019LargeBall} (Theorem 2.1 and corollaries) have derived bounds similar to~\eqref{eq:remark:theorem:Gaussian-Comparison-KolmogorovDistance-1} but based on a completely different approach that involves bounding the density functions of $\|X\|_2^2$ and $\|Y\|_2^2$. In the one- and two-dimensional cases their bounds are strictly tighter than~\eqref{eq:remark:theorem:Gaussian-Comparison-KolmogorovDistance-1}. In the $d$-dimensional case with $d \geq 3$, their bound is (roughly) of the order of $O\left(\|\Sigma^X - \Sigma^Y\|_*\big(\|\mathrm{vec}(\Sigma^X)\|_2 \vee \|\mathrm{vec}(\Sigma^Y)\|_2 \big)^{-1}\right)$ which is just slightly smaller than the square of~\eqref{eq:remark:theorem:Gaussian-Comparison-KolmogorovDistance-1}. In general, neither their nor our bound is clearly better or worse.

For $p = \infty$ and $(\log d)^2\Delta_\infty = o(1)$ inequality~\ref{eq:theorem:Gaussian-Comparison-KolmogorovDistance-2} improves Theorem 2 in~\cite{chernozhukov2015ComparisonAnti} in two ways: First, we improve the rate from $(\log d)^{2/3} \Delta_\infty^{1/3}$ to $(\log d) \Delta_\infty^{1/2}$. Second, our bound depends only on the inverse of $\|\sigma_X\|_\infty$ or $\|\sigma_Y\|_\infty$, whereas the inequality by~\cite{chernozhukov2015ComparisonAnti} depends also on the the inverse of either $\sigma_{X,\min}^2$ or $\sigma_{Y,\min}^2$. The second improvement can be ascribed to our improved anti-concentration inequality.

The proofs of Theorem~\ref{theorem:BerryEsseen-Even-p} and~\ref{theorem:Gaussian-Comparison-KolmogorovDistance} are conceptually very similar and rely on the same anti-concentration and smoothing inequalities. The main difference between the two proofs is that Theorem~\ref{theorem:BerryEsseen-Even-p} uses a second-order Taylor approximation to expand the smoothed Kolmogorov-Smirnov distance and matches the first two moments of $S^X_n$ and $S^Z_n$, whereas Theorem~\ref{theorem:Gaussian-Comparison-KolmogorovDistance} uses only a first-order Taylor approximation and matches only the first moments (because the second moments $\Sigma^X$ and $\Sigma^Y$ differ). Along the way, the proof of Theorem~\ref{theorem:Gaussian-Comparison-KolmogorovDistance} also makes heavily use of $X$ and $Y$ being Gaussian.

\subsection{Smoothing inequalities}\label{subsec:SmoothingInequalities}
Smoothing inequalities allow us to replace probabilities like $\mathrm{P}\left(\|X\|_p \in A \right) = \mathrm{E}[\mathbf{1}\{\|X\|_p \in A\}]$, which are expectations of non-differentiable indicator functions of non-differentiable maps $x \mapsto \|x\|_p$, by expectations of smooth functions. This enables us to approximate these probabilities via first- or second-order Taylor approximations, which is the first step in establishing the abstract Berry-Essen-type and Gaussian comparison inequalities in Sections~\ref{subsec:BerryEsseen} and~\ref{subsec:GaussianComparison}.

\begin{lemma}[$C^\infty_b(\mathbb{R}^d)$-Approximation of $\ell_p$-Norms]\label{lemma:SmoothIndicator-C-Infty-Function-Even-p} Let $X \in \mathbb{R}^d$ be an arbitrary random vector. There exists a family of smooth functions $\mathcal{H} = \{h_{p,d,\beta, \delta, A} \in C^{\infty}_b(\mathbb{R}^d) : p \in 2\mathbb{N}, d, \beta, \delta > 0, A \in  \mathcal{B}(\mathbb{R}) \}$ which satisfies the following:
	\begin{itemize}
		\item[(i)] For $A \in  \mathcal{B}(\mathbb{R})$, $p \in 2\mathbb{N}$, $\tau \in [1, \infty]$,  and $\kappa_p = 3\beta^{-1} pd^{1/(\tau p)}$,
		\begin{align}\label{eq:lemma:SmoothIndicator-C-Infty-Function-Even-p-1}
		\mathrm{P}\left(\|X\|_p \in A \right) \leq \mathrm{E}\left[h_{p, d, \beta, \delta, A^{\kappa_p}}\left(X\right) \right] \leq \mathrm{P}\left(\|X\|_p \in A^{3\delta + 2\kappa_p} \right).
		\end{align}
		\item[(ii)] The functions in $\mathcal{H}$ have support set $\big\{x \in \mathbb{R}^d: M_{\beta}(x) \in A^{3\delta} \setminus A\big\}$, where $M_{\beta}$ is defined in eq.~\eqref{eq:DefM-2}.
		\item[(iii)] For $p \in 2\mathbb{N}$, $\tau \in [1, \infty]$, and $q= \frac{\tau p}{\tau p-1}$,
		\begin{align}\label{eq:lemma:SmoothIndicator-C-Infty-Function-Even-p-2}
		\begin{split}
		\sup_{A \in \mathcal{B}(\mathbb{R})}\left\|\left(\sum_{|\alpha| = 2} \left|D^\alpha h_{p, d, \beta, \delta, A}\right|^q\right)^{1/q}\right\|_\infty &\lesssim \left(\frac{1}{\delta^2} + \frac{\beta}{\delta} \right)d^{2(\tau-1)/(\tau p)},\\
		\sup_{A \in \mathcal{B}(\mathbb{R})}\left\| \left(\sum_{|\alpha| = 3} \left|D^\alpha  h_{p, d, \beta, \delta, A} \right|^q\right)^{1/q}\right\|_\infty &\lesssim \left(\frac{1}{\delta^3}+ \frac{\beta}{\delta^2} + \frac{\beta^2}{\delta}\right)d^{3(\tau-1)/(\tau p)}
		\end{split}
		\end{align}
		\item[(iv)] Let $\mathcal{I} = \{A \subseteq \mathbb{R} : A = [0,t], t \geq 0\}$. Then, for $A \in \mathcal{I}$,  $\tau \in [1, \infty]$, and $p \in [1, \infty)$,
		\begin{align}\label{eq:lemma:SmoothIndicator-C-Infty-Function-Even-p-3}
		\mathrm{P}\left(\|X\|_p \in A \right) \leq \mathrm{E}\left[h_{p_+, d, \beta, \delta, A^{\kappa_{p_+}}}\left(X\right) \right] \leq \mathrm{P}\left(\|X\|_p \in A^{3\delta + 4\kappa_p} \right),
		\end{align}
		where $p_+ = 2 \lceil\frac{p}{2}\rceil$ is the smallest even integer larger than $p$ and $\kappa_p = 3\beta^{-1} pd^{1/(\tau p)}$.
	\end{itemize}
\end{lemma}
The key observation behind this result is that $\ell_p$-norms with even exponents $p \in 2\mathbb{N}$ are (multivariate) polynomials of degree $p$ and continuously differentiable (except at 0) with self-normalizing derivatives.

Lemma~\ref{lemma:SmoothIndicator-C-Infty-Function-Even-p} $(iv)$ holds for all $p \in [1, \infty)$ but for $p \rightarrow \infty$ it is impossible to simultaneously control the (probability of the) enlarged set $A^{3\delta + 4\kappa_p}$ in~\eqref{eq:lemma:SmoothIndicator-C-Infty-Function-Even-p-3} and the partial derivatives in~\eqref{eq:lemma:SmoothIndicator-C-Infty-Function-Even-p-2}: For one, if we set $\beta = O(p)$, we can control $\kappa_p$ and $\mathrm{P}(\|X\|_p \in A^{3\delta + 4\kappa_p})$ but the bounds on the partial derivatives diverge. For another, if we set $\beta = O(1)$, we can control the partial derivatives but $\mathrm{P}(\|X\|_p \in A^{3\delta + 4\kappa_p}) \rightarrow 1$. The reason for this is that $\ell_p$-norms with large exponents (relative to dimension $d$) behave essentially like the non-differentiable maximum norm ($\ell_\infty$-norm). We therefore have to smooth $\ell_p$-norms with large exponents $p \geq \log d$ differently. This is content of the next result.

\begin{lemma}[$C^\infty_b(\mathbb{R}^d)$-Approximation of $\ell_p$-Norms for $p \geq \log d$]\label{lemma:SmoothIndicator-C-Infty-Function-Large-p} Let $X \in \mathbb{R}^d$ be an arbitrary random vector. There exists a family of smooth functions $\mathcal{H} = \{h_{p,d,\beta, \delta, A} \in C^{\infty}_b(\mathbb{R}^d) : p \in [\log d, \infty], d, \beta, \delta > 0, A \in  \mathcal{B}(\mathbb{R}) \}$ which satisfies the following:
	\begin{itemize}
		\item[(i)] For $A \in  \mathcal{B}(\mathbb{R})$, $p \in [\log d , \infty]$, and $\kappa = e\beta^{-1} \log (2d)$,
		\begin{align}\label{eq:lemma:SmoothIndicator-C-Infty-Function-Large-p-1s}
		\mathrm{P}\left(\|X\|_p \in A \right) \leq \mathrm{E}\left[h_{p, d, \beta, \delta, A^{\kappa}}\left(X\right) \right] \leq \mathrm{P}\left(\|X\|_p \in A^{3\delta + 2\kappa} \right).
		\end{align}
		\item[(ii)] The functions in $\mathcal{H}$ have support set $\big\{x \in \mathbb{R}^d: F_{\beta}(x) \in A^{3\delta} \setminus A\big\}$, where $F_{\beta}$ is the smooth-max function and defined in eq.~\eqref{eq:DefF}.
		\item[(iii)] For $p \in [\log d, \infty]$,
		\begin{align}\label{eq:lemma:SmoothIndicator-C-Infty-Function-Large-p-2s}
		\begin{split}
		\sup_{A \in \mathcal{B}(\mathbb{R})}\left\|\sum_{|\alpha| = 2} \left|D^\alpha h_{p, d, \beta, \delta, A}\right|\right\|_\infty &\lesssim \frac{1}{\delta^2} + \frac{\beta}{\delta},\\
		\sup_{A \in \mathcal{B}(\mathbb{R})}\left\|\sum_{|\alpha| = 3} \left|D^\alpha  h_{p, d, \beta, \delta, A} \right|\right\|_\infty &\lesssim \frac{1}{\delta^3}+ \frac{\beta}{\delta^2} + \frac{\beta^2}{\delta}.
		\end{split}
		\end{align}	
	\end{itemize}
\end{lemma}
The smooth approximation in this lemma is based on the smooth-max function that was first introduced by~\cite{chernozhukov2013GaussianApproxVec}. Thus, the modest novelty of this result is that the smooth-max function can be used to approximate not just the maximum-norm ($\ell_\infty$-norm) but also all $\ell_p$-norms with exponent $p \geq \log d$. Smooth approximations based on the smooth-max function satisfy other useful stability properties beyond the bounds on the second and third derivative in~\eqref{eq:lemma:SmoothIndicator-C-Infty-Function-Large-p-2s}. While these other stability properties are crucial to the proofs of the Berry-Esseen-type CLTs in~\cite{chernozhukov2013GaussianApproxVec, chernozhukov2015ComparisonAnti, chernozhukov2017CLTHighDim, deng2020beyond, koike2019notes}, these properties are not essential to our proofs.
\subsection{Auxiliary results I (Bootstrap consistency)}\label{subsec:AuxResults}
In this section we collect bounds on tail probabilities and moments of vectors and covariance matrices in $\ell_p$-norms. These results are used in Sections~\ref{subsec:BootstrapConsistency} and~\ref{subsec:BootstrapConsistencyStructure}.

Throughout this section we write $\bar{X}_n = n^{-1}\sum_{i=1}^n X_i$ for the sample mean and $\widehat{\Sigma}_{\mathrm{naive}} = (\widehat{\sigma}_{kj})_{k,j=1}^d = n^{-1}\sum_{i=1}^n (X_i - \bar{X}_n)(X_i - \bar{X}_n)'$ for the sample covariance matrix.

\begin{lemma}\label{lemma:MomentsOfLpNorm}
	Let $X \in \mathbb{R}^d$ be a random vector that satisfies Assumption~\ref{assumption:FiniteMoments} with $s \geq t \vee p$ for some $t \geq 0$. Then,
	\begin{align*}
	\left(\mathrm{E}\|X\|^t_p\right)^{1/t} \lesssim K_s \|\sigma\|_p,
	\end{align*}
	where $\sigma = (\sigma_k)_{k=1}^d$ and $\sigma_k^2 = \mathrm{E}[X_k^2]$ for $1 \leq k \leq d$.
\end{lemma}

\begin{lemma}[Product of subgaussian random variables]\label{lemma:ProductSubgaussian}
	Let $X_1, \ldots, X_K \in \mathbb{R}$ be sub-gaussian random variables. Then, for $K \geq 1$,
	\begin{align*}
	\left\|\prod_{k=1}^K X_k\right\|_{\psi_{2/K}} \leq \prod_{k=1}^K \left\|X_k \right\|_{\psi_2}.
	\end{align*}
\end{lemma}

\begin{lemma}[Sub-Gaussian]\label{lemma:BoundsCovariance} Let $X=\{X_i\}_{i=1}^n$ be a sequence of i.i.d. mean zero random vectors in $\mathbb{R}^d$ with covariance matrix $\Sigma = (\sigma_{jk})_{j,k=1}^d$ which satisfies Assumption~\ref{assumption:SubGaussian}. Set $\sigma^2 = (\sigma_{kk})_{k=1}^d$ and $\widehat{\sigma}^2 = (\widehat{\sigma}_{kk})_{k=1}^d$. Let $\zeta \in (0,1)$ be arbitrary.
	\begin{itemize}
		\item[(i)] With probability at least $1 - 2 \zeta$, for all $p \in [1, \infty]$,
		\begin{align*}
		\|\mathrm{vec}(\widehat{\Sigma}_{\mathrm{naive}} - \Sigma)\|_p \lesssim \|\sigma\|_p^2\left(\sqrt{\frac{\log d + \log (2/\zeta)}{n}} \bigvee \frac{\log d + \log (2/\zeta)}{n} \right).
		\end{align*}
		\item[(ii)] Let $\mathrm{r}(\Sigma) = \mathrm{tr}(\Sigma)/ \|\Sigma\|_{op}$ be the effective rank of $\Sigma$. With probability at least $1 - 2\zeta$,
		\begin{align*}
		\|\widehat{\Sigma}_{\mathrm{naive}} - \Sigma\|_{op} \lesssim \|\Sigma\|_{op} \left(\sqrt{\frac{\mathrm{r}(\Sigma)\log d + \log (2/\zeta)}{n}} \bigvee \frac{\mathrm{r}(\Sigma)\log d + \log (2/\zeta)}{n} \right).
		\end{align*}
		\item[(iii)] With probability at least $1 - 2 \zeta$,
		\begin{align*}
		\max_{1 \leq k \leq d} \big|(\widehat{\sigma}_k/\sigma_k)^2 - 1\big|\lesssim  \left(\sqrt{\frac{\log d + \log (2/\zeta)}{n}} \bigvee \frac{\log d + \log (2/\zeta)}{n} \right).
		\end{align*}
	\end{itemize}
\end{lemma}
\begin{remark}
	For $p= 2$  case (i) is qualitatively (i.e up to log-factors) identical to Theorem 2.1 in~\cite{bunea2015on}. Cases (ii) and (iii) are folklore. Note that these results are usually given for the matrix of second moments $n^{-1}\sum_{i=1}^nX_iX_i'$ only, whereas we provide bounds for the sample covariance matrix $\widehat{\Sigma}_{\mathrm{naive}}$.
\end{remark}

\begin{lemma}[Finite Moments]\label{lemma:BoundsCovariance-FiniteMoments} Let $X=\{X_i\}_{i=1}^n$ be a sequence of i.i.d. mean zero random vectors in $\mathbb{R}^d$ with covariance matrix $\Sigma = (\sigma_{jk})_{j,k=1}^d$. Set $\sigma^2 = (\sigma_{kk})_{k=1}^d$ and $\widehat{\sigma}^2 = (\widehat{\sigma}_{kk})_{k=1}^d$.
	\begin{itemize}
		\item[(i)] Suppose Assumption~\ref{assumption:FiniteMoments} holds with $s \geq  p \vee 4$. For $p \in [1, \infty]$,
		\begin{align*}
		\|\mathrm{vec}(\widehat{\Sigma}_{\mathrm{naive}} - \Sigma)\|_p = O_p\left( K_s^2 \|\sigma\|_p^2 \sqrt{\frac{ p \wedge \log d}{n}}\right).
		\end{align*}
		\item[(ii)] Suppose Assumption~\ref{assumption:FiniteMoments} holds with $s = 2$. Set $\mathrm{m}(\Sigma) = \mathrm{E}[\max_{1\leq i\leq n}\|X_i\|_2^2]/\|\Sigma\|_{op}$.
		\begin{align*}
		\|\widehat{\Sigma}_{\mathrm{naive}} - \Sigma\|_{op} = O_p\left(\|\Sigma\|_{op}\left(\sqrt{\frac{\mathrm{m}(\Sigma)\log(d \wedge n)}{n}} \bigvee \frac{\mathrm{m}(\Sigma)\log(d \wedge n)}{n} \right) \right).
		\end{align*}
		\item[(iii)] Suppose Assumption~\ref{assumption:FiniteMoments} holds for $s \geq 2$. Set $ \widetilde{\mathrm{m}}(\Sigma) = \mathrm{E}[\max_{1\leq i\leq n}\|\mathrm{diag}(\Sigma)^{-1}X_i\|_2^2]$.
		\begin{align*}
		\max_{1 \leq k \leq d} \big|(\widehat{\sigma}_k/\sigma_k)^2 - 1\big| =
		\begin{cases}
		O_p\left( K_s^2 d^{1/s} \sqrt{\frac{ s \wedge \log d}{n}}\right) & \mathrm{for}\:\: s \geq 4,\\
		O_p\left(\sqrt{\frac{\widetilde{\mathrm{m}}(\Sigma) \log(d \wedge n)}{n}} \bigvee \frac{\widetilde{\mathrm{m}}(\Sigma) \log(d \wedge n)}{n} \right)& \mathrm{for}\:\: s \geq 2.
		\end{cases}
		\end{align*}
	\end{itemize}
\end{lemma}

Covariance matrix estimators that can exploit Assumption~\ref{assumption:ApproxSparseCovariance} are the so-called thresholding estimator. In the following, we consider a generic thresholding operator $T_\lambda: \mathbb{R} \rightarrow \mathbb{R}$ with thresholding parameter $\lambda$, which satisfies
\begin{align*}
(i) \:\:\: |T_\lambda(u)| \leq |u|; \hspace{35pt} (ii)\:\:\: T_\lambda(u) = 0 \:\:\mathrm{for}\:\: |u| \leq \lambda; \hspace{35pt} (iii) \:\:\:|T_\lambda(u) - u| \leq \lambda.
\end{align*}
Thesholding operators satisfying these three properties include the hard-thresholding operator, $T_\lambda(u) = u \mathbf{1}\{|u| > \lambda\}$ (from Section~\ref{sec:MainResults}) as well as the soft-thresholding operator, $T_\lambda(u) = \mathrm{sign}(u)\big( (|u| - \lambda) \vee 0\big)$.

\begin{lemma}[Thresholded covariance estimators]\label{lemma:BoundsThresholdEstimators} Let $X=\{X_i\}_{i=1}^n$ be a sequence of i.i.d. mean zero random vectors in $\mathbb{R}^d$ with covariance matrix $\Sigma = (\sigma_{jk})_{j,k=1}^d$. Set $\sigma^2 = (\sigma_{kk})_{k=1}^d$ and $\widehat{\sigma}^2 = (\widehat{\sigma}_{kk})_{k=1}^d$. Let $\zeta \in (0,1)$ be arbitrary and set $\lambda_n \asymp \sqrt{\frac{\log d + \log (2/\zeta)}{n}} \bigvee \frac{\log d + \log (2/\zeta)}{n}$.
	\begin{itemize}
		\item[(i)] Suppose Assumption~\ref{assumption:SubGaussian} holds. Let $A \in \mathbb{R}^{d \times d}$ be the adjacency matrix of $\Sigma$, i.e. $A_{jk} = \mathbf{1}\{\sigma_{jk}\neq 0\}$. With probability at least $1 - 2 \zeta$, for all $p \in [1, \infty]$,
		\begin{align*}
		\big\|\mathrm{vec}\big(T_{\lambda_n}(\widehat{\Sigma}_{\mathrm{naive}}) - \Sigma\big)\big\|_p &\lesssim \|\mathrm{vec}(A)\|_p \|\sigma\|_\infty^2 \lambda_n, \\
		\big\|T_{\lambda_n}(\widehat{\Sigma}_{\mathrm{naive}}) - \Sigma\big\|_{op} &\lesssim \|\mathrm{vec}(A)\|_{op} \|\sigma\|_\infty^2\lambda_n.
		\end{align*}
		
		\item[(ii)] Suppose Assumptions~\ref{assumption:SubGaussian} and~\ref{assumption:ApproxSparseCovariance} hold. With probability at least $1 - 2 \zeta$, for all $p \in [\theta, \infty]$,
		\begin{align*}
		\big\|\mathrm{vec}\big(T_{\lambda_n}(\widehat{\Sigma}_{\mathrm{naive}}) - \Sigma\big)\big\|_p  &\lesssim d^{1/p} R_{\gamma, p} \|\sigma\|_\infty^{2(1-\gamma)}\lambda_n^{1-\gamma},\\
		\big\|T_{\lambda_n}(\widehat{\Sigma}_{\mathrm{naive}}) - \Sigma\big\|_{op} &\lesssim R_{\gamma, \theta} \|\sigma\|_\infty^{2(1-\gamma)}\lambda_n^{1-\gamma}.
		\end{align*}
		
		\item[(iii)] Suppose Assumptions~\ref{assumption:FiniteMoments} holds with $s \geq (p \wedge \log d) \vee 4$. Then, for all $p \in [2 \vee \theta, \infty]$,
		\begin{align*}
		\big\|\mathrm{vec}\big(T_{\lambda_n}(\widehat{\Sigma}) - \Sigma\big)\big\|_p  &= O_p\left( \|\mathrm{vec}(A)\|_p K_s^2 \|\sigma\|_s^2 \sqrt{\frac{ s \wedge \log d}{n}}\right),\\
		\big\|T_{\lambda_n}(\widehat{\Sigma}) - \Sigma\big\|_{op} &=  O_p\left( \|A\|_{op} K_s^2 \|\sigma\|_s^2 \sqrt{\frac{ s \wedge \log d}{n}}\right).
		\end{align*}
		
		\item[(iv)] Suppose Assumptions~\ref{assumption:FiniteMoments} and~\ref{assumption:ApproxSparseCovariance} hold with $s \geq (p \wedge \log d) \vee 4$. Then, for all $p \in [2 \vee \theta, \infty]$,
		\begin{align*}
		\big\|\mathrm{vec}\big(T_{\lambda_n}(\widehat{\Sigma}) - \Sigma\big)\big\|_p  &= O_p\left( d^{1/p} R_{\gamma, p}K_s^{2(1-\gamma)} \|\sigma\|_s^{2(1-\gamma)} \left(\frac{ s \wedge \log d}{n}\right)^{(1-\gamma)/2}\right),\\
		\big\|T_{\lambda_n}(\widehat{\Sigma}) - \Sigma\big\|_{op} &=   O_p\left(R_{\gamma, \theta}K_s^{2(1-\gamma)} \|\sigma\|_s^{2(1-\gamma)} \left(\frac{ s \wedge \log d}{n}\right)^{(1-\gamma)/2}\right).
		\end{align*}
	\end{itemize}
\end{lemma}
\begin{remark}
	Cases (i) and (ii) generalize Theorems 6.23 and 6.27 in~\cite{wainwright2019high} to the sample covariance matrix $\widehat{\Sigma}_{\mathrm{naive}}$ and the vectorized $\ell_p$-norm with $p \in [\theta, \infty]$. Note that for $\gamma= 0$ case (ii) reduces to case (i). Cases (iii) and (iv) bounds are only better than the naive bounds from Lemma~\ref{lemma:BoundsCovariance-FiniteMoments} for $2p < s$. For $s \asymp \log d$ they match the bounds of cases (i) and (ii).
\end{remark}

\begin{remark}
	For $\max_{1 \leq j\leq d} \sum_{k=1}^d |\sigma_{jk}|^\gamma\leq R_{1,\gamma}$ and $p = 2$ we obtain the same rate as Theorem 2 in~\cite{bickel2008regularization}
\end{remark}

\begin{lemma}[Banded covariance estimators]\label{lemma:BoundsBandedEstimators}  Let $X=\{X_i\}_{i=1}^n$ be a sequence of i.i.d. mean zero random vectors in $\mathbb{R}^d$ with covariance matrix $\Sigma = (\sigma_{jk})_{j,k=1}^d$. Set $\sigma^2 = (\sigma_{kk})_{k=1}^d$ and $\widehat{\sigma}^2 = (\widehat{\sigma}_{kk})_{k=1}^d$. Let $\zeta \in (0,1)$ be arbitrary and set $\ell_n = B_p^{p/(1 + p\alpha)} \|\sigma\|_\infty^{-2p/(1 + p\alpha)}\lambda_n^{-p/(1 + p\alpha)}$ for $\lambda_n, \alpha > 0$.
	\begin{itemize}
		\item[(i)] Suppose Assumptions~\ref{assumption:SubGaussian} and~\ref{assumption:ApproxBandableCovariance} hold. Set $\lambda_n \asymp \sqrt{\frac{\log d + \log (2/\zeta)}{n}} \bigvee \frac{\log d + \log (2/\zeta)}{n}$. With probability at least $1 - 2 \zeta$, for all $p \in [\theta, \infty]$,
		\begin{align*}
		\|\mathrm{vec}(B_{\ell_n}(\widehat{\Sigma}_{\mathrm{naive}}) - \Sigma)\|_p &\lesssim d^{1/p}B_p^{1/(1 + p\alpha)} \|\sigma\|_\infty^{2p\alpha/(1+p\alpha)}  \lambda_n^{p\alpha/(1 + p\alpha)}, \\
		\big\|B_{\ell_n}(\widehat{\Sigma}_{\mathrm{naive}}) - \Sigma\big\|_{op} &\lesssim B_\theta^{1/(1 + \alpha\theta)} \|\sigma\|_\infty^{2\alpha\theta/(1+\alpha\theta)}  \lambda_n^{\alpha\theta/(1 + \alpha\theta)}.
		\end{align*}
		
		\item[(ii)] Suppose Assumptions~\ref{assumption:FiniteMoments} and~\ref{assumption:ApproxBandableCovariance} hold with $s \geq (p \wedge \log d) \vee 4$. Set $\lambda_n \asymp \sqrt{\frac{s \wedge \log d}{n}}$. For all $p \in [2 \vee \theta, \infty]$,
		\begin{align*}
		\big\|\mathrm{vec}\big(B_{\ell_n}(\widehat{\Sigma}_{\mathrm{naive}}) - \Sigma\big)\big\|_p &= O_p\left(d^{1/p}B_p^{1/(1 + p\alpha)} K_s^{2p\alpha/(1+p\alpha)} \|\sigma\|_s^{2p\alpha/(1+p\alpha)} \lambda_n^{p\alpha/(1 + p\alpha)} \right), \\
		\big\|B_{\ell_n}(\widehat{\Sigma}_{\mathrm{naive}}) - \Sigma\big\|_{op} &= O_p\left(B_\theta^{1/(1 + \alpha\theta)} K_s^{2\alpha\theta/(1+\alpha\theta)} \|\sigma\|_s^{2\alpha\theta/(1+\alpha\theta)} \lambda_n^{\alpha\theta/(1 + \alpha\theta)} \right).
		\end{align*}
	\end{itemize}
\end{lemma}
\begin{remark}
	Analogous results also hold for covariance estimates based on the tapering operator as defined in~\cite{bickel2008covariance}.
\end{remark}

\subsection{Auxiliary results II (Testing high-dimensional mean vectors)}\label{subsec:AuxResultsApplication}
In this section we present results that are used in Section~\ref{sec:SimultaneousTesting}. Recall that the Gaussian parametric bootstrap estimate of the $\alpha$-quantile of test statistic $S_{n,p}$ is
\begin{align*}
c^*_{n,p}(\alpha) := \inf\left\{t \in \mathbb{R} : \mathrm{P}(S_{n,p}^* \leq t \mid X)  \geq \alpha \right\}.
\end{align*}
We now introduce the Gaussian approximation of the $\alpha$-quantile of test statistic $S_{n,p}$,
\begin{align*}
\tilde{c}_p(\alpha) := \inf\left\{t \in \mathbb{R} : \mathrm{P}(\widetilde{S}_p \leq t )  \geq \alpha \right\},
\end{align*}
where $\widetilde{S}_p := \|V^1\|_p$ with $V^1 \sim N(0, \Omega)$. Given the theoretical results in Section~\ref{sec:MainResults} we can expect that these two quantiles are close. The next lemma formalizes this intuition. It is a straightforward adaptation of Lemma 3.2 in~\cite{chernozhukov2013GaussianApproxVec} to our setup.
\begin{lemma}[Comparison of quantiles]\label{lemma:ComparisonQuantiles}	For all $p \in [1, \infty)$ and all $\delta > 0$,
	\begin{align}\label{eq:lemma:ComparisonQuantiles-1}
	\begin{split}
	&\sup_{\alpha \in (0,1)} \mathrm{P}\Big(c^*_{n,p}(\alpha) \leq \tilde{c}_{n,p}\big(\pi_p(\delta) + \alpha\big)\Big) \geq 1 - \mathrm{P}\left(\Pi_p > \delta\right), \\
	&\sup_{\alpha \in (0,1)} \mathrm{P}\Big(\tilde{c}_{n,p}(\alpha) \leq c^*_{n,p}\big(\pi_p(\delta) + \alpha\big)\Big) \geq 1 - \mathrm{P}\left(\Pi_p > \delta\right),
	\end{split}
	\end{align}
	where
	\begin{align*}
	\pi_p^2(\delta) = \begin{cases}
	\frac{p^2r_\omega^{1/p}}{{d'}^{1/p}}\frac{\delta}{\omega_{\min}^2} &\mathrm{if} \:\: p \in [1, \infty),\\
	(\log^2 d') \frac{ \delta }{\omega_{\max}^2} &\mathrm{if} \:\:  p \geq \log d'
	\end{cases}
	\hspace{20pt} \mathrm{and} \hspace{20pt}
	\Pi_p = \begin{cases}
	\widehat{\Gamma}_p &\mathrm{if} \:\: p \in [1, \infty),\\
	\widehat{\Gamma}_{op} \wedge \widehat{\Gamma}_\infty  &\mathrm{if} \:\:  p \geq \log d'.
	\end{cases}
	\end{align*}
\end{lemma}

The following lemma provides bounds on the upper quantiles of the Gaussian proxy statistic $\widetilde{S}_p$ in terms of its expected value and the covariance matrix.

\begin{lemma}[Bounds on (upper) quantiles of $\widetilde{S}_p$]\label{lemma:UpperBoundQuantiles}
	For all $\alpha \in (0,1/2]$,
	\begin{align*}
	 \mathrm{E}[\widetilde{S}_p] - \|\Omega^{1/2}\|_{2 \rightarrow p} \wedge \sqrt{\mathrm{Var}[\widetilde{S}_{p}]} \leq \tilde{c}_p(1-\alpha) \leq \mathrm{E}[\widetilde{S}_p] + \sqrt{2\log(1/\alpha)} \|\Omega^{1/2}\|_{2 \rightarrow p} \wedge \sqrt{(1/\alpha) \mathrm{Var}[\widetilde{S}_p]}.
	\end{align*}
	In fact, the upper bound holds for all $\alpha \in (0, 1)$.
\end{lemma}
\begin{remark}\label{remark:lemma:UpperBoundQuantiles}
	Note that one may combine this lemma with Lemma~\ref{lemma:ComparisonQuantiles} to obtain bounds on the bootstrap critical values $c^*_{n,p}(1- \alpha)$ for $\alpha \in (0, 1/2]$.
\end{remark}

\subsection{Auxiliary results III (Partial derivatives of $\ell_p$-norms)}\label{subsec:AuxResultsDerivatives}
Here we collect three lemmata on the partial derivatives of $\ell_p$-norms.

\begin{lemma}[Partial Derivatives]\label{lemma:Lp-Norm-Derivatives} Let $p > 1$ and $x \in \big\{z \in \mathbb{R}^d: z_i \geq 0, i=1, \ldots, d \big\} \setminus \{0\}$. Set $M_p(x)= \big(\sum_{j=1}^d x_j^p\big)^{1/p}$. The following are the partial derivatives of $M_p$ up to order three:
	\begin{align*}
	\frac{\partial M_p(x)}{\partial x_k} &= \frac{x_k^{p-1}}{\big(M_p(x)\big)^{p-1}}\\
	\frac{\partial^2 M_p(x)}{\partial x_k^2} &= \frac{(p-1)x_k^{p-2}}{M_p(x) \big)^{p-1}} - \frac{(p-1) x_k^{2p-2}}{M_p(x)\big)^{2p-1}}\\
	\frac{\partial^2 M_p(x)}{\partial x_k \partial x_\ell} &= -\frac{(p-1) x_k^{p-1} x_\ell^{p-1}  }{\big(M_p(x)\big)^{2p -1}}\\
	\frac{\partial^3 M_p(x)}{\partial x_k \partial x_\ell \partial x_m} &= \frac{(2p-1)(p-1) x_k^{p-1} x_\ell^{p-1} x_m^{p-1}  }{\big(M_p(x)\big)^{3p-1}}\\
	\frac{\partial^3 M_p(x)}{\partial x_k^2 \partial x_\ell} &= -\frac{(p-1)^2 x_k^{p-2} x_\ell^{p-1} }{\big(M_p(x)\big)^{2p-1}} + \frac{(2p-1)(p-1) x_k^{2p-2} x_\ell^{p-1}  }{\big(M_p(x)\big)^{3p-1}}\\
	\frac{\partial^3M_p(x)}{\partial x_k^3} &= \frac{(p-1)(p-2) x_k^{p-3} }{\big(M_p(x)\big)^{p-1}} - \frac{3 (p-1)^2 x_k^{2p-3} }{\big(M_p(x)\big)^{2p-1}} + \frac{(2p-1)(p-1)x_k^{3p-3}  }{\big(M_p(x)\big)^{3p-1}}.
	\end{align*}
\end{lemma}

\begin{lemma}[Stability of Partial Derivatives (Conjugate Norm)]\label{lemma:Lp-Norm-Derivatives-Stability} Let $p,q > 1$ be conjugate exponents such that $1/p + 1/q =1$. Let $x \in \big\{z \in \mathbb{R}^d: z_i \geq 0, i=1, \ldots, d \big\} \setminus \{0\}$ and set $M_p(x)= \big(\sum_{j=1}^d x_j^p\big)^{1/p}$. We have the following:
	\begin{align*}
	\left(\sum_{k=1}^d \left|\frac{\partial M_p(x)}{\partial x_k}\right|^q\right)^{1/q} &= 1\\
	\left(\sum_{k=1}^d \left|\frac{\partial^2 M_p(x)}{\partial x_k^2}\right|^q\right)^{1/q} &\leq \frac{2(p-1)d^{1/p}}{M_p(x)} + \frac{2(p-1)}{M_p(x)} \hspace{20pt} \forall p \geq 2\\
	\left(\sum_{k, \ell} \left|\frac{\partial^2 M_p(x)}{\partial x_k \partial x_\ell}\right|^q\right)^{1/q} &\leq \frac{(p-1)}{M_p(x)}\\
	\left(\sum_{k, \ell, m} \left|\frac{\partial^3 M_p(x)}{\partial x_k \partial x_\ell \partial x_m}\right|^q\right)^{1/q} &\leq \frac{(2p-1)(p-1)}{M_p^2(x)}\\
	\left(\sum_{k, \ell} \left|\frac{\partial^3 M_p(x)}{\partial x_k^2 \partial x_\ell}\right|^q\right)^{1/q} &\leq \frac{2(p-1)^2d^{1/p}}{M_p^2(x)} + \frac{2(2p-1)(p-1)}{M_p^2(x)} \hspace{20pt} \forall p \geq 2\\
	\left(\sum_{k=1}^d \left|\frac{\partial^3 M_p(x)}{\partial x_k^3}\right|^q\right)^{1/q} &\leq \frac{4(p-1)(p-2)d^{2/p}}{M_p^2(x)} + \frac{12(p-1)}{M_p^2(x)} + \frac{4(2p-1)(p-1)}{M_p^2(x)}\hspace{20pt} \forall p \geq 3\\
	\left(\sum_{k=1}^d \left|\frac{\partial^3 M_2(x)}{\partial x_k^3}\right|^2\right)^{1/2} &\leq \frac{24}{M_2^2(x)}.
	\end{align*}
\end{lemma}

\begin{lemma}[Stability of Partial Derivatives (Transformed Conjugate Norm)]\label{lemma:Lp-Norm-Derivatives-Stability-3} Let $p > 1, \tau \geq 1$, and $q'= \frac{\tau p}{\tau p -1}$. Let $x \in \big\{z \in \mathbb{R}^d: z_i \geq 0, i=1, \ldots, d \big\} \setminus \{0\}$ and set $M_p(x)= \big(\sum_{j=1}^d x_j^p\big)^{1/p}$. We have the following:
	\begin{align*}
	\left(\sum_{k=1}^d \left|\frac{\partial M_p(x)}{\partial x_k}\right|^{q'}\right)^{1/q'} &\leq d^{(\tau-1)/(\tau p)}\\
	\left(\sum_{k=1}^d \left|\frac{\partial^2 M_p(x)}{\partial x_k^2}\right|^{q'}\right)^{1/q'} &\leq \frac{2(p-1)d^{(2\tau-1)/(\tau p)}}{M_p(x)} + \frac{2(p-1)d^{(\tau-1)/(\tau p)}}{M_p(x)} \hspace{30pt} \forall p \geq 2\\
	\left(\sum_{k, \ell} \left|\frac{\partial^2 M_p(x)}{\partial x_k \partial x_\ell}\right|^{q'}\right)^{1/q'} &\leq \frac{(p-1)d^{2(\tau-1)/(\tau p)}}{M_p(x)}\\
	\left(\sum_{k, \ell, m} \left|\frac{\partial^3 M_p(x)}{\partial x_k \partial x_\ell \partial x_m}\right|^{q'}\right)^{1/q'} &\leq \frac{(2p-1)(p-1)d^{3(\tau-1)/(\tau p)}}{M_p^2(x)}\\
	\left(\sum_{k, \ell} \left|\frac{\partial^3 M_p(x)}{\partial x_k^2 \partial x_\ell}\right|^{q'}\right)^{1/q'} &\leq \frac{2(p-1)^2d^{(3\tau-2)/(\tau p)}}{M_p^2(x)} + \frac{2(2p-1)(p-1)d^{2(\tau-1)/(\tau p)}}{M_p^2(x)} \hspace{30pt} \forall p \geq 2\\
	\begin{split}
	\left(\sum_{k=1}^d \left|\frac{\partial^3 M_p(x)}{\partial x_k^3}\right|^{q'}\right)^{1/q'} &\leq \frac{4(p-1)(p-2)d^{(3\tau-1)/(\tau p)}}{M_p^2(x)} + \frac{12(p-1)d^{(\tau-1)/(\tau p)}}{M_p^2(x)}\\
	&\quad{} + \frac{4(2p-1)(p-1)d^{(\tau-1)/(\tau p)}}{M_p^2(x)}
	\end{split}
	\hspace{30pt} \forall p \geq 3\\
	\left(\sum_{k=1}^d \left|\frac{\partial^3 M_2(x)}{\partial x_k^3}\right|^{\frac{2\tau}{2\tau-1}}\right)^{\frac{2\tau-1}{2\tau}} &\leq \frac{24d^{(\tau-1)/(2\tau)}}{M_2^2(x)}.
	\end{align*}
\end{lemma}

\section{Proofs}\label{apx:Proofs}
\subsection{Proofs for Section~\ref{sec:MainResults}}
\begin{proof}[\textbf{Proof of Theorem~\ref{theorem:GaussianApprox}}]
	The proof of this theorem follows from Theorem~\ref{theorem:BerryEsseen-Even-p} and Proposition~\ref{proposition:Berry-Esseen-Large-p}. We first establish Case (iv). Cases (i) and (ii) are especial cases of (iv). Case (iii) has a standalone proof.
		
	\textbf{Proof of Case (iv).}
	For $s \geq 3$ define
	\begin{align*}
	\overline{L}_{n, \tau p}(s) := \left(\frac{1}{n}\sum_{i=1}^n\mathrm{E}\left[\|X_i\|_{\tau p}^s + \|Z_i\|_{\tau p}^s\right]\right)^{3/s}
	\end{align*}
	and note that by two applications H{\"o}lder's inequality $\overline{L}_{n, \tau p}(s) \geq L_{n, \tau p}$. Recall that Gaussian random vectors $Z \in \mathbb{R}^d$ satisfy Assumption~\ref{assumption:FiniteMoments} with $K_s = \sqrt{s}$ for all $s \geq 1$. Therefore, by Assumption~\ref{assumption:FiniteMoments}, Lemma~\ref{lemma:MomentsOfLpNorm}, and H{\"o}lder's inequality there exists an absolute constant $C > 1$ such that for all $s \geq 3$, 
	\begin{align*}
	\overline{L}_{n, \tau p}^{1/3}(s) \leq C \big(K_{\tau p \vee s} \vee \sqrt{\tau p \vee s}\big) d^{1/{\tau p}} \sigma_{n,\max}.
	\end{align*}
	Hence,
	\begin{align}\label{eq:theorem-Gaussian-Approximation-1}
	\frac{(pd^{1/p})^{1-1/(3\tau)}\overline{L}_{n, \tau p}^{1/3}(s)}{n^{1/6}\omega_p^{-1}(d, r_n)\|\sigma_n\|_p} \leq C \big(K_{\tau p \vee s} \vee \sqrt{\tau p \vee s}\big) \frac{(pd^{1/p})^{1-1/(3\tau)} \omega_p(d, r_n)}{n^{1/6}}\frac{d^{1/{\tau p}} \sigma_{n,\max}}{d^{1/p} \sigma_{n,\min}}.
	\end{align}
	Moreover, observe that for any real-valued random variable $Z$ and any $t > 0$ and $s \geq 3$, $\mathrm{E}[|Z|^3\mathbf{1}\{|Z| > t\}] \leq \mathrm{E}[|Z|^3(|Z|/t)^{s-3}\mathbf{1}\{|Z| > t\} ] \leq t^{3-s} \mathrm{E}[|Z|^s]$. Hence, Assumption~\ref{assumption:FiniteMoments} and Lemma~\ref{lemma:MomentsOfLpNorm} imply that for all $s \geq 3$,
	\begin{align}\label{eq:theorem-Gaussian-Approximation-2}
	\frac{M_{n, \tau p}\big(p^{1-1/(3\tau)}n^{1/3}\overline{L}_{n,\tau p}^{1/3}(s)\big)}{p^{1 - 1/\tau} \overline{L}_{n,\tau p}(s)} & \leq \frac{p^{3-1/\tau}n\overline{L}_{n,\tau p}(s)}{p^{s-s/(3\tau)}n^{s/3}\overline{L}_{n,\tau p}^{s/3}(s)} \frac{n^{-1}\sum_{i=1}^n\mathrm{E}\left[\|X_i\|_{\tau p}^s + \|Z_i\|_{\tau p}^s\right]}{p^{1 - 1/\tau} \overline{L}_{n,\tau p}(s)} \nonumber\\
	&=	\frac{p^{-s(1-1/\tau)/3}}{(p^2n)^{(s-3)/3}}.
	\end{align}
	Combining eq.~\eqref{eq:theorem-Gaussian-Approximation-1} and \eqref{eq:theorem-Gaussian-Approximation-2} with Theorem~\ref{theorem:BerryEsseen-Even-p} we conclude that for $s \geq 3$, $p \in [1, \infty)$ and $\tau \in [1, \infty]$,
	\begin{align}\label{eq:theorem-Gaussian-Approximation-3}
	\varrho_{n,p} &\lesssim \big(K_{\tau p \vee s} \vee \sqrt{\tau p \vee s}\big) \frac{(pd^{1/p})^{1-1/(3\tau)} \omega_p(d, r_n)}{n^{1/6}}\frac{d^{1/{\tau p}} \sigma_{n,\max}}{d^{1/p} \sigma_{n,\min}} +  \frac{p^{-s(1-1/\tau)/3}}{(p^2n)^{(s-3)/3}} \nonumber\\
	&\lesssim \big(K_{\tau p \vee s} \vee \sqrt{\tau p \vee s}\big) \sqrt{\frac{p^3 d^{4/(3\tau p)} r_n^{1/p}}{p^{2/(3\tau)} n^{1/3}} \frac{\sigma_{n,\max}^2}{\sigma_{n,\min}^2}} +  \frac{p^{-s(1-1/\tau)/3}}{(p^2n)^{(s-3)/3}}.
	\end{align}
	Since the $X$'s have only $s \geq 3$ moments, $K_{\tau p \vee s} < \infty$ only if $\tau p \leq s$. Since $\tau \geq 1$, this implies that $K_{\tau p \vee s} < \infty$ only if $p \leq s$. Hence, we deduce from~\eqref{eq:theorem-Gaussian-Approximation-3} that for all $p \in [1, s]$ and $s \geq 3$,
	\begin{align*}
	\varrho_{n,p} &\lesssim \big(K_s \vee \sqrt{s}\big) \sqrt{\frac{p^3 d^{4/(3s)} r_n^{1/p}}{n^{1/3}} \frac{\sigma_{n,\max}^2}{\sigma_{n,\min}^2}} + \frac{1}{(p^2n)^{(s-3)/3}}.
	\end{align*}
	For $s \geq 4$ this can be further simplified as in the statement of the theorem. This completes the proof of case (vi).

	\textbf{Proof of Case (i).} This is a special case of statement (iv); more precisely eq.~\eqref{eq:theorem-Gaussian-Approximation-3}. Recall that that if the $X$'s are sub-gaussian, then Assumption~\ref{assumption:FiniteMoments} holds for all $\tau p, s \geq 1$ and $K_{\tau p\vee s} = \sqrt{\tau p \vee s}$. Take $\tau = (\log d)/p \vee 1$ and $s= 6$ to obtain 
	\begin{align*}
	\big(K_{\tau p \vee s} \vee \sqrt{\tau p \vee s}\big) \sqrt{\frac{p^3 d^{4/(3\tau p)} r_n^{1/p}}{p^{2/(3\tau)} n^{1/3}} \frac{\sigma_{n,\max}^2}{\sigma_{n,\min}^2}} &\lesssim \sqrt{\log d} \sqrt{\frac{p^3 d^{4/(3\log d)} r_n^{1/p}}{p^{2p/(3\log d)} n^{1/3}} \frac{\sigma_{n,\max}^2}{\sigma_{n,\min}^2}} \\
	&\lesssim \sqrt{\frac{p^3 (\log d) r_n^{1/p}}{n^{1/3}} \frac{\sigma_{n,\max}^2}{\sigma_{n,\min}^2}},
	\end{align*}
	and
	\begin{align*}
	\frac{p^{-s(1-1/\tau)/3}}{(p^2n)^{(s-3)/3}}&\lesssim \frac{1}{p^2 n}.
	\end{align*}
	Thus, by eq.~\eqref{eq:theorem-Gaussian-Approximation-3},
	\begin{align*}
	\varrho_{n,p} \lesssim  \sqrt{\frac{p^3 (\log d) r_n^{1/p}}{n^{1/3}} \frac{\sigma_{n,\max}^2}{\sigma_{n,\min}^2}}.
	\end{align*}
	This completes the proof of case (i).
	
	\textbf{Proof of Case (ii).} This is again a special case of statement (iv). Note that if the $X$'s are sub-exponential, then Assumption~\ref{assumption:FiniteMoments} holds for all $\tau p, s \geq 1$ and $K_{\tau p\vee s} = \tau p \vee s$. Therefore, in eq.~\eqref{eq:theorem-Gaussian-Approximation-3} take $\tau = (\log d)/p \vee 1$ and $s= 6$ to obtain
	\begin{align*}
	\varrho_{n,p} \lesssim \sqrt{\frac{p^3 (\log d)^2 r_n^{1/p}}{n^{1/3}} \frac{\sigma_{n,\max}^2}{\sigma_{n,\min}^2}},
	\end{align*}
	This completes the proof of case (ii).
	
	\textbf{Proof of Case (iii).}
	For $s \geq 3$ define
	\begin{align*}
	L_{n,\max}(s) : = \max_{1 \leq k \leq d} \left(\frac{1}{n}\sum_{i=1}^n \mathrm{E}\left[|X_{ik}|^s + |Z_{ik}|^s\right]\right)^{3/s},
	\end{align*}
	and observe that by two applications of Jensen's inequality
	\begin{align*}
	L_{n,\max}(s) \geq \max_{1 \leq k \leq d} \frac{1}{n}\sum_{i=1}^n \mathrm{E}\left[|X_{ik}|^3\right].
	\end{align*}
	As in the proof of case (iv), we have
	\begin{align}\label{eq:theorem-Gaussian-Approximation-4}
	\frac{M_{n,\infty}\big(n^{1/3} (\log d)^{-1/3}  L_{n, \max}^{1/3}(s)\big)}{L_{n, \max}(s)} &\leq \frac{n (\log d)^{-1}L_{n, \max}(s)}{n^{s/3} (\log d)^{-s/3}  L_{n, \max}^{s/3}(s) } \frac{n^{-1}\sum_{i=1}^n\mathrm{E}\left[\|X_i\|_\infty^s + \|Z_i\|_\infty^s\right]}{L_{n, \max}(s)}\nonumber\\
	&\lesssim \frac{ (\log d)^{(4s-3)/3}}{n^{(s-3)/3}},
	\end{align}
	where the last inequality follows from Lemma 2.2.2 in~\cite{vandervaart1996weak} and Assumption~\ref{assumption:SubGaussian} or~\ref{assumption:SubExponential}. By Lemma~\ref{lemma:MomentsOfLpNorm}, Jensen's inequality, and Assumption~\ref{assumption:SubGaussian} or~\ref{assumption:SubExponential} there exists an absolute constant $C \geq 1$ (independent of $s, n, d$) such that
	\begin{align*}
	L_{n,\max}^{1/3}(s) &\leq C \max_{1 \leq k \leq d} \left(\frac{1}{n}\sum_{i=1}^n \mathrm{E}\left[X_{ik}^2\right]^{s/2}\right)^{1/s}\\
	& \leq C\max_{1 \leq k \leq d} \left(\frac{1}{n}\sum_{i=1}^n \mathrm{E}\left[X_{ik}^2\right]\right)^{1/s} \max_{1 \leq i \leq n}  \mathrm{E}[X_{ik}^2]^{(s-2)/(2s)},
	\end{align*}
	and, by H{\"o}lder's inequality,
	\begin{align*}
	\|\sigma_n\|_\infty = \max_{1 \leq k \leq d} \left(\frac{1}{n}\sum_{i=1}^n \mathrm{E}\left[X_{ik}^2\right]\right)^{1/2} \geq \max_{1 \leq k \leq d} \left(\frac{1}{n}\sum_{i=1}^n \mathrm{E}\left[X_{ik}^2\right]\right)^{1/s} \min_{1 \leq i \leq n} \mathrm{E}[X_{ik}^2]^{(s-2)/(2s)}.
	\end{align*}
	Combine the preceding two inequalities to obtain
	\begin{align}\label{eq:theorem-Gaussian-Approximation-5}
	\frac{(\log d)^{7/6}}{n^{1/6}}\frac{L_{n, \max}^{1/3}}{\|\sigma_n\|_\infty} \leq C \frac{(\log d)^{7/6}}{n^{1/6}} \kappa_n^{(s-2)/(2s)}
	\end{align}	
	Set $s = 6$, combine eq.~\eqref{eq:theorem-Gaussian-Approximation-4} and~\eqref{eq:theorem-Gaussian-Approximation-5} with Proposition~\ref{proposition:Berry-Esseen-Large-p} and conclude that
	\begin{align*}
	\varrho_{n,p} \lesssim \frac{\log^7 d}{n} + \left(\frac{\kappa_n^2 \log^7 d}{n}\right)^{1/6} \lesssim  \left(\frac{\kappa_n^2 \log^7 d}{n}\right)^{1/6},
	\end{align*}
	where the last inequality follows since without loss of generality we may assume that $C\kappa_n^2 \\ (\log d)^7/n < 1$ (otherwise the bound is trivially true).
\end{proof}
\begin{proof}[\textbf{Proof of Theorem~\ref{theorem:ConsistencyGPB}}]
	The result follows from the triangle inequality and Theorem~\ref{theorem:GaussianApprox} and Theorem~\ref{theorem:Gaussian-Comparison-KolmogorovDistance} as described in the main text.
\end{proof}
\begin{proof}[\textbf{Proof of Corollary~\ref{corollary:ConsistencyGPB-ApproxSparsity}}]
	Combine Theorem~\ref{theorem:ConsistencyGPB} and Lemma~\ref{lemma:BoundsThresholdEstimators}.	
\end{proof}
\begin{proof}[\textbf{Proof of Corollary~\ref{corollary:ConsistencyGPB-ApproxBandable}}]
	Combine Theorem~\ref{theorem:ConsistencyGPB} and Lemma~\ref{lemma:BoundsBandedEstimators}.	
\end{proof}

\subsection{Proofs for Section~\ref{sec:SimultaneousTesting}}
\begin{proof}[\textbf{Proof of Theorem~\ref{theorem:SizeAlphaTest}}]
	The proof is an adaptation of the proof of Theorem 3.1 in~\cite{chernozhukov2013GaussianApproxVec} to our setup.	Note that
	\begin{align}\label{eq:theorem:SizeAlphaTest-1-1}
	&\sup_{\alpha \in (0,1)}  \sup_{\mu \in \mathcal{H}_0} \left|\mathrm{P}_\mu\left(S_{n,p} +  \xi\leq c_{n,p}^*(\alpha)\right) - \alpha\right| \nonumber \\
	\begin{split}
	&\quad{}\leq\sup_{\alpha \in (0,1)}  \sup_{\mu \in \mathcal{H}_0}\left|\mathrm{P}_\mu\left(S_{n,p} + \xi \leq c_{n,p}^*(\alpha)\right) - \mathrm{P}_\mu\left(S_{n,p} + \xi \leq \tilde{c}_p(\alpha) \right) \right|\\
	&\quad{}\quad{}+ \sup_{\alpha \in (0,1)}  \sup_{\mu \in \mathcal{H}_0} \left|\mathrm{P}_\mu\left(S_{n,p} + \xi \leq \tilde{c}_p(\alpha) \right) - \mathrm{P}_\mu\left(\widetilde{S}_{n,p} \leq \tilde{c}_p(\alpha) \right) \right|.
	\end{split}
	\end{align}
	For $\delta > 0$ arbitrary the first term can be upper bounded by Lemma~\ref{lemma:ComparisonQuantiles} as
	\begin{align}\label{eq:theorem:SizeAlphaTest-2-1}
	&\sup_{\alpha \in (0,1)} \sup_{\mu \in \mathcal{H}_0}\mathrm{P}_\mu \Big(  \tilde{c}_p\big(\alpha - \pi_p(\delta) \big) < S_{n,p} + \xi \leq \tilde{c}_p\big(\alpha + \pi_p(\delta)\big)\Big)  + 2 \mathrm{P}_\mu\left(\Pi_p > \delta\right)\nonumber\\
	&\quad{}\leq \sup_{\alpha \in (0,1)} \sup_{\mu \in \mathcal{H}_0}\mathrm{P}_\mu \Big(\tilde{c}_p\big(\alpha - \pi_p(\delta)\big) < \widetilde{S}_p + \xi \leq \tilde{c}_p\big(\alpha + \pi_p(\delta)\big)\Big) \nonumber\\
	&\quad{}\quad{}+ 2\sup_{t \geq 0} \sup_{\mu \in \mathcal{H}_0} \left|\mathrm{P}_\mu\big(S_{n,p} \leq t\big) - \mathrm{P}_\mu\big(\widetilde{S}_p \leq t\big) \right| + 2 \mathrm{P}_\mu\left(\Pi_p > \delta\right) \nonumber\\
	\begin{split}
	&\quad{}\leq \sup_{\alpha \in (0,1)} \sup_{\mu \in \mathcal{H}_0} \left\{\mathrm{P}_\mu \Big(\tilde{c}_p\big(\alpha - \pi_p(\delta)\big) < \widetilde{S}_p + \xi \leq \tilde{c}_p\big(\alpha + \pi_p(\delta)\big)\Big) \right.\\
	&\quad{}\quad{}\quad{} \quad{} \quad{}\quad{}\quad{} \quad{} \left.- \mathrm{P}_\mu \Big(\tilde{c}_p\big(\alpha - \pi_p(\delta)\big) < \widetilde{S}_p\leq \tilde{c}_p\big(\alpha + \pi_p(\delta)\big)\Big) \right\}\\
	&\quad{}\quad{}+ 2 \pi_p(\delta) + 2\sup_{t \geq 0} \sup_{\mu \in \mathcal{H}_0} \left|\mathrm{P}_\mu\big(S_{n,p} \leq t\big) - \mathrm{P}_\mu\big(\widetilde{S}_p \leq t\big) \right| + 2\sup_{\mu \in \mathcal{H}_0} \mathrm{P}_\mu\left(\Pi_p > \delta\right),
	\end{split}
	\end{align}
	where the second inequality follows by definition of quantiles and because $\widetilde{S}_p$ has no point masses. Let $\eta > 0$ be arbitrary. The first term on the right hand side of eq.~\eqref{eq:theorem:SizeAlphaTest-2-1} can be bounded in the following way:
	\begin{align}\label{eq:theorem:SizeAlphaTest-3-1}
	&\sup_{\alpha \in (0,1)} \sup_{\mu \in \mathcal{H}_0} \left\{\mathrm{P}_\mu \Big(\tilde{c}_p\big(\alpha - \pi_p(\delta)\big) < \widetilde{S}_p + \xi \leq \tilde{c}_p\big(\alpha + \pi_p(\delta)\big)\Big) \right. \nonumber\\
	&\quad{}\quad{}\quad{} \left. - \mathrm{P}_\mu \Big(\tilde{c}_p\big(\alpha - \pi_p(\delta)\big) < \widetilde{S}_p\leq \tilde{c}_p\big(\alpha + \pi_p(\delta)\big)\Big) \right\} \nonumber\\
	&\quad{} \leq  \sup_{\alpha \in (0,1)} \sup_{\mu \in \mathcal{H}_0}\left|\mathrm{P}_\mu \left(\widetilde{S}_p + \xi \leq \tilde{c}_p\big(\alpha + \pi_p(\delta) \big)\right) - \mathrm{P}_\mu \left(\widetilde{S}_p \leq \tilde{c}_p\big(\alpha + \pi_p(\delta) \big)\right) \right|\nonumber \\
	&\quad{}\quad{} + \sup_{\alpha \in (0,1)} \sup_{\mu \in \mathcal{H}_0}\left|\mathrm{P}_\mu \left(\widetilde{S}_p \leq \tilde{c}_p\big(\alpha - \pi_p(\delta)\big) \right) - \mathrm{P}_\mu \left(\widetilde{S}_p +  \xi \leq \tilde{c}_p\big(\alpha - \pi_p(\delta) \big)\right) \right|\nonumber\\
	&\quad{} \lesssim \sup_{t \geq 0} \sup_{\mu \in \mathcal{H}_0}\left| \mathrm{P}_\mu \left(\widetilde{S}_p \leq t\right) - \mathrm{P}_\mu \left(\widetilde{S}_p +  \xi \leq t \right) \right|\nonumber\\
	&\quad{} \lesssim \sup_{\mu \in \mathcal{H}_0} \mathrm{P}_\mu \left(|\xi| > \eta \right) + \sup_{t \geq 0} \sup_{\mu \in \mathcal{H}_0} \mathrm{P}_\mu \left(t - \eta \leq \widetilde{S}_p  \leq t + \eta \right)\nonumber\\
	&\quad{} \lesssim  \sup_{\mu \in \mathcal{H}_0} \mathrm{P}_\mu \left(|\xi| > \eta \right) +  \eta \frac{\omega_p(d', r_\omega)}{\|\omega\|_p},
	\end{align}
	where the last inequality follows from Theorem~\ref{theorem:AntiConcentration-LpNorm}.
	
	We now bound the second term on the right hand side of eq.~\eqref{eq:theorem:SizeAlphaTest-1-1} by
	\begin{align}\label{eq:theorem:SizeAlphaTest-4-1}
	&\sup_{\alpha \in (0,1)}  \sup_{\mu \in \mathcal{H}_0} \left|\mathrm{P}_\mu\left(S_{n,p} + \xi \leq \tilde{c}_p(\alpha) \right) - \mathrm{P}_\mu\left(S_{n,p} \leq \tilde{c}_p(\alpha) \right) \right|  \nonumber\\
	&\quad{}+ \sup_{t \geq 0} \sup_{\mu \in \mathcal{H}_0}\left| \mathrm{P}_\mu \left(S_{n,p} \leq t\right) - \mathrm{P}_\mu \left(\widetilde{S}_p \leq t \right) \right|\nonumber\\
	&\quad{} \lesssim \sup_{\mu \in \mathcal{H}_0} \mathrm{P}_\mu \left(|\xi| > \eta \right) + \sup_{t \geq 0} \sup_{\mu \in \mathcal{H}_0} \mathrm{P}_\mu \left(t - \eta \leq S_{n,p}  \leq t + \eta \right) \nonumber\\
	&\quad{}\quad{} + \sup_{t \geq 0} \sup_{\mu \in \mathcal{H}_0}\left| \mathrm{P}_\mu \left(S_{n,p} \leq t\right) - \mathrm{P}_\mu \left(\widetilde{S}_p \leq t \right) \right|\nonumber\\
	&\quad{} \lesssim \sup_{\mu \in \mathcal{H}_0} \mathrm{P}_\mu \left(|\xi| > \eta \right) + \sup_{t \geq 0} \sup_{\mu \in \mathcal{H}_0} \mathrm{P}_\mu \left(t - \eta \leq \widetilde{S}_p  \leq t + \eta \right) \nonumber\\
	&\quad{}\quad{} + \sup_{t \geq 0} \sup_{\mu \in \mathcal{H}_0}\left| \mathrm{P}_\mu \left(S_{n,p} \leq t\right) - \mathrm{P}_\mu \left(\widetilde{S}_p \leq t \right) \right|\nonumber\\
	&\quad{}\lesssim  \sup_{\mu \in \mathcal{H}_0} \mathrm{P}_\mu \left(|\xi| > \eta \right) +  \eta \frac{\omega_p(d', r_\omega)}{\|\omega\|_p} + \sup_{t \geq 0} \sup_{\mu \in \mathcal{H}_0}\left| \mathrm{P}_\mu \left(S_{n,p} \leq t\right) - \mathrm{P}_\mu \left(\widetilde{S}_p \leq t \right) \right|,
	\end{align}	
	where the third inequality follows from Theorem~\ref{theorem:AntiConcentration-LpNorm}.
	
	Combine eq.~\eqref{eq:theorem:SizeAlphaTest-1-1}--\eqref{eq:theorem:SizeAlphaTest-4-1} to obtain
	\begin{align*}
	&\sup_{\alpha \in (0,1)}  \sup_{\mu \in \mathcal{H}_0} \left|\mathrm{P}_\mu\left(S_{n,p} +  \xi\leq c_{n,p}^*(\alpha)\right) - \alpha\right|\nonumber\\
	&\quad{} \lesssim  \sup_{t \geq 0} \sup_{\mu \in \mathcal{H}_0}\left| \mathrm{P}_\mu \left(S_{n,p} \leq t\right) - \mathrm{P}_\mu \left(\widetilde{S}_p \leq t \right) \right| + \inf_{\delta > 0}\left\{ \pi_p(\delta) + \sup_{\mu \in \mathcal{H}_0} \mathrm{P}_\mu\left(\Pi_p > \delta\right) \right\} \nonumber\\
	&\quad{}\quad{}+ \inf_{\eta > 0} \left\{ \eta \frac{\omega_p(d', r_\omega)}{\|\omega\|_p} +  \sup_{\mu \in \mathcal{H}_0} \mathrm{P}_\mu \left(|\xi| > \eta \right)\right\}.
	\end{align*}
	To complete the proof bound the first term on the right hand side by Theorem~\ref{theorem:GaussianApprox}.
\end{proof}
\begin{proof}[\textbf{Proof of Theorem~\ref{theorem:ConsistencyLocalAlternatives}}]
	We prove a slightly sharper result than the one given in the main text. Let $\mathcal{L} \subseteq \mathbb{R}_+$ be arbitrary and define the collection of alternatives
	\begin{align}\label{eq:seubsec:ConsistencyTesting-1-1}
	\mathcal{A}_p(\mathcal{L}) := \left\{ \left(\mu_n\right)_{n \in \mathbb{N}}, \mu_n \in \mathbb{R}^{d_n} : \lim_{n \rightarrow \infty} \frac{ \mathrm{E}\|\Omega^{1/2}Z\|_p \vee \sqrt{\mathrm{Var}\|\Omega^{1/2}Z\|_p}}{ \sqrt{n}\|M\mu_n - m_0\|_p} \in \mathcal{L}\right\},
	\end{align}
	where $Z \sim N(0, I_{d'})$. Note that $\mathcal{A}_p = \mathcal{A}_p(\{0\})$ and $\mathcal{Z}_p \subset \mathcal{A}_p((1, \infty])$.
	
	\textbf{Prove of Case (i).} Recall $\pi_p(\cdot)$ and $\Pi_p$ from Lemma~\ref{lemma:ComparisonQuantiles}. Let $\delta_\alpha > 0$ be such that $1/ (\alpha - \pi_p(\delta_\alpha) \leq 2/\alpha$. Fix a sequence $(\mu_n)_{n \in \mathbb{N}} \in \mathcal{A}_p(\{0\})$. Then, by Lemma~\ref{lemma:ComparisonQuantiles} and Assumption~\ref{assumption:SufficientConditions},
	\begin{align}\label{eq:theorem:ConsistencyLocalAlternatives-1}
	&\mathrm{P}_{\mu_n}\left(S_{n,p} > c_{n,p}^*(1-\alpha)\right) \nonumber\\
	&\quad{} \geq \mathrm{P}_{\mu_n}\left(S_{n,p} > c_{n,p}^*(1-\alpha), \:  c_{n,p}^*(1-\alpha) \leq \tilde{c}_p\big(\pi_p(\delta_\alpha) + 1- \alpha\big)\right)  \nonumber\\
	&\quad{} \geq \mathrm{P}_{\mu_n}\left(S_{n,p} > \tilde{c}_p\big(\pi_p(\delta_\alpha) + 1- \alpha\big) \right)  + \mathrm{P} \left(c_{n,p}^*(1-\alpha) \leq \tilde{c}_p\big(\pi_p(\delta_\alpha) + 1- \alpha\big)\right) -1\nonumber\\
	&\quad{} \geq \mathrm{P}_{\mu_n}\left(S_{n,p} > \tilde{c}_p\big(\pi_p(\delta_\alpha) + 1- \alpha\big) \right) - \mathrm{P} \big( \Pi_p >\delta_\alpha\big) \nonumber\\
	&\quad{} \geq \mathrm{P}_{\mu_n}\left(S_{n,p} > \tilde{c}_p\big(\pi_p(\delta_\alpha) + 1- \alpha\big) \right)  + o(1).	
	\end{align}
	We now lower bound the first factor on the far right hand side in above display. By the reverse triangle inequality and Lemma~\ref{lemma:UpperBoundQuantiles},
	\begin{align}\label{eq:theorem:ConsistencyLocalAlternatives-2}
	&\mathrm{P}_{\mu_n}\left(S_{n,p} > \tilde{c}_p\big(\pi_p(\delta_\alpha) + 1- \alpha\big) \right) \nonumber\\
	&\quad{} = \mathrm{P}_{\mu_n}\left(\left\|\frac{1}{\sqrt{n}}\sum_{i=1}^n M(X_i - \mu_n) + \sqrt{n}(M \mu_n - m_0)  \right\|_p > \tilde{c}_p\big(\pi_p(\delta_\alpha) + 1- \alpha\big)\right) \nonumber\\
	&\quad{} \geq \mathrm{P}_{\mu_n}\left(\left\|\frac{1}{\sqrt{n}}\sum_{i=1}^n M(X_i - \mu_n) \right\|_p < \sqrt{n}\|M \mu_n- m_0\|_p - \tilde{c}_p\big(\pi_p(\delta_\alpha) + 1- \alpha\big)\right)\nonumber\\
	&\quad{} \geq \mathrm{P}_0\left(\widetilde{S}_p < \sqrt{n}\|M \mu_n- m_0\|_p -  \tilde{c}_p\big(\pi_p(\delta_\alpha) + 1- \alpha\big) \right) - \sup_{t \geq 0} \left|\mathrm{P}_0\left(S_{n,p}\leq t \right) - \mathrm{P}\left(\widetilde{S}_p \leq t \right)  \right|\nonumber\\
	\begin{split}
	&\quad{} \geq \mathrm{P}_0\left(\widetilde{S}_p < \sqrt{n}\|M \mu_n- m_0\|_p - \mathrm{E}[\widetilde{S}_p]  - \sqrt{1/\big(\alpha- \pi_p(\delta_\alpha)\big)\mathrm{Var}[\widetilde{S}_p]}\right)\\
	&\quad{}\quad{}\quad{} - \sup_{t \geq 0} \left|\mathrm{P}_0\left(S_{n,p} \leq t \right) - \mathrm{P}\left(\widetilde{S}_p \leq t \right)  \right|.
	\end{split}
	\end{align}
	By Theorem~\ref{theorem:GaussianApprox} and Assumption~\ref{assumption:SufficientConditions} the second term in the last line is of order $o(1)$. By Markov's inequality the first term can be bounded by
	\begin{align}\label{eq:theorem:ConsistencyLocalAlternatives-3}
	&\inf_{\mu_n\in \mathcal{A}_p}\mathrm{P}_0\left(\widetilde{S}_p < \sqrt{n}\|M \mu_n- m_0\|_p - \mathrm{E}[\widetilde{S}_p]  - \sqrt{1/\big(\alpha- \pi_p(\delta_\alpha)\big)\mathrm{Var}[\widetilde{S}_p]} \right) \nonumber\\
	&\quad{}\geq 1 - \frac{2\mathrm{E}[\widetilde{S}_p] + \sqrt{(2/\alpha)\mathrm{Var}[\widetilde{S}_p]}}{\sqrt{n}\|M \mu_n- m_0\|_p} \nonumber\\
	&\quad{}= 1 + o(1),
	\end{align}
	where we have used that $(\mu_n)_{n \in \mathbb{N}} \in \mathcal{A}_p(\{0\})$.
	
	To conclude the proof combine eq.~\eqref{eq:theorem:ConsistencyLocalAlternatives-1}--\eqref{eq:theorem:ConsistencyLocalAlternatives-3}.
	
	\textbf{Prove of Case (ii).} Recall $\pi_p(\cdot)$ and $\Pi_p$ from Lemma~\ref{lemma:ComparisonQuantiles}. Let $\delta_\alpha > 0$ be such that $1/ (\alpha - \pi_p(\delta_\alpha) \leq 2/\alpha$.
	Let $(\mu_n)_{n \in \mathbb{N}} \in \mathcal{A}_p((1, \infty])$, and compute
	\begin{align}\label{eq:theorem:ConsistencyLocalAlternatives-4}
	&\mathrm{P}_{\mu_n}\left(S_{n,p} > c_{n,p}^*(1-\alpha)\right) \nonumber\\
	&\quad{}\leq \mathrm{P}_{\mu_n}\left(\left\|\frac{1}{\sqrt{n}}\sum_{i=1}^n M(X_i - \mu_n) \right\|_p + \sqrt{n} \left\| M \mu_n - m_0\right\|_p > c_{n,p}^*(1-\alpha)\right) \nonumber\\
	\begin{split}
	&\quad{}\leq \mathrm{P}_0\left(\widetilde{S}_p + \sqrt{n} \left\| M \mu_n - m_0\right\|_p > \tilde{c}_p(1-\alpha)\right)\\
	&\quad{} + \left\{\mathrm{P}_0\left(S_{n,p} +  \sqrt{n} \left\| M \mu_n - m_0\right\|_p\leq c_{n,p}^*(1-\alpha)\right) \right.\\
	&\quad{}\quad{}\quad{}\quad{} \left. - \mathrm{P}_{0}\left(S_{n,p} + \sqrt{n} \left\| M \mu_n - m_0\right\|_p \leq  \tilde{c}_p(1-\alpha)\right) \right\}.
	\end{split}
	\end{align}
		
	The second term can be bounded as the first term in eq.~\eqref{eq:theorem:SizeAlphaTest-1-1}: Set $\xi = \sqrt{n} \left\| M \mu_n - m_0\right\|_p$. Then, as in eq.~\eqref{eq:theorem:SizeAlphaTest-2-1} and by the equivariance of the quantile function,
	\begin{align}\label{eq:theorem:ConsistencyLocalAlternatives-5}
	&\mathrm{P}_0\left(S_{n,p} + \xi \leq c_{n,p}^*(1-\alpha)\right) - \mathrm{P}_{0}\left(S_{n,p} + \xi \leq  \tilde{c}_p(1-\alpha)\right) \nonumber\\
	&\quad{} \leq \mathrm{P}_0 \Big(\tilde{c}_p\big(1- \alpha - \pi_p(\delta_\alpha)\big) < \widetilde{S}_p + \xi \leq \tilde{c}_p\big(1- \alpha + \pi_p(\delta_\alpha)\big)\Big) \nonumber\\
	&\quad{}\quad{} - \mathrm{P}_0 \Big(\tilde{c}_p\big(1-\alpha - \pi_p(\delta_\alpha)\big) + \xi < \widetilde{S}_p + \xi\leq \tilde{c}_p\big(1-\alpha + \pi_p(\delta_\alpha)\big) + \xi\Big)\nonumber\\
	&\quad{}\quad{}+ 2 \pi_p(\delta_\alpha) + 2\sup_{t \geq 0} \left|\mathrm{P}_0\big(S_{n,p} \leq t\big) - \mathrm{P}_0\big(\widetilde{S}_p \leq t\big) \right| + 2 \mathrm{P}_0\left(\Pi_p > \delta_\alpha\right) \nonumber\\
	&\quad{}\leq  2\pi_p(\delta_\alpha) + 2\sup_{t \geq 0} \left|\mathrm{P}_0\big(S_{n,p} \leq t\big) - \mathrm{P}_0\big(\widetilde{S}_p \leq t\big) \right| + 2 \mathrm{P}_0\left(\Pi_p > \delta_\alpha\right) + o(1),
	\end{align}
	where we have used that
	\begin{align*}
	&\mathrm{P}_0 \Big(\tilde{c}_p\big(1- \alpha - \pi_p(\delta_\alpha)\big) < \widetilde{S}_p + \xi \leq \tilde{c}_p\big(1- \alpha + \pi_p(\delta_\alpha)\big)\Big) \nonumber\\
	&\quad{} - \mathrm{P}_0 \Big(\tilde{c}_p\big(1-\alpha - \pi_p(\delta_\alpha)\big) + \xi < \widetilde{S}_p + \xi\leq \tilde{c}_p\big(1-\alpha + \pi_p(\delta_\alpha)\big) + \xi\Big) \nonumber\\
	& =  - \mathrm{P}_0 \Big(\tilde{c}_p\big(1- \alpha + \pi_p(\delta_\alpha)\big) < \widetilde{S}_p + \xi \leq \tilde{c}_p\big(1- \alpha + \pi_p(\delta_\alpha)\big) + \xi\Big) \nonumber\\
	&\quad{}+ \mathrm{P}_0 \Big(\tilde{c}_p\big(1-\alpha - \pi_p(\delta_\alpha)\big) < \widetilde{S}_p + \xi\leq \tilde{c}_p\big(1-\alpha - \pi_p(\delta_\alpha)\big) + \xi \Big)  \nonumber\\
	&= o(1),
	\end{align*}
	since by Assumption~\ref{assumption:SufficientConditions} $\pi_p(\delta_\alpha) = o(1)$.
	
	By another application of Assumption~\ref{assumption:SufficientConditions} and Theorem~\ref{theorem:GaussianApprox} we conclude that the remaining terms in eq.~\eqref{eq:theorem:ConsistencyLocalAlternatives-5} are negligible as well.
	
	Since $\widetilde{S}_{p}$ has no point mass, the first term on the far right hand side of eq.~\eqref{eq:theorem:ConsistencyLocalAlternatives-4} is strictly less than 1 whenever
	\begin{align}\label{eq:theorem:ConsistencyLocalAlternatives-6}
	&\sqrt{n} \left\| M \mu_n - m_0\right\|_p < \mathrm{E}[\widetilde{S}_p] - \sqrt{\mathrm{Var}[\widetilde{S}_p]} \leq \tilde{c}_p(1-\alpha),
	\end{align}
	where the second inequality follows from Lemma~\ref{lemma:UpperBoundQuantiles}. This inequality holds by definition of the set $\mathcal{A}_p((1, \infty])$
	
	To conclude the proof combine eq.~\eqref{eq:theorem:ConsistencyLocalAlternatives-3}--\eqref{eq:theorem:ConsistencyLocalAlternatives-6}.
\end{proof}

\subsection{Proofs for Appendix~\ref{apx:MainTheory}}
\subsubsection{Proofs for Appendix~\ref{subsec:BerryEsseen}}
\begin{proof}[\textbf{Proof of Theorem~\ref{theorem:BerryEsseen-Even-p}}]
	
	\noindent
	
	\textbf{Proof of Case (i).}
	
	\textbf{Step 1. Fundamental smoothing inequality.} Let $Y=\{Y_i\}_{i=1}^n$ be an independent copy of $Z=\{Z_i\}_{i=1}^n$ and define
	\begin{align*}
	W_n(s) := \sum_{i=1}^n \left(\sqrt{\frac{s}{n}}X_i + \sqrt{\frac{1-s}{n}} Y_i\right), \hspace{10pt} s \in [0,1].
	\end{align*}
	
	Consider the family of sets $\mathcal{I} = \{A \subseteq \mathbb{R} : A = [0,t], t \geq 0 \}$.  Let $p \in [1, \infty)$ be arbitrary. Define $p_+ = 2\lceil \frac{p}{2}\rceil$ to be the smallest even integer larger than (or equal to) $p$. By Lemma~\ref{lemma:SmoothIndicator-C-Infty-Function-Even-p} for $A \in \mathcal{I}$, we have
	\begin{align*}
	\mathrm{P}\left(\|W_n(s)\|_p  \in A \right) - \mathrm{P}\left(\|S_n^Z\|_p   \in A^{12\kappa_p + 3\delta}  \right) \leq \mathrm{E}\left[h_{p_+, d, \beta, A^{3\kappa_{p_+}}}\big(W_n(s)\big)- h_{p_+, d, \beta, A^{3\kappa_{p_+}}}(S_n^Z)\right].
	\end{align*}
	Re-arrange the terms in above inequality and take the supremum over $A \in \mathcal{I}$ to obtain
	\begin{align}\label{eq:theorem:BerryEsseen-Even-p-1}
	\begin{split}
	&\sup_{A \in \mathcal{I}} \Big(\mathrm{P}\left(\|W_n(s)\|_p  \in A \right) - \mathrm{P}\left(\|S_n^Z\|_p   \in A\right) \Big)\\
	&\quad{}\leq \sup_{A \in \mathcal{I}}\mathrm{P}\left(\|S_n^Z\|_p   \in A^{12\kappa_p + 3\delta} \setminus A \right) + \sup_{A \in \mathcal{I}} \Big|\mathrm{E}\left[h_{p_+, d, \beta,\delta, A}\big(W_n(s)\big)- h_{p_+, d, \beta,\delta, A}(S_n^Z)\right]\Big|,
	\end{split}
	\end{align}
	By Lemma~\ref{lemma:SmoothIndicator-C-Infty-Function-Even-p} we also have for $A \in \mathcal{I}$,
	\begin{align*}
	&\mathrm{P}\left(\|S_n^Z\|_p  \in A^{-(12\kappa_p + 3\delta)} \right) - \mathrm{P}\left(\|W_n(s)\|_p   \in A \right) \\
	&\quad{}\leq \mathrm{E}\left[h_{p_+, d, \beta, A^{-(12\kappa_{p_+}+ 3\delta)}}(S_n^Z)- h_{p_+, d, \beta, A^{-(12\kappa_{p_+} + 3\delta)}}\big(W_n(v)\big)\right].
	\end{align*}
	Observe that $\sup_{A \in \mathcal{I}}\mathrm{P}\left( \|S_n^Z\|_p \in A \setminus A^{-(12\kappa_p + 3\delta)}\right) \leq \sup_{A \in \mathcal{I}}\mathrm{P}\left( \|S_n^Z\|_p \in A^{12\kappa_p + 3\delta}\setminus A\right)$. Together with the preceding inequality this yields
	\begin{align}\label{eq:theorem:BerryEsseen-Even-p-2}
	\begin{split}
	&\sup_{A \in \mathcal{I}} \Big( \mathrm{P}\left(\|S_n^Z\|_p  \in A \right) - \mathrm{P}\left(\|W_n(s)\|_p   \in A \right) \Big)\\
	&\quad{}\leq \sup_{A \in \mathcal{I}}\mathrm{P}\left(\|S_n^Z\|_p   \in A^{12\kappa_p + 3\delta} \setminus A \right) + \sup_{A \in \mathcal{I}} \Big|\mathrm{E}\left[h_{p_+, d, \beta, \delta,A}\big(W_n(s)\big)- h_{p_+, d, \beta,\delta, A}(S_n^Z)\right]\Big|,
	\end{split}
	\end{align}
	Combine eq.~\eqref{eq:theorem:BerryEsseen-Even-p-1} and eq.~\eqref{eq:theorem:BerryEsseen-Even-p-2} to obtain
	\begin{align*}
	&\sup_{s \in [0,1]}\sup_{A \in \mathcal{I}} \Big| \mathrm{P}\left(\|S_n^Z\|_p  \in A \right) - \mathrm{P}\left(\|W_n(s)\|_p   \in A \right) \Big|\\
	&\quad{}\leq \sup_{A \in \mathcal{I}}\mathrm{P}\left(\|S_n^Z\|_p \in A^{12\kappa_p + 3\delta} \setminus A \right) + \sup_{s \in [0,1]}\sup_{A \in \mathcal{I}} \Big|\mathrm{E}\left[h_{p_+, d, \beta, \delta,A}\big(W_n(s)\big)- h_{p_+, d, \beta,\delta, A}(S_n^Z)\right]\Big|
	\end{align*}
	Note that above inequality holds also for $p = p_+$. Thus, we have the following fundamental smoothing inequality
	\begin{align}\label{eq:theorem:BerryEsseen-Even-p-3}
	\begin{split}
	&\sup_{s \in [0,1]}\sup_{A \in \mathcal{I}} \sup_{r \in \{p, p_+\}}\Big| \mathrm{P}\left(\|S_n^Z\|_r \in A \right) - \mathrm{P}\left(\|W_n(s)\|_r   \in A \right) \Big| \\
	&\quad{}\quad{}\quad{}\quad{}\leq \sup_{A \in \mathcal{I}} \sup_{r \in \{p, p_+\}}\mathrm{P}\left(\|S_n^Z\|_r \in A^{12\kappa_r + 3\delta} \setminus A \right) \\
	&\quad{}\quad{}\quad{}\quad{}\quad{}\quad{}\quad{}\quad{} + \sup_{s \in [0,1]}\sup_{A \in \mathcal{I}} \Big|\mathrm{E}\left[h_{p_+, d, \beta, \delta,A}\big(W_n(s)\big)- h_{p_+, d, \beta,\delta, A}(S_n^Z)\right]\Big|
	\end{split}
	\end{align}
	We now bound the second term on the right hand side of eq.~\eqref{eq:theorem:BerryEsseen-Even-p-3}.
	
	\textbf{Step 2. Slepian-Stein interpolation.} Define the Slepian interpolant as
	\begin{align*}
	V(t;s) := \sum_{i=1}^n V_i(t;s), \hspace{10pt} \mathrm{where} \hspace{10pt} V_i(t;s) := \sqrt{\frac{st}{n}}X_i + \sqrt{\frac{(1-s)t}{n}} Y_i + \sqrt{\frac{1 - t}{n}}Z_i , \hspace{10pt} s,t \in [0,1],
	\end{align*}
	the Stein leave-one-out term as
	\begin{align*}
	V^{(i)}(t;s) := V(t;s) - V_i(t;s), \hspace{10pt} i=1, \ldots, n,
	\end{align*}
	and denote the derivative of the $i$th summand $V_i(t;s)$ with respect to $t$ by
	\begin{align*}
	\dot{V}_i(t;s) :=\frac{d}{dt}V_i(t;s) = \frac{1}{2} \left[\frac{1}{\sqrt{t}}\left(\sqrt{\frac{s}{n}}X_i + \sqrt{\frac{1-s}{n}} Y_i\right)  - \frac{1}{\sqrt{1-t}}\frac{1}{\sqrt{n}} Z_i \right].
	\end{align*}
	Since $V(0;s) = S_n^Z$ and $V(1;s) = W_n(s)$, by expressing the difference as integration of the derivative function, we have
	\begin{align}\label{eq:theorem:BerryEsseen-Even-p-4}
	\begin{split}
	&\mathrm{E}\left[h_{p_+, d, \beta, \delta,A}\big(W_n(s)\big)- h_{p_+, d, \beta,\delta, A}(S_n^Z) \right]\\
	&\hspace{40pt}= \sum_{i=1}^n\sum_{|\alpha|=1} \int_0^1 \mathrm{E}\left[\dot{V}_i^\alpha(t;s)\big(D^\alpha h_{p_+, d, \beta, \delta,A} \big) \big(V(t;s)\big)\right] dt.
	\end{split}
	\end{align}
	For brevity of notation, we now drop the subscripts $p_+,d,\beta, \delta, A$ and write $h$ instead of $h_{p_+, d, \beta,\delta, A}$. We also write $X_{ni}$, $Y_{ni}$, and $Z_{ni}$ instead of $\frac{1}{\sqrt{n}}X_i$, $\frac{1}{\sqrt{n}}Y_i$, and $\frac{1}{\sqrt{n}}Z_i$, respectively. In above display, expanding the summands over $i=1, \ldots, n$ via a first-order Taylor expansion around $V^{(i)}(t;s)$ in direction $V_i(t;s)$ yields, for all $s \in [0,1]$,
	\begin{align}\label{eq:theorem:BerryEsseen-Even-p-5}
	&\mathrm{E}\left[h\big(W_n(s)\big)- h(S_n^Z) \right] \nonumber\\
	&\quad{} =\sum_{i=1}^n \sum_{|\alpha|=1}\int_0^1  \mathrm{E}\left[\dot{V}_i^\alpha (t;s)\big(D^\alpha h\big) \big(V^{(i)}(t;s)\big) \right] dt \nonumber\\
	&\quad{} \quad{} +\sum_{i=1}^n\sum_{|\alpha'|=1} \sum_{|\alpha|=1} \int_0^1 \mathrm{E}\left[V_i^{\alpha'}(t;s)\dot{V}_i^\alpha(t;s)\big(D^{\alpha + \alpha'} h \big) \big(V^{(i)}(t;s)\big) \right] dt\nonumber \\
	&\quad{} \quad{} +\sum_{i=1}^n\sum_{|\alpha'|=2} \sum_{|\alpha|=1} \int_0^1 \int_0^1 (1-u) \mathrm{E}\left[V_i^{\alpha'}(t;s)\dot{V}_i^\alpha(t;s)\big(D^{\alpha + \alpha'} h\big) \big(V^{(i)}(t;s) + uV_i(t;s)\big) \right] dt du\nonumber\\
	&\quad{} = \mathrm{\mathbf{I}} + \mathrm{\mathbf{II}} + \mathrm{\mathbf{III}} .
	\end{align}
	It is standard to verify that $\mathrm{\mathbf{I}} = \mathrm{\mathbf{II}} = 0$ (because $\mathbb{E}[\dot{V}_i (t;s)] =0$ and $\mathbb{E}[V_i^{\alpha'} (t;s) \dot{V}_i^\alpha(t;s)] = 0$ for $|\alpha'| = |\alpha|=1$, and $\dot{V}_i(t;s)$ and $V_i(t;s)$ are independent of $V^{(i)}(t;s)$; see p. 2327 in~\cite{chernozhukov2017CLTHighDim}). Thus, we only need to bound the third term. Let $\xi > 0$, $\tau \geq 1$, set $\chi_i = \mathbf{1}\{\|X_{ni}\|_{\tau p_+} \vee \|Y_{ni}\|_{\tau p_+} \vee \|Z_{ni}\|_{\tau p_+}\leq \xi\}$, and compute
	\begin{align}\label{eq:theorem:BerryEsseen-Even-p-6}
	&\left|\mathrm{\mathbf{III}} \right| =\nonumber\\
	& = \sum_{i=1}^n\sum_{|\alpha'|=2} \sum_{|\alpha|=1} \int_0^1 \int_0^1 (1-u) \mathrm{E}\left[\chi_i \big|V_i^{\alpha'}(t;s)\dot{V}_i^\alpha(t;s)\big(D^{\alpha + \alpha'} h\big) \big(V^{(i)}(t;s) + uV_i(t;s)\big)\big| \right] dt du \nonumber\\
	&+ \sum_{i=1}^n\sum_{|\alpha'|=2} \sum_{|\alpha|=1} \int_0^1 \int_0^1 (1-u) \mathrm{E}\left[(1-\chi_i)\big|V_i^{\alpha'}(t;s)\dot{V}_i^\alpha(t;s)\big(D^{\alpha + \alpha'} h\big) \big(V^{(i)}(t;s) + uV_i(t;s)\big)\big| \right] dt du\nonumber\\
	&= \left|\mathrm{\mathbf{III}}_1\right| + \left|\mathrm{\mathbf{III}}_2\right|.
	\end{align}
	
	\textbf{Step 3. Bound on $\mathrm{\mathbf{III}}_2$.}
	Set $q = \frac{\tau p_+}{ \tau p_+ -1}$ and $B = \left\|\left(\sum_{|\alpha| = 3} \left| D^\alpha h\right|^q\right)^{1/q}\right\|_\infty$. By repeated applications of H{\"o}lder's inequality,
	\begin{align}\label{eq:theorem:BerryEsseen-Even-p-7}
	|\mathrm{\mathbf{III}_2}| &\leq  \mathrm{E}\left[\sum_{i=1}^n\int_0^1 \int_0^1 (1-u)(1-\chi_i)\left(\sum_{|\alpha'|=2} \sum_{|\alpha|=1} |V_i^{\alpha'}(t;s)\dot{V}_i^\alpha(t;s)|^{\tau p_+}\right)^{1/(\tau p_+)} \right. \nonumber\\
	&\quad{} \times \left.  \left(\sum_{|\alpha'|=2} \sum_{|\alpha|=1} \Big|\big(D^{\alpha + \alpha'} h\big) \big(V^{(i)}(t;s) + uV_i(t;s)\big)\Big|^q \right)^{1/q}  dt du\right] \nonumber\\
	&\leq B \mathrm{E}\left[\sum_{i=1}^n\int_0^1 (1-\chi_i)\left(\sum_{|\alpha'|=2} \sum_{|\alpha|=1} |V_i^{\alpha'}(t;s)\dot{V}_i^\alpha(t;s)|^{\tau p_+}\right)^{1/(\tau p_+)}dt  \right] \nonumber\\
	&\leq B \left(\int_0^1 \frac{dt}{\sqrt{t} \wedge \sqrt{1-t}} \right)
	\mathrm{E}\left[\sum_{i=1}^n (1-\chi_i) \left(\sum_{|\alpha'|=2} \sum_{|\alpha|=1} 9^{\tau p_+}\Big|\big( |X_{ni}|^{\alpha'} \vee |Y_{ni}|^{\alpha'} \vee |Z_{ni}|^{\alpha'}\big) \right.\right. \nonumber\\
	&\left.\left. \phantom{B \left(\int_0^1 \frac{dt}{\sqrt{t} \wedge \sqrt{1-t}} \right)
	\mathrm{E}\left[ (1-\chi_i) \left(\sum_{|\alpha'|=2} \right.\right.}\times \big( |X_{ni}|^{\alpha} + |Y_{ni}|^{\alpha}+ |Z_{ni}|^{\alpha}\big)\Big|^{\tau p_+} \right)^{1/ \tau p_+}  \right] \nonumber\\
	&\lesssim B \mathrm{E}\left[\sum_{i=1}^n (1-\chi_i)\left(\sum_{|\alpha'|=2} \sum_{|\alpha|=1} \Big| |X_{ni}|^{\alpha + \alpha'} \vee |Y_{ni}|^{\alpha + \alpha'} \vee |Z_{ni}|^{\alpha + \alpha'}\Big|^{\tau p_+}\right)^{1/(\tau p_+)}\right]
	\end{align}
	We now apply inequality (20) of Lemma B.1 in~\cite{chernozhukov2017CLTHighDim} to the far right hand side of~\eqref{eq:theorem:BerryEsseen-Even-p-7} and obtain
	\begin{align}
	&|\mathrm{\mathbf{III}_2}| \nonumber\\
	&\lesssim  B \mathrm{E}\left[\sum_{i=1}^n \left(\sum_{|\alpha|=3}  \Big| |X_{ni}|^{\alpha}\mathbf{1}\{\|X_{ni}\|_{\tau p_+} > \xi\}\Big|^{\tau p_+}\right)^{1/(\tau p_+)}\right. \nonumber\\
	&\quad{}\left. +  \left(\sum_{|\alpha|=3}  \Big| |Y_{ni}|^{\alpha}\mathbf{1}\{\|Y_{ni}\|_{\tau p_+} > \xi\}\Big|^{\tau p_+}\right)^{1/(\tau p_+)}\hspace{-20pt}+ \left(\sum_{|\alpha|=3}  \Big| |Z_{ni}|^{\alpha}\mathbf{1}\{\|Z_{ni}\|_{\tau p_+} > \xi\}\Big|^{\tau p_+}\right)^{1/(\tau p_+)}\right] \nonumber\\
	&\lesssim  B \mathrm{E}\left[\sum_{i=1}^n \|X_{ni}\|_{\tau p_+}^3\mathbf{1}\{\|X_{ni}\|_{\tau p_+} > \xi\} + \|Z_{ni}\|_{\tau p_+}^31\{\|Z_{ni}\|_{\tau p_+} > \xi\}\right],	
	\end{align}
	where the last inequality holds because $Y_{ni} \overset{d}{=} Z_{ni}$ for all $i=1, \ldots, n$.
	
	\textbf{Step 4. Bound on $\mathrm{\mathbf{III}}_1$.} Recall from Lemma~\ref{lemma:SmoothIndicator-C-Infty-Function-Even-p} that $h \equiv h_{p_+, d, \beta, \delta, A}$ is non-constant on the set $\big\{z \in \mathbb{R}^d: M_{p_+,\kappa_{p_+}}(z) \in A^{3\delta} \setminus A\big\}$ only. By construction of $M_{p_+,\kappa_{p_+}}$ it holds that $\big\{z \in \mathbb{R}^d: \|z\|_{p_+} \in A^{3\delta } \setminus A^{-\kappa_{p_+}}  \big\} \supseteq \big\{z \in \mathbb{R}^d: M_{p_+,\kappa_{p_+}}(z) \in A^{3\delta} \setminus A\big\}$. Thus, for $\varphi(x) = \mathbf{1}\big\{ \|x\|_{p_+}  \in A^{3\delta} \setminus A^{-\kappa_{p_+}} \big\}$ we have
	\begin{align}\label{eq:stability-1}
	\begin{split}
	&\chi_i\sum_{|\alpha'|=2} \sum_{|\alpha|=1} \int_0^1 (1-u) \big(D^{\alpha + \alpha'} h\big) \big(V^{(i)}(t;s) + uV_i(t;s)\big) du\\
	&\quad{}= \chi_i \varphi\big(V^{(i)}(t;s)+ uV_i(t;s)\big)\sum_{|\alpha'|=2} \sum_{|\alpha|=1} \int_0^1 (1-u) \big(D^{\alpha + \alpha'} h\big) \big(V^{(i)}(t;s) + uV_i(t;s)\big) du.
	\end{split}
	\end{align}
	Set $\gamma_{p_+} = d^{1/p_+-1/(\tau p_+)} \xi$. Note that for all $t, s, u \in [0,1]$,
	\begin{align*}
	\|V^{(i)}(t;s)\|_{p_+}\chi_i &\leq \|V^{(i)}(t;s) + u V_i(t; s)\|_{p_+}\chi_i + \|V_i(t;s)\|_{p_+}\chi_i \nonumber\\
	&\leq \|V^{(i)}(t;s) + u V_i(t; s)\|_{p_+}\chi_i + 3 \gamma_{p_+},
	\end{align*}
	and, similarly,
	\begin{align*}
	\|V^{(i)}(t;s)\|_{p_+}\chi_i \geq \|V^{(i)}(t;s) + u V_i(t; s)\|_{p_+}\chi_i - 3 \gamma_{p_+}.
	\end{align*}
	Thus, for all $t, s, u \in [0,1]$,
	\begin{align}\label{eq:stability-3}
	\left| \|V^{(i)}(t;s) + u V_i(t; s)\|_{p_+} - \|V^{(i)}(t;s)\|_{p_+}\right|\chi_i  \leq 3 \gamma_{p_+}.
	\end{align}
	Define $\phi(x) = \mathbf{1}\big\{ \|x\|_{p_+}  \in A^{3\delta + 3\gamma_{p_+}} \setminus A^{-(\kappa_{p_+} + 3\gamma_{p_+})} \big\}$. Now, eq.~\eqref{eq:stability-1} and~\eqref{eq:stability-3} imply
	\begin{align}\label{eq:stability-4}
	\begin{split}
	&\chi_i\sum_{|\alpha'|=2} \sum_{|\alpha|=1} \int_0^1 (1-u) \big(D^{\alpha + \alpha'} h\big) \big(V^{(i)}(t;s) + uV_i(t;s)\big) du\\
	&\quad{}\leq \chi_i \phi\big(V^{(i)}(t;s)\big)\sum_{|\alpha'|=2} \sum_{|\alpha|=1} \int_0^1 (1-u) \big(D^{\alpha + \alpha'} h\big) \big(V^{(i)}(t;s) + uV_i(t;s)\big) du.
	\end{split}
	\end{align}
	
	Recall that $q = \frac{\tau p_+}{\tau p_+ -1}$ and $B = \left\|\left(\sum_{|\alpha| = 3} \left| D^\alpha h\right|^q\right)^{1/q}\right\|_\infty$. Two applications of H{\"o}lder's inequality and above inequality~\eqref{eq:stability-4} give
	\begin{align}\label{eq:theorem:BerryEsseen-Even-p-8}
	|\mathrm{\mathbf{III}_1}|&\lesssim \mathrm{E}\left[\sum_{i=1}^n\int_0^1 \int_0^1 (1-u)\chi_i \phi\big(V^{(i)}(t;s)\big)\left(\sum_{|\alpha' |=2,\: |\alpha|=1} \Big|V_i^{\alpha'}(t;s)\dot{V}_i^\alpha(t;s)\Big|^{\tau p_+}\right)^{1/(\tau p_+)} \right. \nonumber\\
	&\quad{}\quad{} \times \left.  \left(\sum_{|\alpha'|=2} \sum_{|\alpha|=1} \Big|\big(D^{\alpha + \alpha'} h\big) \big(V^{(i)}(t;s) + uV_i(t;s)\big)\Big|^q\right)^{1/q}  dt du\right] \nonumber\\
	&\lesssim  B \mathrm{E}\left[\sum_{i=1}^n\int_0^1 \chi_i \phi\big(V^{(i)}(t;s)\big)\left(\sum_{|\alpha' |=2,\: |\alpha|=1} |V_i^{\alpha'}(t;s)\dot{V}_i^\alpha(t;s)|^{\tau p_+}\right)^{1/(\tau p_+)}  dt\right] \nonumber\\
	&\lesssim  B \mathrm{E}\left[\sum_{i=1}^n \chi_i \left( \int_0^1  \frac{\phi\big(V^{(i)}(t;s)\big)}{\sqrt{t} \wedge \sqrt{1-t}} dt\right) \left(\sum_{|\alpha|=3} \Big(|X_{ni}|^{\alpha} \vee |Y_{ni}|^{\alpha} \vee |Z_{ni}|^{\alpha} \Big)^{\tau p_+}\right)^{1/(\tau p_+)} \right],
	\end{align}
	where the last inequality follows as the third inequality in the bound on $\mathrm{\mathbf{III}_2}$.
	
	To bound the expected value in the expression in eq.~\eqref{eq:theorem:BerryEsseen-Even-p-8} we plan to apply Harris' association inequality~\citep[e.g.][Theorem 2.15]{boucheron2013concentration}. We note the following: First, the map $(x',y',z')' \mapsto (\sum_{|\alpha|=3} (|X_{ni}|^{\alpha} \vee |Y_{ni}|^{\alpha} \vee |Z_{ni}|^{\alpha}|^{\tau p_+})^{1/(\tau p_+)}$ is non-decreasing in each coordinate of $(x',y' ,z')' \in \mathbb{R}^{3d}$ (while keeping all other coordinates fixed at any value). Second, the map $(x',y',z')' \mapsto \chi_i(x,y,z) = \mathbf{1}\{\|x\|_{\tau p_+} \vee \|y\|_{\tau p_+} \vee \|z\|_{\tau p_+}\leq \xi\}$ is non-increasing in each coordinate of $(x',y' ,z')' \in \mathbb{R}^{3d}$ (while keeping all other coordinates fixed). Third, $V^{(i)}(t;s)$ and $(X_{ni}, Y_{ni}, Z_{ni})$ are independent. Therefore, Fubini's theorem and Harris' association inequality applied conditionally on $V^{(i)}(t;s)$ imply that
	\begin{align}\label{eq:theorem:BerryEsseen-Even-p-10}
	&\mathrm{E}\left[\sum_{i=1}^n \chi_i \left( \int_0^1  \frac{\phi\big(V^{(i)}(t;s)\big)}{\sqrt{t} \wedge \sqrt{1-t}} dt\right) \left(\sum_{|\alpha|=3} \Big(|X_{ni}|^{\alpha} \vee |Y_{ni}|^{\alpha} \vee |Z_{ni}|^{\alpha} \Big)^{\tau p_+}\right)^{1/(\tau p_+)} \right] \nonumber\\
	&\quad{}\leq \mathrm{E}\left[ \int_0^1\sum_{i=1}^n\mathrm{E}\left[ \chi_i \left(  \frac{\phi\big(V^{(i)}(t;s)\big)}{\sqrt{t} \wedge \sqrt{1-t}}\right) \mid V^{(i)}(t;s) \right] \right. \nonumber\\
	&\left. \quad{}\quad{} \times \mathrm{E}\left[\left(\sum_{|\alpha|=3} \Big||X_{ni}|^{\alpha} \vee |Y_{ni}|^{\alpha} \vee |Z_{ni}|^{\alpha} \Big|^{\tau p_+}\right)^{1/(\tau p_+)}\mid V^{(i)}(t;s)\right]  dt\right]\nonumber\\
	&\quad{}= \sum_{i=1}^n\mathrm{E}\left[ \chi_i \left( \int_0^1  \frac{\phi\big(V^{(i)}(t;s)\big)}{\sqrt{t} \wedge \sqrt{1-t}} dt\right) \right] \mathrm{E}\left[\left(\sum_{|\alpha|=3} \Big(|X_{ni}|^{\alpha} \vee |Y_{ni}|^{\alpha} \vee |Z_{ni}|^{\alpha} \Big)^{\tau p_+} \right)^{1/(\tau p_+)}\right].
	\end{align}
	
	To bound the first factor in eq.~\eqref{eq:theorem:BerryEsseen-Even-p-10} define $\psi(x) = \mathbf{1}\big\{ \|x\|_{p_+}  \in A^{3\delta + 6 \gamma_{p_+}} \setminus A^{-(\kappa_{p_+} + 6\gamma_{p_+})} \big\}$ and compute
	\begin{align*}
	\chi_i \left( \int_0^1  \frac{\phi\big(V^{(i)}(t;s)\big)}{\sqrt{t} \wedge \sqrt{1-t}} dt\right) = \chi_i \left( \int_0^1  \frac{\psi\big(V(t;s)\big)\phi\big(V^{(i)}(t;s)\big)}{\sqrt{t} \wedge \sqrt{1-t}} dt\right) \leq \int_0^1  \frac{\psi\big(V(t;s)\big)}{\sqrt{t} \wedge \sqrt{1-t}} dt.
	\end{align*}
	Hence,
	\begin{align}\label{eq:theorem:BerryEsseen-Even-p-12}
	&\mathrm{E}\left[\chi_i \left( \int_0^1  \frac{\phi\big(V^{(i)}(t;s)\big)}{\sqrt{t} \wedge \sqrt{1-t}} dt\right) \right]\leq \mathrm{E}\left[\int_0^1  \frac{\psi\big(V(t;s)\big)}{\sqrt{t} \wedge \sqrt{1-t}} dt \right] \nonumber\\
	&\quad{}\leq 2 \sqrt{2} \sup_{t \in [0,1]} \mathrm{P}\left(\left\|\sqrt{t}W_n(s) + \sqrt{1-t}S_n^Z\right\|_{p_+} \in A^{3\delta + 6 \gamma_{p_+}} \setminus A^{-(\kappa_{p_+} + 6\gamma_{p_+})} \right) \nonumber\\
	&\quad{} \leq 2 \sqrt{2} \sup_{s \in [0,1]} \mathrm{P}\left(\left\|W_n(s)\right\|_{p_+} \in A^{3\delta + 6 \gamma_{p_+}} \setminus A^{-(\kappa_{p_+} + 6\gamma_{p_+})} \right),
	\end{align}
	where the last inequality follows from $\sqrt{t} W_n(s) + \sqrt{1-t} S_n^Z \overset{d}{=} W_n(st)$. We bound the probability in eq.~\eqref{eq:theorem:BerryEsseen-Even-p-12} by
	\begin{align}\label{eq:theorem:BerryEsseen-Even-p-13}
	& \sup_{s \in [0,1]}\mathrm{P}\left(\|W_n(s)\|_{p_+} \in  A^{3\delta + 6 \gamma_{p_+}} \setminus A^{-(\kappa_{p_+} + 6\gamma_{p_+})} \right) \nonumber\\
	& \quad{}= 	 \sup_{s \in [0,1]} \Big\{\mathrm{P}\left(\|W_n(s)\|_{p_+} \in A^{3\delta + 6\gamma_{p_+}} \right) - \mathrm{P}\left(\|S_n^Z\|_{p_+} \in A^{3\delta + 6\gamma_{p_+}} \right)\nonumber\\
	& \phantom{\quad{}=\sup_{s \in [0,1]} \Big\{ } + \mathrm{P}\left(\|S_n^Z\|_{p_+} \in A^{-(\kappa_{p_+} + 6\gamma_{p_+})} \right) - \mathrm{P}\left(\|W_n(s)\|_{p_+} \in  A^{- (\kappa_{p_+} + 6\gamma_{p_+})} \right)\nonumber\\
	&\phantom{\quad{}=\sup_{s \in [0,1]} \Big\{ }   + \mathrm{P}\left(\|S_n^Z\|_{p_+} \in A^{3\delta + 6\gamma_{p_+}} \setminus A^{-(\kappa_{p_+} + 6\gamma_{p_+})} \right) \Big\} \nonumber\\
	\begin{split}
	&\quad{}\leq 2 \sup_{A \in \mathcal{I}}  \sup_{s \in [0,1]} \sup_{r \in \{p, p_+\}}\Big| \mathrm{P}\left(\|W_n(s)\|_r  \in A \right) - \mathrm{P}\left(\|S_n^Z\|_r   \in A \right) \Big|\\
	&\quad{}\quad{}+   \sup_{A \in \mathcal{I}}\sup_{r \in \{p, p_+\}} \mathrm{P}\left(\|S_n^Z\|_r \in A^{3\delta + \kappa_r + 12\gamma_r} \setminus A \right).
	\end{split}
\end{align}
	Combine eq.~\eqref{eq:theorem:BerryEsseen-Even-p-8}--\eqref{eq:theorem:BerryEsseen-Even-p-13} to conclude that
	\begin{align}\label{eq:theorem:BerryEsseen-Even-p-14}
	|\mathrm{\mathbf{III}_1}| &\lesssim  B \sum_{i=1}^n \mathrm{E}\left[ \left(\sum_{|\alpha|=3} \Big(|X_{ni}|^{\alpha} \vee |Y_{ni}|^{\alpha} \vee |Z_{ni}|^{\alpha} \Big)^{\tau p_+}\right)^{1/(\tau p_+)}\right] \nonumber\\
	&\phantom{B \sum_{i=1}^n \mathrm{E}}\quad{} \times \sup_{s \in [0,1]} \mathrm{P}\left(\left\|W_n(s)\right\|_{p_+} \in A^{3\delta + 6 \gamma_{p_+}} \setminus A^{-(\kappa_{p_+}+ 6\gamma_{p_+})} \right) \nonumber\\
	\begin{split}
	&\lesssim B \left(\sum_{i=1}^n \mathrm{E}[\|X_{ni}\|_{\tau p_+}^3] +  \mathrm{E}[\|Z_{ni}\|_{\tau p_+}^3]\right)\\
	&\quad{}\quad{} \times \left(\sup_{s \in [0,1]} \sup_{A \in \mathcal{I}} \sup_{r\in \{p, p_+\}}\Big| \mathrm{P}\left(\|W_n(s)\|_r  \in A \right) - \mathrm{P}\left(\|S_n^Z\|_r   \in A \right) \Big| \right.\\
	&\quad{}\quad{}\quad{}\quad{} \left. + \sup_{A \in \mathcal{I}}\sup_{r \in \{p, p_+\}} \mathrm{P}\left(\|S_n^Z\|_r \in A^{3\delta + \kappa_r + 12\gamma_r} \setminus A \right) \right),
	\end{split}
	\end{align}
	where the second inequality follows from $Y_{ni} \overset{d}{=} Z_{ni}$ for all $i=1, \ldots, n$.
	
	\textbf{Step 5. Recursive bound on eq.~\eqref{eq:theorem:BerryEsseen-Even-p-3}.} To simplify notation, let us write
	\begin{align*}
	\varrho_{n,p} &= \sup_{s \in [0,1]} \sup_{A \in \mathcal{I}} \sup_{r \in \{p, p_+\}}\Big| \mathrm{P}\left(\|W_n(s)\|_r  \in A \right) - \mathrm{P}\left(\|S_n^Z\|_r   \in A \right) \Big|.
	\end{align*}
	Recall that $q = \frac{\tau p_+}{\tau p_+-1}$ and $B = \left\|\left(\sum_{|\alpha| = 3} \left| D^\alpha h\right|^q\right)^{1/q}\right\|_\infty$. Now, the bounds from Step 1 through 4 imply that
	\begin{align}\label{eq:theorem:BerryEsseen-Even-p-15}
	\begin{split}
	\varrho_{n,p} &\lesssim B \frac{M_{n,\tau p_+}(\xi\sqrt{n})}{\sqrt{n}} + B \frac{L_{n,\tau p_+} }{\sqrt{n}} \varrho_{n,p}  \\
	&\quad{} + B\frac{L_{n,\tau p_+} }{\sqrt{n}} \left(\sup_{A \in \mathcal{I}} \sup_{r \in \{p, p_+\}} \mathrm{P}\left(\|S_n^Z\|_r \in A^{3\delta + \kappa_r + 12\gamma_r} \setminus A \right) \right)\\
	&\quad{} + \sup_{A \in \mathcal{I}} \sup_{r \in \{p, p_+\}}\mathrm{P}\left(\|S_n^Z\|_r \in A^{3\delta + 12\kappa_r} \setminus A \right).
	\end{split}
	\end{align}
	
	By Lemma~\ref{lemma:SmoothIndicator-C-Infty-Function-Even-p} and Theorem~\ref{theorem:AntiConcentration-LpNorm} we can simplify eq.~\eqref{eq:theorem:BerryEsseen-Even-p-15} to
	\begin{align}\label{eq:theorem:BerryEsseen-Even-p-16}
	\varrho_{n,p} &\leq C \left(\frac{1}{\delta^3}+ \frac{\beta}{\delta^2} + \frac{\beta^2}{\delta}\right) d^{\frac{3(\tau-1)}{\tau p_+}} \frac{M_{n,\tau p_+}(\xi\sqrt{n})}{\sqrt{n}} + C  \left(\frac{1}{\delta^3}+ \frac{\beta}{\delta^2} + \frac{\beta^2}{\delta}\right) d^{\frac{3(\tau-1)}{\tau p_+}} \frac{L_{n,\tau p_+} }{\sqrt{n}} \varrho_{n,p} \nonumber \\
	&\quad{} + C  \left(\frac{1}{\delta^3}+ \frac{\beta}{\delta^2} + \frac{\beta^2}{\delta}\right) d^{\frac{3(\tau-1)}{\tau p_+}}\frac{L_{n,\tau p_+} }{\sqrt{n}} \sup_{A \in \mathcal{I}} \sup_{r \in \{p, p_+\}}\mathrm{P}\left(\|S_n^Z\|_r \in A^{3\delta + \kappa_r + 12\gamma_r} \setminus A \right)\nonumber \\
	&\quad{} + \sup_{A \in \mathcal{I}} \sup_{r \in \{p, p_+\}} \mathrm{P}\left(\|S_n^Z\|_r \in A^{3\delta + 12\kappa_r} \setminus A \right)\nonumber \\
	\begin{split}
	&\leq C_1 \left(\frac{1}{\delta^3}+ \frac{\beta}{\delta^2} + \frac{\beta^2}{\delta}\right) d^{\frac{3(\tau-1)}{\tau p}}  \frac{M_{n,\tau p}(\xi\sqrt{n})}{\sqrt{n}} + C_1  \left(\frac{1}{\delta^3}+ \frac{\beta}{\delta^2} + \frac{\beta^2}{\delta}\right) d^{\frac{3(\tau-1)}{\tau p}}\frac{L_{n,\tau p} }{\sqrt{n}} \varrho_{n,p} \\
	&\quad{} + C_2\frac{\delta + \beta^{-1} pd^{1/(\tau p)} + \gamma_p}{\omega_p^{-1}(d, r_n)\|\sigma_n\|_p}\left(1 + C_1   \left(\frac{1}{\delta^3}+ \frac{\beta}{\delta^2} + \frac{\beta^2}{\delta}\right) d^{\frac{3(\tau-1)}{\tau p}} \frac{L_{n,\tau p} }{\sqrt{n}}\right),
	\end{split}
	\end{align}
	where $C_1, C_2 \geq 1$ are absolute constants and the last inequality holds because $p \leq p_+$ and by Remark~\ref{remark:theorem:AntiConcentration-LpNorm}. Now, set
	\begin{align*}
	\beta = p^{1/(3\tau)}d^{-(\tau-1)/(\tau p)} d^{1/(3\tau p)} n^{1/6}\overline{L}_{n,\tau p}^{-1/3} \quad{}\mathrm{and}\quad{}\delta = 6C_1p^{1- 1/(3\tau)}d^{(\tau-1)/(\tau p)} d^{2/(3\tau p)}n^{-1/6}\overline{L}_{n,\tau p}^{1/3},
	\end{align*}
	and note that
	\begin{align*}
	\left(\frac{1}{\delta^3}+ \frac{\beta}{\delta^2} + \frac{\beta^2}{\delta}\right) d^{\frac{3(\tau-1)}{\tau p}} \leq \frac{3d^{\frac{3(\tau-1)}{\tau p}}\beta^2}{\delta} =\frac{\sqrt{n}}{2C_1p^{1-1/\tau} \overline{L}_{n,\tau p}}.
	\end{align*}
	Therefore, eq.~\eqref{eq:theorem:BerryEsseen-Even-p-16} reduces to
	\begin{align*}
	\varrho_{n,p} &\leq \frac{M_{n,\tau p}(\xi\sqrt{n})}{2 p^{1 - 1/\tau} \overline{L}_{n,\tau p}} + \frac{\varrho_n}{2p^{1-1/\tau}} + \frac{3C_2(pd^{1/p})^{1-1/(3\tau)} \overline{L}_{n, \tau p}^{1/3} }{2n^{1/6}\omega_p^{-1}(d, r_n)\|\sigma_n\|_p} + \frac{3C_1\gamma_p}{2\omega_p^{-1}(d, r_n)\|\sigma_n\|_p}.
	\end{align*}
	Recall that $\gamma_p = d^{1/p(1-1/\tau)} \xi$. Set $\xi = p^{1- 1/(3\tau)}n^{-1/6}\overline{L}_{n,\tau p}^{1/3}$ and solve above inequality for $\varrho_{n,p}$ to obtain
	\begin{align*}
	\varrho_{n,p} &\lesssim \frac{M_{n, \tau p}\big(p^{1-1/(3\tau)}n^{1/3}\overline{L}_{n,\tau p}^{1/3}\big)}{p^{1 - 1/\tau} \overline{L}_{n,\tau p}} + \frac{(pd^{1/p})^{1-1/(3\tau)}\overline{L}_{n, \tau p}^{1/3}}{n^{1/6}\omega_p^{-1}(d, r_n)\|\sigma_n\|_p}.
	\end{align*}
	
	Lastly, note that $\sup_{t \geq 0} \Big| \mathrm{P}\big(\|S_n^X\|_p  \leq t \big) - \mathrm{P}\big(\|S_n^Z\|_p \leq t \big) \Big| \leq \varrho_{n,p}$. This concludes the proof of the first statement.
	
	\textbf{Proof of Case (ii).} The proof is identical to the proof of the first statement except for the following four changes: First, the fundamental smoothing inequality is the inequality that directly precedes inequality~\eqref{eq:theorem:BerryEsseen-Even-p-3}; there is no need to introduce $p_+$, the smallest even exponent larger than $p$. Second, replace arguments involving H{\"o}lder's inequality with the conjugate exponents $(q, \tau p_+)$ by arguments based on H{\"o}lder's inequality with the conjugate exponents $(1,\infty)$. Third, replace Lemma~\ref{lemma:SmoothIndicator-C-Infty-Function-Even-p} by Lemma~\ref{lemma:SmoothIndicator-C-Infty-Function-Large-p} throughout. Lastly, set
	\begin{align*}
	\beta = n^{1/6} (\log d)^{1/3} \overline{L}_{n,\infty}^{-1/3}, \hspace{25pt}
	\delta = 6C_1 (\log d)^{2/3} n^{-1/6} \overline{L}_{n,\infty}^{1/3},\hspace{25pt}
	\xi = (\log d)^{2/3} n^{-1/6} \overline{L}_{n,\infty}^{1/3},
	\end{align*}
	and proceed as in Step 5. This concludes the proof of the second statement.
\end{proof}

\subsubsection{Proofs for Appendix~\ref{subsec:AntiConcentration}}
\begin{proof}[\textbf{Proof of Theorem~\ref{theorem:AntiConcentration-LpNorm-Even-p}}]
	Observe that for $p \in \mathbb{N}$ even, $\|X\|_p^p - t = \sum_{j=1}^dX_j^p - t$ is a (multivariate) polynomial of degree $p$ in $X \in \mathbb{R}^d$. Therefore, by Theorem 8 in~\cite{carbery2001Distributional} uniformly in $t \geq 0$, $q \geq 1$, and $p \in \mathbb{N}$ even,
	\begin{align}\label{eq:lemma:AntiConcentration-LpNorm-Even-p-1}
	q^{-1} \left(\mathrm{E}\big|\|X\|_p^p - t\big|^{q/p}\right)^{1/q} \mathrm{P}\left(\big|\|X\|_p^p - t\big| \leq \varepsilon^p\right) \lesssim \varepsilon.
	\end{align}
	Furthermore, note that for all $p,q \geq 1$ and any pair $Z, Z'$ of independent and identically distributed random variables,
	\begin{align*}
	\left(\mathrm{E}\big|Z - Z' \big|^{q/p}\right)^{p/q} \leq (2^{p/q} \vee 2 ) \: \left(\mathrm{E}|Z|^{q/p}\right)^{p/q}.
	\end{align*}
	Thus, for all $t \geq 0$, $q \geq 1$, and $p \in \mathbb{N}$ even,
	\begin{align}\label{eq:lemma:AntiConcentration-LpNorm-Even-p-2}
	\left(\mathrm{E}\big|\|X\|_p^p - t\big|^{q/p}\right)^{1/q} &\gtrsim \left(\mathrm{E}\big|\|X\|_p^p - \|X' \|_p^p\big|^{q/p}\right)^{1/q} \gtrsim \left(\mathrm{E}\big|\|X\|_p - \|X' \|_p\big|^q\right)^{1/q},
	\end{align}
	where the second inequality follows from the reverse triangle inequality applied to $| \cdot |^{1/p}$.
	
	Combine eq.~\eqref{eq:lemma:AntiConcentration-LpNorm-Even-p-1} and eq.~\eqref{eq:lemma:AntiConcentration-LpNorm-Even-p-2} to conclude that
	\begin{align}\label{eq:lemma:AntiConcentration-LpNorm-Even-p-3}
	\sup_{q \geq 1} q^{-1}\left(\mathrm{E}\big|\|X\|_p - \|X' \|_p\big|^q \right)^{1/q} \sup_{t \geq 0}\mathrm{P}\left(\big|\|X\|_p^p - t\big| \leq \varepsilon^p\right) \lesssim \varepsilon.
	\end{align}
	This is a statement about the polynomial $\|X\|_p^p - t = \sum_{j=1}^dX_j^p - t$. We reduce it to a statement about $\|X\|_p$ by lower bounding the probability as follows:
	\begin{align}\label{eq:lemma:AntiConcentration-LpNorm-Even-p-4}
	\begin{split}
	\sup_{t \geq0} \mathrm{P}\left( t \leq  \|X\|_p \leq t + \varepsilon\right)	&= \sup_{t \geq 0} \mathrm{P}\Big(t^p \leq \|X\|_p^p \leq (t + \varepsilon)^p \Big) \\
	&\leq \sup_{t \geq 0} \mathrm{P}\left(t^p - 2^{p-1} \varepsilon^p \leq \|X\|_p^p \leq 2^{p-1} t^p + 2^{p-1} \varepsilon^p\right) \\
	& = \sup_{t \geq 0}  \mathrm{P}\left( \big\{ \left|\|X\|_p^p - t^p \right| \leq 2^{p-1}\varepsilon^p \big\} \cup \big\{ \left|\|X\|_p^p - 2^{p-1}t^p \right| \leq 2^{p-1}\varepsilon^p \big\}\right) \\
	&\leq 2 \sup_{t \geq 0}  \mathrm{P}\left(\left| \|X\|_p^p - t \right| \leq 2^{p-1}\varepsilon^p \right),
	\end{split}
	\end{align}
	where the first inequality holds since $(a + b)^p \leq 2^{p-1}a^p + 2^{p-1}b^p$ for all $a, b \geq 0$ and $p > 1$. Hence, by eq.~\eqref{eq:lemma:AntiConcentration-LpNorm-Even-p-3},
	\begin{align}\label{eq:lemma:AntiConcentration-LpNorm-Even-p-5}
	\sup_{q \geq 1} q^{-1} \left(\mathrm{E}\big|\|X\|_p - \|X' \|_p\big|^q\right)^{1/q}\sup_{t \geq0} \mathrm{P}\left( t \leq  \|X\|_p \leq t + \varepsilon\right) \lesssim \varepsilon.
	\end{align}
\end{proof}
\begin{proof}[\textbf{Proof of Corollary~\ref{corollary:AntiConcentration-L2Norm}}]
	Note that by eq.~\eqref{eq:lemma:AntiConcentration-LpNorm-Even-p-2},
	\begin{align*}
	\sup_{q \geq 1} q^{-1}\left(\mathrm{E}\big|\|X\|_2^2 - \|X'\|_2^2\big|^{q/2}\right)^{1/q}\sup_{t \geq0} \mathrm{P}\left( t \leq  \|X\|_2 \leq t + \varepsilon\right) \lesssim \varepsilon.
	\end{align*}
	The expected value on the left hand side can be lower bounded as
	\begin{align*}
	\sup_{q \geq 1} q^{-1}\left(\mathrm{E}\big|\|X\|_2^2 - \|X'\|_2^2\big|^{q/2}\right)^{1/q} &\geq \frac{1}{4} \left(\mathrm{E}\big|\|X\|_2^2 - \|X'\|_2^2\big|^2\right)^{1/4}\\
	& = \frac{1}{2}\mathrm{Var}[\|X\|_2^2]^{1/4}\\
	& = \frac{1}{2} \left(\mathrm{tr}(\Sigma^2)  + \mu'\Sigma \mu  \right)^{1/4},
	\end{align*}
	where the last line follows from direct calculations. This completes the proof.
\end{proof}
\begin{proof}[\textbf{Proof of Theorem~\ref{theorem:AntiConcentration-LpNorm}}]
	We split the proof into two parts. First, for $p \in [1, \infty)$ arbitrary we show that $\sup_{t \geq0} \mathrm{P}( t \leq  \|X\|_p \leq t +  \varepsilon p^{-1/2} r^{-1/(2p)} \|\sigma\|_p)  \lesssim \varepsilon$. Then, for $p \geq \log d$, we show that $\sup_{t \geq0} \mathrm{P}( t \leq  \|X\|_p \leq t +  \varepsilon \|\sigma\|_p/\sqrt{\log d})  \lesssim \varepsilon$.
	
	\textbf{Step 1.} Let $p \geq 1$ be arbitrary and define $p_+ = 2\lceil \frac{p}{2}\rceil$ to be the smallest even integer larger than (or equal to) $p$. Then, $\|x\|_{p_+}\leq \|x\|_p$ for all $x \in \mathbb{R}^d$ and  therefore
	\begin{align}\label{eq:theorem:AntiConcentration-LpNorm-General-p-1}
	\sup_{t \geq0} \mathrm{P}\left( t \leq  \|X\|_p \leq t + \varepsilon\right) \leq \sup_{t \geq0} \mathrm{P}\left( t \leq \|X\|_{p_+} \leq t + \varepsilon\right) \leq 2 \sup_{t \geq 0}  \mathrm{P}\left(\left| \|X\|_{p_+}^{p_+} - t \right| \leq 2^{p_+-1}\varepsilon^{p_+}\right),
	\end{align}
	where the second inequality follows as in eq.~\eqref{eq:lemma:AntiConcentration-LpNorm-Even-p-4}. Thus, by eq.~\eqref{eq:lemma:AntiConcentration-LpNorm-Even-p-2} and~\eqref{eq:lemma:AntiConcentration-LpNorm-Even-p-3}, we have
	\begin{align}\label{eq:theorem:AntiConcentration-LpNorm-General-p-2}
	\sup_{q \geq 1} q^{-1}\left(\mathrm{E}\big|\|X\|_{p_+}^{p_+} - \|X' \|_{p_+}^{p_+}\big|^{q/p_+}\right)^{1/q} \sup_{t \geq0} \mathrm{P}\left( t \leq  \|X\|_p \leq t + \varepsilon\right)  \lesssim \varepsilon.
	\end{align}
	Set $q = 2p_+$ and lower bound the expected value on the left hand side in above display as follows
	\begin{align}\label{eq:theorem:AntiConcentration-LpNorm-General-p-3}
	\begin{split}
	\sup_{q \geq 1} q^{-1}\left(\mathrm{E}\big|\|X\|_{p_+}^{p_+} - \|X' \|_{p_+}^{p_+}\big|^{q/p_+}\right)^{1/q} &\geq \frac{1}{2p_+}\left( \mathrm{E}\big|\|X\|_{p_+}^{p_+} - \|X' \|_{p_+}^{p_+}\big|^{2}\right)^{1/(2p_+)}\\
	&= p_+^{-1}\mathrm{Var}\left[\|X\|_{p_+}^{p_+}\right]^{1/(2p_+)}.
	\end{split}
	\end{align}
	
	Since $\Sigma$ is positive semi-definite and has rank $r$, there exists a lower triangular matrix $\Gamma$ such that $\Gamma\Gamma' = \Sigma$ and $\gamma_{kj} = 0$ for all $(k,j)$ which satisfy $k < j$ or $j > r$  (via $LDL'$ decomposition). Let $Z \sim N(0, I_d)$ be a standard normal random vector in $\mathbb{R}^d$ and set $X \overset{d}{=} \Gamma Z$. By Lemma 11.1 in~\cite{chatterjee2014Superconcentration} and Cauchy-Schwarz we have
	\begin{align}\label{eq:theorem:AntiConcentration-LpNorm-General-p-4}
	\begin{split}
	\mathrm{Var}\left[\|X\|_{p_+}^{p_+}\right] &\geq \frac{1}{2} \sum_{j=1}^d \left(\mathrm{E}\left[ p_+ \sum_{k=1}^d \gamma_{kj} Z_{j} (\gamma_k'Z)^{p_+-1} \right]\right)^2 \\
	&= \frac{1}{2} \sum_{j=1}^{r} \left(\mathrm{E}\left[ p_+ \sum_{k=1}^d \gamma_{kj} Z_{j} (\gamma_k'Z)^{p_+-1} \right]\right)^2 \\
	&\geq \frac{p_+^2}{2r} \left(\mathrm{E}\left[\sum_{k=1}^d(\gamma_k'Z)^{p_+}\right] \right)^2\\
	&=  \frac{p_+^2}{2r} \left(\mathrm{E}\|X\|_{p_+}^{p_+} \right)^2.
	\end{split}
	\end{align}
	By Stirling's formula the central moments of a standard normal random variable $Z_1$ satisfy the following asymptotic estimate:
	\begin{align*}
	\mathbb{E}[|Z_1|^{p}] = \frac{2^{p/2}}{\sqrt{\pi}} \Gamma\left(\frac{p +1}{2}\right) \asymp 2^{p/2} \left(\frac{p}{e}\right)^{p/2} \hspace{10pt} \mathrm{as} \hspace{10pt} p \rightarrow \infty.
	\end{align*}
	Thus, since $p \leq p_+ \leq 2p$, it follows that
	\begin{align}\label{eq:theorem:AntiConcentration-LpNorm-General-p-5}
	\left(\mathrm{E}\|X\|_{p_+}^{p_+}\right)^{1/p_+} \gtrsim \left(\mathrm{E} \|X\|_{p}^{p}\right)^{1/p} \gtrsim \left(\sum_{j=1}^d \sigma_j^{p}\right)^{1/p} \left(\frac{p}{e}\right)^{1/2} \gtrsim p^{1/2} \|\sigma\|_p.
	\end{align}
	Combine eq.~\eqref{eq:theorem:AntiConcentration-LpNorm-General-p-3}--\eqref{eq:theorem:AntiConcentration-LpNorm-General-p-5} to conclude that
	\begin{align*}
	\sup_{q \geq 1} q^{-1}\left(\mathrm{E}\big|\|X\|_{p_+}^{p_+} - \|X' \|_{p_+}^{p_+}\big|^{q/p_+}\right)^{1/q} &\gtrsim p_+^{-1} p^{1/2} r^{-1/(2p)} \|\sigma\|_p \gtrsim p^{-1/2} r^{-1/(2p)} \|\sigma\|_p.
	\end{align*}
	Hence, by eq.~\eqref{eq:theorem:AntiConcentration-LpNorm-General-p-2}, for any $p \in [1, \infty)$,
	\begin{align}\label{eq:theorem:AntiConcentration-LpNorm-General-p-6}
	\sup_{t \geq0} \mathrm{P}\left( t \leq  \|X\|_p \leq t +  \varepsilon p^{-1/2} r^{-1/(2p)} \|\sigma\|_p\right)  \lesssim \varepsilon.
	\end{align}
	
	\textbf{Step 2.} Let $p \geq \log d$ be arbitrary. Note that $r^{1/(2\log d)} \leq e^{1/2}$. Also, $\|x\|_p \leq \|x\|_{\log d} \leq e \|x\|_p$. Therefore, by Step 1,
	\begin{align}\label{eq:theorem:AntiConcentration-LpNorm-General-p-7}
	\sup_{ t \geq 0 }\mathrm{P}\left( t \leq \|X\|_p \leq t + \varepsilon \frac{\|\sigma\|_p}{\sqrt{\log d}}\right) &= \sup_{ t \geq 0 }\mathrm{P}\left( et \leq e\|X\|_p \leq et + e\varepsilon \frac{\|\sigma\|_p}{\sqrt{\log d}}\right) \nonumber\\
	&\leq \sup_{ t \geq 0 }\mathrm{P}\left( et \leq \|X\|_{\log d} \leq et + \varepsilon \frac{ e^{3/2} \|\sigma\|_{\log d}}{\sqrt{\log d} \: r^{1/(2\log d)}}\right) \nonumber\\
	&  \lesssim  \varepsilon.
	\end{align}
	
	To combine Step 1 and 2 as in the statement of the theorem, simply note that for $p \leq \log d$, we have
	\begin{align*}
	r^{1/p} p \leq d^{1/p} p \leq e \log d.
	\end{align*}
	Hence, $\|\sigma\|_p p^{-1/2} r^{1/(2p)} \geq \|\sigma\|_p /\sqrt{\log d}$. For $p \geq \log d$, the inequality is reversed.
\end{proof}

\subsubsection{Proofs for Appendix~\ref{subsec:GaussianComparison}}
\begin{proof}[\textbf{Proof of Theorem~\ref{theorem:Gaussian-Comparison-KolmogorovDistance}}]
	\noindent
	
	\textbf{Proof of Case (i).}
	
	\textbf{Step 1. Fundamental smoothing inequality.} Let $Z$ be an independent copy of $Y$ and define
	\begin{align*}
	W(s) := \sqrt{s}X + \sqrt{1-s} Z, \hspace{10pt} s \in [0,1].
	\end{align*}
	Consider the family of sets $\mathcal{I} = \{A \subseteq \mathbb{R} : A = [0,t], t \geq 0 \}$.  Let $p \in [1, \infty)$ be arbitrary. Define $p_+ = 2\lceil \frac{p}{2}\rceil$ to be the smallest even integer larger than (or equal to) $p$.	By Lemma~\ref{lemma:SmoothIndicator-C-Infty-Function-Even-p} for $A \in \mathcal{I}$, we have
	\begin{align*}
	\mathrm{P}\left(\|W(s)\|_p  \in A \right) - \mathrm{P}\left(\|Y\|_p   \in A^{12\kappa_p + 3\delta}  \right) \leq \mathrm{E}\left[h_{p_+, d, \beta, A^{3\kappa_p}}\big(W(s)\big)- h_{p++, d, \beta, A^{3\kappa_p}}(Y)\right].
	\end{align*}
	Re-arrange the terms in above inequality and take the supremum over $A \in \mathcal{A}$ to obtain
	\begin{align}\label{eq:theorem:GaussianComparison-1}
	\begin{split}
	&\sup_{A \in \mathcal{A}} \Big(\mathrm{P}\left(\|W(s)\|_p  \in A \right) - \mathrm{P}\left(\|Y\|_p   \in A\right) \Big)\\
	&\quad{}\leq \sup_{A \in \mathcal{A}}\mathrm{P}\left(\|Y\|_p   \in A^{12\kappa_p + 3\delta} \setminus A \right) + \sup_{A \in \mathcal{I}} \Big|\mathrm{E}\left[h_{p_+, d, \beta,\delta, A}\big(W(s)\big)- h_{p_+, d, \beta,\delta, A}(Y)\right]\Big|,
	\end{split}
	\end{align}
	By Lemma~\ref{lemma:SmoothIndicator-C-Infty-Function-Even-p} we also have
	\begin{align*}
	&\mathrm{P}\left(\|Y\|_p  \in A^{-(12\kappa_p + 3\delta)} \right) - \mathrm{P}\left(\|W(s)\|_p   \in A \right)\\
	&\quad{}\leq \mathrm{E}\left[h_{p_+, d, \beta, A^{-(12\kappa_p + 3\delta)}}(T)- h_{p_+, d, \beta, A^{-(12\kappa_p + 3\delta)}}\big(W(s)\big)\right].
	\end{align*}
	Observe that $\sup_{A \in \mathcal{I}}\mathrm{P}\left( \|Y\|_p \in A \setminus A^{-(12\kappa_p + 3\delta)}\right) \leq \sup_{A \in \mathcal{I}}\mathrm{P}\left( \|Y\|_p \in A^{12\kappa_p + 3\delta}\setminus A\right)$. Together with the preceding inequality this yields
	\begin{align}\label{eq:theorem:GaussianComparison-2}
	\begin{split}
	&\sup_{A \in \mathcal{I}} \Big( \mathrm{P}\left(\|Y\|_p  \in A \right) - \mathrm{P}\left(\|W(s)\|_p   \in A \right) \Big)\\
	&\quad{}\quad{}\leq \sup_{A \in \mathcal{I}}\mathrm{P}\left(\|Y\|_p  \in A^{12\kappa_p + 3\delta} \setminus A \right) + \sup_{A \in \mathcal{I}} \Big|\mathrm{E}\left[h_{p_+, d, \beta, \delta,A}\big(W(s)\big)- h_{p_+, d, \beta,\delta, A}(Y)\right]\Big|,
	\end{split}
	\end{align}
	Combine eq.~\eqref{eq:theorem:GaussianComparison-1} and eq.~\eqref{eq:theorem:GaussianComparison-2} to obtain
	\begin{align*}
	&\sup_{s \in [0,1]}\sup_{A \in \mathcal{I}} \Big| \mathrm{P}\left(\|Y\|_p  \in A \right) - \mathrm{P}\left(\|W(s)\|_p   \in A \right) \Big|\\
	&\quad{}\leq \sup_{A \in \mathcal{I}}\mathrm{P}\left(\|Y\|_p \in A^{12\kappa_p + 3\delta} \setminus A \right) + \sup_{s \in [0,1]}\sup_{A \in \mathcal{I}} \Big|\mathrm{E}\left[h_{p_+, d, \beta, \delta,A}\big(W(s)\big)- h_{p_+, d, \beta,\delta, A}(Y)\right]\Big|
	\end{align*}
	Note that above inequality also holds for $p = p_+$. Thus, we have the following fundamental smoothing inequality
	\begin{align}\label{eq:theorem:GaussianComparison-3}
	\begin{split}
	&\sup_{s \in [0,1]}\sup_{A \in \mathcal{I}} \sup_{r \in \{p, p_+\}}\Big| \mathrm{P}\left(\|Y\|_r \in A \right) - \mathrm{P}\left(\|W(s)\|_r   \in A \right) \Big| \\
	&\quad{}\leq \sup_{A \in \mathcal{I}} \sup_{r \in \{p, p_+\}}\mathrm{P}\left(\|Y\|_r \in A^{12\kappa_r + 3\delta} \setminus A \right) + \sup_{s \in [0,1]}\sup_{A \in \mathcal{I}} \Big|\mathrm{E}\left[h_{p_+, d, \beta, \delta,A}\big(W(s)\big)- h_{p_+, d, \beta,\delta, A}(Y)\right]\Big|
	\end{split}
	\end{align}
	We now bound the second term on the right hand side of eq.~\eqref{eq:theorem:GaussianComparison-3}.
	
	We now bound the second term.
	
	\textbf{Step 2. Stein's Lemma.}
	Define the Slepian-Stein (double) interpolant as
	\begin{align*}
	V(t;s) := \sqrt{t} W(s) + \sqrt{1 - t}Y, \hspace{10pt} s,t \in [0,1],
	\end{align*}
	and its derivative with respect to $t$ by
	\begin{align*}
	\dot{V}(t;s) :=\frac{d}{dt}V(t;s) = \frac{1}{2} \left[\frac{1}{\sqrt{t}} \left(\sqrt{s}X + \sqrt{1-s} Z\right)  - \frac{1}{\sqrt{1-t}} Y \right].
	\end{align*}
	Since $V(0;s) = Y$ and $V(1;s) = W(s)$, the mean value theorem gives
	\begin{align}\label{eq:theorem:GaussianComparison-4}
	\begin{split}
	&\mathrm{E}\left[h_{p_+, d, \beta, \delta,A}\big(W(s)\big)- h_{p_+, d, \beta,\delta, A}(Y) \right]\\
	&\quad{} = \sum_{|\alpha|=1} \int_0^1 \mathrm{E}\left[\dot{V}^\alpha(t;s)\big(D^\alpha h_{p_+, d, \beta, \delta,A} \big) \big(V(t;s)\big)\right] dt.
	\end{split}
	\end{align}
	For brevity of notation, we now drop the subscripts $p_+,d,\beta, \delta, A$ and write $h$ instead of $h_{p_+, d, \beta,\delta, A}$. By Stein's identity, for multi-indices $\alpha$, $\alpha'$ with $|\alpha| = |\alpha'| = 1$,
	\begin{align*}
	\mathrm{E}\big[W^\alpha(s) (D^\alpha h)\big(V(t;s)\big)\big] &= \sqrt{t}\sum_{|\alpha'|=1} \mathrm{E}\big[W^\alpha(s) W^{\alpha'}(s)\big] \mathrm{E}\big[(D^{\alpha + \alpha'}h )\big(V(t;s)\big)\big],\\
	\mathrm{E}\big[Y^\alpha (D^\alpha h)\big(V(t;s)\big)\big] &= \sqrt{1-t}\sum_{|\alpha'|=1} \mathrm{E}\big[Y^\alpha Y^{\alpha'}\big] \mathrm{E}\big[(D^{\alpha + \alpha'}h)\big(V(t;s)\big)\big].
	\end{align*}
	Hence, eq.~\eqref{eq:theorem:GaussianComparison-4} simplifies to
	\begin{align}\label{eq:theorem:GaussianComparison-5}
	&\mathrm{E}\left[h\big(W(s)\big) - h(Y) \right]\nonumber\\
	&\quad{}= \frac{1}{2}\sum_{|\alpha|=1} \sum_{|\alpha'|=1} \int_0^1 \mathrm{E}\big[W^\alpha(s) W^{\alpha'}(s) - Y^\alpha Y^{\alpha'}\big] \mathrm{E}\big[(D^{\alpha + \alpha'}h)\big(V(t;s)\big)\big] dt \nonumber\\
	&\quad{}= \frac{s}{2}\sum_{|\alpha|=1} \sum_{|\alpha'|=1} \int_0^1 \mathrm{E}\big[X^\alpha X^{\alpha'} - Y^\alpha Y^{\alpha'}\big] \mathrm{E}\big[(D^{\alpha + \alpha'}h)\big(V(t;s)\big)\big] dt,
	\end{align}
	where the second equality follows since $Z$ is an independent copy of $Y$.
	
	Recall from Lemma~\ref{lemma:SmoothIndicator-C-Infty-Function-Even-p} that $h$ is non-constant on the set $\big\{z \in \mathbb{R}^d: M_{p_+, \kappa_{p_+}}(z) \in A^{3\delta} \setminus A\big\}$ only. Set $\phi(x) = \mathbf{1}\big\{ \|x\|_{p_+}  \in A^{3\delta} \setminus A^{-\kappa_{p_+}} \big\}$ and note that $\phi(x)=1$ if $x \in \big\{z \in \mathbb{R}^d: M_{p_+, \kappa_{p_+}}(z) \in A^{3\delta} \setminus A\big\}$. Therefore, for all $s \in [0,1]$, the term on the far right hand side in eq.~\eqref{eq:theorem:GaussianComparison-5} is not larger than
	\begin{align}\label{eq:theorem:GaussianComparison-6}
	\frac{s}{2}\sum_{|\alpha|=1} \sum_{|\alpha'|=1} \int_0^1 \mathrm{E}\big[X^\alpha X^{\alpha'} - Y^\alpha Y^{\alpha'}\big] \mathrm{E}\big[\phi\big(V(t;s)\big)(D^{\alpha + \alpha'}h)\big(V(t;s)\big)\big] dt.
	\end{align}
	By H{\"o}lder's inequality, for $1/p_+ + 1/q = 1$,
	\begin{align}\label{eq:theorem:GaussianComparison-7}
	&\frac{s}{2}\sum_{|\alpha|=1} \sum_{|\alpha'|=1} \int_0^1 \mathrm{E}\big[X^\alpha X^{\alpha'} - Y^\alpha Y^{\alpha'}\big] \mathrm{E}\big[\phi\big(V(t;s)\big)(D^{\alpha + \alpha'}h)\big(V(t;s)\big)\big] dt \nonumber\\
	&\hspace{5pt}\leq \frac{s}{2}\left(\sum_{|\alpha|=1} \sum_{|\alpha'|=1} \big|\mathrm{E}\big[X^\alpha X^{\alpha'} - Y^\alpha Y^{\alpha'}\big]\big|^{p_+} \right)^{1/p_+} \nonumber\\
	&\hspace{45pt} \times  \int_0^1\left(\sum_{|\alpha|=1} \sum_{|\alpha'|=1} \big|\mathrm{E}\big[\phi\big(V(t;s)\big)(D^{\alpha + \alpha'}h)\big(V(t;s)\big)\big] \big|^q\right)^{1/q} dt \nonumber\\
	\begin{split}
	&\hspace{5pt}\leq \frac{1}{2}\big\|\mathrm{vec}(\Sigma_X - \Sigma_Y)\big\|_{p_+} \left\|\left(\sum_{|\alpha| = 2} \left|D^\alpha  h\right|^q\right)^{1/q}\right\|_\infty \\
	& \hspace{45pt} \times\mathrm{P}\left(\left\|\sqrt{t} W(s) + \sqrt{1 -t}Y \right\|_{p_+}  \in A^{3\delta}\setminus A^{-\kappa_{p_+}} \right).
	\end{split}
	\end{align}
	
	We bound the last factor in eq.~\eqref{eq:theorem:GaussianComparison-7} by
	\begin{align}\label{eq:theorem:GaussianComparison-8}
	&\sup_{s \in [0,1]}\mathrm{P}\left(\left\|\sqrt{t} W(s) + \sqrt{1 -t}Y \right\|_{p_+}  \in A^{3\delta}\setminus A^{-\kappa_{p+}} \right) \nonumber\\
	&\quad{} \leq \sup_{s \in [0,1]} \mathrm{P}\left(\|W(s)\|_{p_+}  \in A^{3\delta }\setminus A^{-\kappa_{p_+}} \right) \nonumber\\
	& \quad{}= 	\sup_{s \in [0,1]} \Big(\mathrm{P}\left(\|W(s)\|_{p_+}  \in A^{3\delta} \right) - \mathrm{P}\left(\|Y\|_{p_+}  \in A^{3\delta} \right) - \mathrm{P}\left(\|W(s)\|_{p_+} \in  A^{- \kappa_{p_+}}\right) \Big) \nonumber\\
	&\quad{}\quad{} + \mathrm{P}\left(\|Y\|_{p_+}  \in A^{- \kappa_{p_+}} \right)  + \mathrm{P}\left(\|Y\|_{p_+} \in A^{3\delta } \setminus A^{- \kappa_{p_+}} \right)\nonumber\\
	\begin{split}
	&\quad{}\leq 2 \sup_{A \in \mathcal{I}} \sup_{r \in \{p, p_+\}}\Big| \mathrm{P}\left(\|W(s)\|_r  \in A \right) - \mathrm{P}\left(\|Y\|_r   \in A \right) \Big|\\
	&\quad{}\quad{} +   \sup_{A \in \mathcal{I}} \sup_{r \in \{p, p_+\}}\mathrm{P}\left(\|Y\|_r \in A^{3\delta + \kappa_r} \setminus A \right),
	\end{split}
	\end{align}
	where the first inequality holds since $\sqrt{t} W(s) + \sqrt{1 -t}Y \overset{d}{=} W(st)$.
	
	\textbf{Step 3. Recursive bound on eq.~\eqref{eq:theorem:GaussianComparison-3}.}
	To simplify notation we define
	\begin{align*}
	\varrho_X &= \sup_{s \in [0,1]} \sup_{A \in \mathcal{I}} \sup_{r \in \{p, p_+\}} \Big| \mathrm{P}\left(\|W(s)\|_r  \in A \right) - \mathrm{P}\left(\|Y\|_r   \in A \right) \Big|.
	\end{align*}
	
	Combine eq.~\eqref{eq:theorem:GaussianComparison-3} and eq.~\eqref{eq:theorem:GaussianComparison-5}--\eqref{eq:theorem:GaussianComparison-8}, and conclude that
	\begin{align}\label{eq:theorem:GaussianComparison-10}
	\begin{split}
	\varrho_X &\leq  \Delta_{p_+} \left\|\left(\sum_{|\alpha| = 2} \left|D^\alpha  h\right|^q\right)^{1/q}\right\|_\infty \varrho_X\\
	&\quad{} + \frac{\Delta_{p_+} }{2}\left\|\left(\sum_{|\alpha| = 2} \left|D^\alpha  h\right|^q\right)^{1/q}\right\|_\infty \sup_{A \in \mathcal{I}} \sup_{r \in \{p, p_+\}}\mathrm{P}\big(\|Y\|_r \in A^{3 \delta + \kappa_r} \setminus A\big)\\
	&\quad{} + \sup_{A \in \mathcal{I}} \sup_{r \in \{p, p_+\}} \mathrm{P}\big(\|Y\|_r \in A^{3 \delta + 12 \kappa_r} \setminus A\big).
	\end{split}
	\end{align}
	
	By Lemma~\ref{lemma:SmoothIndicator-C-Infty-Function-Even-p} and Theorem~\ref{theorem:AntiConcentration-LpNorm}, eq.~\eqref{eq:theorem:GaussianComparison-10} reduces to
	\begin{align}\label{eq:theorem:GaussianComparison-11}
	\begin{split}
	\varrho_X &\leq C_1\left(\delta^{-2} + \frac{\beta}{\delta}\right)\Delta_{p_+} \varrho_X  + \sup_{A \in \mathcal{I}} \sup_{r \in \{p, p_+\}}\mathrm{P}\big(\|Y\|_r \in A^{3\delta + 12\kappa_r } \setminus A\big) \\
	&\quad{} + \frac{C_1}{2}\left(\delta^{-2}+ \frac{\beta}{\delta}\right) \Delta_{p_+} \left(\sup_{A \in \mathcal{I}}\sup_{r \in \{p, p_+\}}\mathrm{P}\big(\|Y\|_r \in A^{3 \delta + \kappa_r} \setminus A\big) \right)\\
	&\leq C_1\left(\delta^{-2} + \frac{\beta}{\delta}\right)\Delta_p \varrho_X  + C_2\frac{\delta + \beta^{-1} pd^{1/p}}{\omega_p^{-1}(d, r_X)\|\sigma_X\|_p} \left(1 + C_1  \left(\delta^{-2}+ \frac{\beta}{\delta}\right) \Delta_p\right),
	\end{split}
	\end{align}
	where $C \geq 1$ is an absolute constant and the last inequality holds since $p \leq p_+$ and by Remark~\ref{remark:theorem:AntiConcentration-LpNorm}.
	
	Set $\beta = p^{1/2}d^{1/(2p)}\Delta_p^{-1/2}$ and $\delta = 4 C_1 p^{1/2} d^{1/(2p)} \Delta_p^{1/2}$. Note that $\delta^{-2} + \delta^{-1} \beta \leq 2\delta^{-1}\beta$. Thus, eq.~\eqref{eq:theorem:GaussianComparison-11} simplifies to
	\begin{align*}
	\varrho_X &\leq \frac{\varrho_X}{2} + \frac{3C_2}{2}\frac{p^{1/2}d^{1/(2p)} \Delta_p^{1/2} }{\omega_p^{-1}(d, r_X)\|\sigma_X\|_p},
	\end{align*}
	which implies
	\begin{align}\label{eq:theorem:GaussianComparison-12}
	\varrho_X \lesssim \frac{p^{1/2}d^{1/(2p)} \Delta_p^{1/2} }{\omega_p^{-1}(d, r_X)\|\sigma_X\|_p}.
	\end{align}
	
	Since $X$ and $Y$ are both Gaussian, we can interchange their role in the proof and obtain analogous bounds on $\varrho_Y$ involving $\sigma_Y$ and $r_Y$. Since $\varrho_X$ and $\varrho_Y$ both upper bound $|\mathrm{P}\left(\|X\|_p  \in A \right) - \mathrm{P}\left(\|Y\|_p \in A \right)|$ the first claim of Theorem~\ref{theorem:Gaussian-Comparison-KolmogorovDistance} follows.

	\textbf{Proof of Case (ii).} We split the proof into two parts.
	
	\textbf{Step 1.} We derive the bound involving $\Delta_\infty$. The proof of this result is identical to the proof of the four statement except for the following three changes: First, we do not need to introduce $p_+$, the smallest even integer larger than $p$. Instead, as fundamental smoothing inequality we may take the inequality directly preceding eq.~\eqref{eq:theorem:GaussianComparison-3}. Second, we replace arguments involving H{\"o}lder's inequality with the conjugate exponents $(q,p_+)$ by arguments based on H{\"o}lder's inequality with the conjugate exponents $(1,\infty)$. Third, replace Lemma~\ref{lemma:SmoothIndicator-C-Infty-Function-Even-p} by Lemma~\ref{lemma:SmoothIndicator-C-Infty-Function-Large-p} throughout. Lastly, set
	\begin{align*}
	\beta =  (\log d)^{1/2}\Delta_\infty^{-1/2} \hspace{25pt} \mathrm{and} \hspace{25pt}
	\delta = 4C_1 (\log d)^{1/2} \Delta_\infty^{1/2},
	\end{align*}
	and proceed as in Step 5.
	
	\textbf{Step 2.} To derive the bound involving $\Delta_{op}$ we have to make the following changes: Denote by $\nabla^2h$ the Hessian of $h$. Then, by H{\"o}lder's inequality for matrix inner products (and a rough upper bound following from Gershgorin's circle theorem) we can upper bound eq.~\eqref{eq:theorem:GaussianComparison-6} as follows:
	\begin{align*}
	&\frac{s}{2}\sum_{|\alpha|=1} \sum_{|\alpha'|=1} \int_0^1 \mathrm{E}\big[X^\alpha X^{\alpha'} - Y^\alpha Y^{\alpha'}\big] \mathrm{E}\big[\phi\big(V(t;s)\big)(D^{\alpha + \alpha'}h)\big(V(t;s)\big)\big] dt \nonumber\\
	&\quad{} = \frac{s}{2}\int_0^1 \mathrm{E}\left[\phi\big(V(t;s)\big)\mathrm{tr}\left\{(\Sigma_X - \Sigma_Y)\nabla^2h\big(V(t)\big)\right\} \right]dt \nonumber\\
	&\quad{} \leq \frac{1}{2} \big\| \Sigma_X - \Sigma_Y\big\|_{op} \left\| \left\|\nabla^2 h\right\|_{S_1} \right\|_\infty \mathrm{P}\left(\left\|\sqrt{t} W(s) + \sqrt{1 -t}Y \right\|_p  \in A^{3\delta + \kappa_p}\setminus A \right) \nonumber\\
	&\quad{} \leq \frac{1}{2} \big\| \Sigma_X - \Sigma_Y\big\|_{op} \left\| \sum_{|\alpha| = 2} \left|D^\alpha  h\right| \right\|_\infty \mathrm{P}\left(\left\|\sqrt{t} W(s) + \sqrt{1 -t}Y \right\|_p  \in A^{3\delta + \kappa_p}\setminus A \right),
	\end{align*}
	where $\|\cdot\|_{S_1}$ denotes the Schatten $1$-norm.	Now, proceed as in Step 1 replacing $\Delta_p$ by $\Delta_{op}$.
	
	This concludes the proof of the second statement.
\end{proof}

\subsubsection{Proofs for Appendix~\ref{subsec:SmoothingInequalities}}
\begin{proof}[\textbf{Proof of Lemma~\ref{lemma:SmoothIndicator-C-Infty-Function-Even-p}}]
	\textbf{Step 1. Smooth approximation of an indicator function.} Let $\delta \geq \epsilon > 0$. For $x \in \mathbb{R}$ and $A \in \mathcal{B}(\mathbb{R})$ define
	\begin{align*}
	I_{\delta, A}(x) = \left(1 - \delta^{-1} \inf_{y \in A^{2\delta}} |x-y|\right) \vee 0.
	\end{align*}
	Note that $I_{\delta, A}$ is Lipschitz continuous with Lipschitz constant $\delta^{-1}$. Let $\psi \in C_0^\infty(\mathbb{R})$ be a mollifier with compact support $[-1,1]$, e.g. take
	\begin{align*}
	\psi(t) = \begin{cases}
	c \exp\left(\frac{1}{x^2 - 1}\right) & |x| < 1\\
	0 &	|x| \geq 1,
	\end{cases}
	\end{align*}
	where $c > 0$ is such that $\int \psi(x) dx = 1$, and set
	\begin{align*}
	g_{\delta, \epsilon, A}(x) = \frac{1}{\epsilon} \int_{|x-y| \leq \epsilon}\psi\left(\frac{x-y}{\epsilon}\right)I_{\delta, A}(y) dy.
	\end{align*}
	We observe the following facts: First, $\delta \geq \epsilon$ implies that $g_{\delta, \epsilon, A}(x) =1$ for $x \in A$, and $g_{\delta, \epsilon, A}(x) = 0$ for $x \notin A^{3\delta}$, and hence,
	\begin{align}\label{eq:lemma:SmoothIndicator-C-Infty-Function-1}
	\mathbf{1}_A(x) \leq g_{\delta, \epsilon, A}(x) \leq \mathbf{1}_{A^{3\delta}}(x).
	\end{align}	
	Second, $g_{\delta, \epsilon, A} \in C^\infty_b(\mathbb{R})$ and its derivatives up to order three satisfy
	\begin{align}\label{eq:lemma:SmoothIndicator-C-Infty-Function-2}
	\begin{split}
	|g_{\delta, \epsilon, A}'(x)| &\leq C_0\delta^{-1}\mathbf{1}_{A^{3\delta}\setminus A}(x),\\
	|g_{\delta, \epsilon, A}''(x)| &\leq C_0\epsilon^{-1} \delta^{-1}\mathbf{1}_{A^{3\delta}\setminus A}(x),\\
	|g_{\delta, \epsilon, A}'''(x)| &\leq C_0\epsilon^{-2} \delta^{-1}\mathbf{1}_{A^{3\delta}\setminus A}(x),
	\end{split}
	\end{align}
	where $C_0 > 0$ is an absolute constant that depends only on the mollifier $\psi$.
	
	\textbf{Step 2. Smooth approximation of indicator functions of $\ell_p$-norms for even $p \in 2 \mathbb{N}$.}
	For $\eta > 0$ define
	\begin{align}\label{eq:DefM-2}
	M_{p,\eta}(x) = \left(\eta^p + \sum_{j=1}^d |x_j|^p \right)^{1/p}.
	\end{align}
	Since $\eta > 0$ and $p \in 2 \mathbb{N}$, the map $M_{p,\eta}$ is $C^{\infty}(\mathbb{R}^d)$. Moreover, $M_{p,\eta}(x)$ approximates $\|x\|_p$, i.e.
	\begin{align}\label{eq:lemma:SmoothIndicator-C-Infty-Function-Even-p-Fact-1}
	\|x\|_p \leq M_{p, \eta}(x)\leq \|x\|_p + \eta,
	\end{align}
	and $M_{p, \eta}(x)$ is bounded away from zero, i.e.
	\begin{align}\label{eq:lemma:SmoothIndicator-C-Infty-Function-Even-p-Fact-2}
	\min_{x \in \mathbb{R}^d} M_{p,\eta}(x) \geq\eta > 0.
	\end{align}
	Compose $g_{\delta, \epsilon,A}$ and $M_{p,\eta}$ to obtain a smooth approximation of the map $x \mapsto 1_A(\|x\|_p)$, i.e.
	\begin{align*}
	h_{p, \delta, \epsilon, \eta, A}(x) = \left(g_{\delta, \epsilon, A} \circ M_{p,\eta}\right) (x).
	\end{align*}
	Combine eq.~\eqref{eq:lemma:SmoothIndicator-C-Infty-Function-1} and~\eqref{eq:lemma:SmoothIndicator-C-Infty-Function-Even-p-Fact-1} and take expectation with respect to the law of $X$ to conclude that for $A \in \mathcal{B}(\mathbb{R})$ and $p \in 2 \mathbb{N}$,
	\begin{align*}
	\mathrm{P}\left(\|X\|_p \in A \right) \leq \mathrm{E}\left[h_{p, \delta, \epsilon, \eta,  A^{\eta}}\left(X\right) \right] \leq \mathrm{P}\left(\|X\|_p \in A^{3\delta + 2\eta} \right).
	\end{align*}
	
	\textbf{Step 3. Bounds on partial derivatives of $h_{p, \delta, \epsilon, \eta, A}$ in transformed conjugate norm.} Let $\tau \in [1, \infty]$ and set $q=\frac{\tau p}{\tau p-1}$. By Fa{\`a} di Bruno's chain rule and H{\"o}lder's inequality,
	\begin{align}
	\begin{split}\label{eq:lemma:SmoothIndicator-C-Infty-Function-3}
	&\left(\sum_{|\alpha| = 2}\left|D^\alpha h_{p, \delta, \epsilon, \eta, A}(x)\right|^q\right)^{1/q}\\
	&\hspace{20pt}\leq \|g_{\delta, \epsilon, A}'' \circ M_{p, \eta}\|_\infty \left(\sum_{|\alpha| = 1}\sum_{|\alpha'| = 1}\left|  D^\alpha M_{p, \eta}(x)\right|^q \left|  D^\alpha M_{p, \eta}(x) \right|^q\right)^{1/q}\\
	&\quad{} \hspace{20pt} + \|g_{\delta, \epsilon, A}' \circ M_{p, \eta}\|_\infty  \left(\sum_{|\alpha| = 2}\left|D^\alpha M_{p, \eta}(x)\right|^q \right)^{1/q},
	\end{split}\\
	\begin{split}\label{eq:lemma:SmoothIndicator-C-Infty-Function-4}
	&\left(\sum_{|\alpha| = 3}\left|D^\alpha h_{p, \delta, \epsilon, \eta, A}(x) \right|^q\right)^{1/q}\hspace{-10pt}\\
	&\hspace{20pt}\leq \|g_{\delta, \epsilon, A}''' \circ M_{p, \eta}\|_\infty \left(\sum_{|\alpha| = 1} \sum_{|\alpha'| = 1} \sum_{|\alpha''| = 1}\left| D^\alpha M_{p, \eta}(x) \right|^q \left|  D^{\alpha'} M_{p, \eta}(x) \right|^q \left|  D^{\alpha''} M_{p, \eta}(x)\right|^q \right)^{1/q}\\
	&\quad{}\hspace{20pt} + 3 \|g_{\delta, \epsilon, A}'' \circ M_{p, \eta}\|_\infty \left(\sum_{|\alpha| = 1}\sum_{|\alpha'| = 2}\left| D^\alpha M_{p, \eta}(x) \right|^q \left|D^{\alpha'} M_{p, \eta}(x) \right|^q\right)^{1/q} \\
	&\quad{}\hspace{20pt} + \|g_{\delta, \epsilon, A}' \circ M_{p, \eta}\|_\infty \left( \sum_{|\alpha| = 3}\left| D^\alpha M_{p, \eta}(x) \right|^q\right)^{1/q}.
	\end{split}
	\end{align}
	
	Next, we bound the partial derivatives of $M_{p, \eta}$. By Lemma~\ref{lemma:Lp-Norm-Derivatives-Stability-3}, the chain rule
	\begin{align}
	\begin{split}\label{eq:lemma:SmoothIndicator-C-Infty-Function-Bound-SecondDiff-3}
	&\left(\sum_{|\alpha| = 1}\sum_{|\alpha'| = 1}\left|  D^\alpha M_{p, \eta}(x) \right|^q \left| D^{\alpha'} M_{p, \eta}(x)\right|^q\right)^{1/q}\\
	&\hspace{20pt} = \left(\sum_{k, \ell} \left|\frac{\partial M_{p, \eta}(x)}{\partial x_k}\right|^q \left|\frac{\partial M_{p,\eta}(x)}{\partial x_\ell}\right|^q\right)^{1/q} =  \left(\sum_{k=1}^d \left|\frac{\partial M_{p, \eta}(x)}{\partial x_k}\right|^q \right)^{2/q}\leq d^{2(\tau-1)/(\tau p)},
	\end{split}\\
	\begin{split}\label{eq:lemma:SmoothIndicator-C-Infty-Function-Bound-ThirdDiff-3}
	&\left(\sum_{|\alpha| = 1} \sum_{|\alpha'| = 1} \sum_{|\alpha''| = 1}\left| D^\alpha M_{p, \eta}(x) \right|^q \left|  D^{\alpha'} M_{p, \eta}(x) \right|^q \left|  D^{\alpha''} M_{p, \eta}(x)\right|^q\right)^{1/q}\\
	&\hspace{20pt}= \left(\sum_{k, \ell, m} \left|\frac{\partial M_{p, \eta}(x)}{\partial x_k}\right|^q \left|\frac{\partial M_{p, \eta}(x)}{\partial x_\ell}\right|^q \left|\frac{\partial M_{p, \eta}(x)}{\partial x_m}\right|^q\right)^{1/q}\\
	&\hspace{20pt} = \left(\sum_{k=1}^d \left|\frac{\partial M_{p, \eta}(x)}{\partial x_k}\right|^q \right)^{3/q} \leq d^{3(\tau-1)/(\tau p)},
	\end{split}
	\end{align}
	and
	\begin{align}\label{eq:lemma:SmoothIndicator-C-Infty-Function-Bound-SecondDiff-4}
	\begin{split}
	&\left(\sum_{|\alpha| = 2}\left|D^\alpha M_{p, \eta}(x) \right|^q\right)^{1/q}\\
	&\hspace{20pt}\leq \left(\sum_{k=1}^d \left|\frac{\partial^2 M_{p, \eta}(x)}{\partial x_k^2}\right|^q\right)^{1/q} + \left(\sum_{k, \ell} \left|\frac{\partial^2 M_{p, \eta}(x)}{\partial x_k \partial x_\ell}\right|^q\right)^{1/q} \\
	&\hspace{20pt}\leq \frac{2(p-1)d^{(2\tau-1)/(\tau p)}}{M_{p, \eta}(x)} + \frac{2(p-1)d^{(\tau-1)/(\tau p)}}{M_{p, \eta}(x)} + \frac{(p-1)d^{2(\tau-1)/(\tau p)}}{M_{p, \eta}(x)} \\
	&\hspace{20pt}\lesssim \eta^{-1}pd^{(2\tau-1)/(\tau p)},
	\end{split}
	\end{align}	
	where the third inequality follows from the lower bound~\eqref{eq:lemma:SmoothIndicator-C-Infty-Function-Even-p-Fact-2};  and
	\begin{align}\label{eq:lemma:SmoothIndicator-C-Infty-Function-Bound-ThirdDiff-4}
	\begin{split}
	&\left(\sum_{|\alpha| = 1}\sum_{|\alpha'| = 2}\left| D^\alpha M_{p, \eta}(x) \right|^q \left|D^{\alpha'} M_{p, \eta}(x)\right|^q\right)^{1/q} \\
	&\hspace{20pt}\leq \left(\sum_{k, \ell} \left|\frac{\partial^2 M_{p, \eta}(x)}{\partial x_k^2}\right|^q \left|\frac{\partial M_{p, \eta}(x)}{\partial x_\ell} \right|^q\right)^{1/q} + \left(\sum_{k, \ell, m} \left|\frac{\partial^2 M_{p, \eta}(x)}{\partial x_k \partial x_\ell}\right|^q \left|\frac{\partial M_{p, \eta}(x)}{\partial x_m} \right|^q\right)^{1/q} \\
	&\hspace{20pt} \lesssim  \eta^{-1}p d^{(3\tau-2)/(\tau p)},
	\end{split}
	\end{align}
	where we have used the results from eq.~\eqref{eq:lemma:SmoothIndicator-C-Infty-Function-Bound-SecondDiff-3} and~\eqref{eq:lemma:SmoothIndicator-C-Infty-Function-Bound-SecondDiff-4}; and
	\begin{align}\label{eq:lemma:SmoothIndicator-C-Infty-Function-Bound-ThirdDiff-5}
	&\left(\sum_{|\alpha| = 3}\left| D^\alpha M_{p, \eta}(x)\right|^q\right)^{1/q} \nonumber\\
	&\hspace{20pt}\leq \left(\sum_{k} \left|\frac{\partial^3 M_{p, \eta}(x)}{\partial x_k^3}\right|^q\right)^{1/q}+  \left(\sum_{k, \ell} \left|\frac{\partial^3 M_{p, \eta}(x)}{\partial x_k^2 \partial x_\ell}\right|^q\right)^{1/q} + \left(\sum_{k, \ell, m} \left|\frac{\partial^3 M_{p, \eta}(x)}{\partial x_k \partial x_\ell \partial x_m}\right|^q\right)^{1/q}\nonumber\\
	&\hspace{20pt}\leq \frac{4(p-1)(p-2)d^{(3\tau-1)/(\tau p)}}{ M_{p, \eta}^2(x)} + \frac{12(p-1)d^{(\tau-1)/(\tau p)}}{M_{p, \eta}^2(x)} + \frac{4(2p-1)(p-1)d^{(\tau-1)/(\tau p)}}{ M_{p, \eta}^2(x)}  \nonumber\\
	&\hspace{20pt}\quad{}  + \frac{2(p-1)^2 d^{(3\tau-2)/(\tau p)}}{ M_{p, \eta}^2(x)} + \frac{2(2p-1)(p-1)d^{2(\tau-1)/(\tau p)}}{ M_{p, \eta}^2(x)} + \frac{(2p-1)(p-1)d^{3(\tau-1)/(\tau p)}}{ M_{p, \eta}^2(x)} \nonumber\\
	&\hspace{20pt}\lesssim \eta^{-2} p^2 d^{(3\tau-1)/(\tau p)} ,
	\end{align}
	
	To conclude, set $\delta = \epsilon > 0$ and $\eta = \beta^{-1}p d^{1/(\tau p)}$, $\beta > 0$. Then, the upper bounds~\eqref{eq:lemma:SmoothIndicator-C-Infty-Function-2} and~\eqref{eq:lemma:SmoothIndicator-C-Infty-Function-3}--\eqref{eq:lemma:SmoothIndicator-C-Infty-Function-Bound-ThirdDiff-5} imply, uniformly in $x \in \mathbb{R}^d$ and $A\in \mathcal{A}$,
	\begin{align*}
	&\left(\sum_{|\alpha| = 2} \left|D^\alpha h_{p, \delta, \epsilon, \eta, A}(x)\right|^q\right)^{1/q} \lesssim \left(\delta^{-2} + \delta^{-1}\beta \right) d^{2(\tau-1)/(\tau p)}, \\
	&\left(\sum_{|\alpha| = 3} \left|D^\alpha h_{p, \delta, \epsilon, \eta A}(x) \right|^q\right)^{1/q} \lesssim \left(\delta^{-3}+ \delta^{-2}\beta + \delta^{-1}\beta^2 \right) d^{3(\tau-1)/(\tau p)}.
	\end{align*}
	Note that due to the substitutions $h_{p,  \delta, \epsilon, \eta, A}$ depends only on $p, d, \beta, \delta, A$.
	
	\textbf{Step 4. Smooth approximation of indicator function of $\ell_p$-norms with $p \in [1, \infty)$.}  Let $\mathcal{I} = \{A \subseteq \mathbb{R} : A = [0,t], t \geq 0\}$. Let $p \in [1, \infty)$ be arbitrary and define $p_+ = 2\lceil \frac{p}{2}\rceil$ to be the smallest even integer larger than (or equal to) $p$. Then, $\|x\|_{p_+}\leq \|x\|_p$ for all $x \in \mathbb{R}^d$. We have the following relation between $M_{p_+,\beta}(x)$ and $\|x\|_p$
	\begin{align}\label{eq:lemma:SmoothIndicator-C-Infty-Function-5}
	\|x\|_{p_+} \leq M_{p_+,\eta}(x)\leq \|x\|_{p_+} +\eta \leq \|x\|_p + \eta.
	\end{align}
	Combine eq.~\eqref{eq:lemma:SmoothIndicator-C-Infty-Function-1} and~\eqref{eq:lemma:SmoothIndicator-C-Infty-Function-5} and take expectation with respect to the law of $X$ to conclude that for $A \in \mathcal{I}$ and $p \in [1, \infty)$,
	\begin{align*}
	\mathrm{P}\left(\|X\|_p \in A \right) \leq \mathrm{P}\left(\|X\|_{p_+} \in A \right) \leq \mathrm{E}\left[h_{p_+, \beta, \delta, \epsilon, A^{\eta}}\left(X\right) \right] \leq \mathrm{P}\left(\|X\|_p \in A^{3\delta + 2\eta} \right),
	\end{align*}
	where the first inequality follows from the fact that $0 \leq \|x\|_{p_+}\leq \|x\|_p$ for all $x \in \mathbb{R}^d$ and the fact that $A = [0, t]$ or some $t \geq 0$.
\end{proof}
\begin{proof}[\textbf{Proof of Lemma~\ref{lemma:SmoothIndicator-C-Infty-Function-Large-p}}]
	For $p \geq \log d$ we can approximate any $\ell_p$-norm by the smooth max function. We can therefore sharpen the result from Lemma~\ref{lemma:SmoothIndicator-C-Infty-Function-Even-p}.
	
	\textbf{Step 1. Smooth approximation of indicator functions of $\ell_p$-norms with $p \geq \log d$.} Let $\delta \geq \epsilon > 0$,  $\mathcal{A} = \{A \subseteq \mathbb{R} : A = [0,t], t \geq 0 \}$.  For $x \in \mathbb{R}$ and $A \in \mathcal{A}$ define
	\begin{align*}
	g_{\delta, \epsilon, A}(x) = \frac{1}{\epsilon} \int_{|x-y| \leq \epsilon}\psi\left(\frac{x-y}{\epsilon}\right)I_{\delta, A}(y) dy,
	\end{align*}
	with $I_{\delta, A}$ and $\psi$ as in the proof of Lemma~\ref{lemma:SmoothIndicator-C-Infty-Function-Even-p}. Recall that
	\begin{align}\label{eq:lemma:SmoothIndicator-C-Infty-Function-Large-p-1}
	\mathbf{1}_A(x) \leq g_{\delta, \epsilon, A}(x) \leq \mathbf{1}_{A^{3\delta}}(x).
	\end{align}
	For $\beta > 1$ define the smooth max function
	\begin{align}\label{eq:DefF}
	F_{\beta}(x) := \beta^{-1} \log \left(\sum_{k=1}^p e^{\beta x_kd^{-1/p}} + e^{-\beta x_k d^{-1/p}}\right).
	\end{align}
	Let $x \in \mathbb{R}^d$ be arbitrary. Set $u^* = \arg\max_{\|u\|_q =1}|x'u|$, where $q = p/(p-1)$ is the conjugate exponent to $p$. Note that $d^{-1/p}\|u^*\|_1 \leq 1$ and $1 \leq d^{1/p} \leq e$. Therefore, we have, for $\beta > 0$,
	\begin{align}\label{eq:lemma:SmoothIndicator-C-Infty-Function-Large-p-2}
	\begin{split}
	\|x\|_p &= \beta^{-1} d^{1/p}\sum_{k=1}^d \frac{u_k^*}{d^{1/p}} \log\left(e^{|x_k|\beta}\right)\\
	&\leq \beta^{-1} d^{1/p} \log \left(\sum_{k=1}^d \frac{u_k^*}{d^{1/p}} e^{|x_k|\beta} \right)\\
	&\leq \beta^{-1} d^{1/p} \log \left(\sum_{k=1}^d e^{|x_k|\beta d^{-1/p}} \right)\\
	&\leq \beta^{-1} d^{1/p} \log \left(\sum_{k=1}^d e^{x_k\beta d^{-1/p}} + e^{-x_k\beta d^{-1/p}} \right)\\
	&\leq \|x\|_\infty + d^{1/p} \beta^{-1} \log (2d)\\
	&\leq \|x\|_p + e\beta^{-1} \log (2d),
	\end{split}
	\end{align}
	where the first inequality follows from Jensen's inequality, and the second and third inequalities from elementary calculations.
	
	We define
	\begin{align*}
	h_{p, \beta, \delta, \epsilon, A}(x) = \left(g_{\delta, \epsilon, A} \circ F_\beta \right) (x).
	\end{align*}
	Combine eq.~\eqref{eq:lemma:SmoothIndicator-C-Infty-Function-Large-p-1} and~\eqref{eq:lemma:SmoothIndicator-C-Infty-Function-Large-p-2} and take expectation with respect to the law of $X$ to conclude that
	\begin{align*}
	\mathrm{P}\left(\|X\|_p \in A \right) \leq \mathrm{E}\left[h_{p, \beta, \delta, \epsilon, A^{e\beta^{-1}\log (2d)}}\left(X\right) \right] \leq \mathrm{P}\left(\|X\|_p \in A^{3\delta + 2e\beta^{-1} \log(2d) } \right).
	\end{align*}
	
	\textbf{Step 2. Bounds on partial derivatives of $h_{p, \beta, \delta, \epsilon, A}$ in $\ell_1$-norm.} By Lemma A.2--A.6 in~\cite{chernozhukov2013GaussianApproxVec} we have
	\begin{align*}
	\sup_{A \in \mathcal{A}}\left\|\sum_{|\alpha| = 2} \left|D^\alpha h_{p, d, \beta, \delta, A}\right|\right\|_\infty &\lesssim \frac{1}{\delta^2} + \frac{\beta}{\delta},\\
	\sup_{A \in \mathcal{A}}\left\|\sum_{|\alpha| = 3} \left|D^\alpha  h_{p, d, \beta, \delta, A} \right|\right\|_\infty &\lesssim \frac{1}{\delta^3}+ \frac{\beta}{\delta^2} + \frac{\beta^2}{\delta}.
	\end{align*}
	This concludes the proof.
\end{proof}

\subsubsection{Proofs for Appendix~\ref{subsec:AuxResults}}
\begin{proof}[\textbf{Proof of Lemma~\ref{lemma:MomentsOfLpNorm}}]
	For $t \geq p$, it follows from Minkowski's integral inequality that
	\begin{align*}
	\left(\mathrm{E}\|X\|^t_p\right)^{p/t} \leq \sum_{k=1}^d \left(\mathrm{E}|X_k|^t\right)^{p/t} \lesssim K_s^p \|\sigma\|_p^p.
	\end{align*}
	While for $t \leq p$, H{\"o}lder's inequality yields
	\begin{align*}
	\left(\mathrm{E}\|X\|^t_p\right)^{p/t} \leq \mathrm{E}\left[\sum_{k=1}^d |X_k|^p\right] = \sum_{k=1}^d \mathrm{E}|X_k|^p \lesssim K_s^p \|\sigma\|_p^p.
	\end{align*}
	Combine both inequalities to conclude.
\end{proof}	
\begin{proof}[\textbf{Proof of Lemma~\ref{lemma:ProductSubgaussian}}]
	Recall Young's inequality: $\prod_{i=1}^K x_i^{\alpha_i} \leq \sum_{i=1}^K \alpha_i x_i$ for all $x_i, \alpha_i \geq 0$, $i=1, \ldots, K$ with $\sum_{i=1}^K \alpha_i = 1$. Without loss of generality, we can assume that $\|X_i\|_{\psi_2} = 1$ for all $i=1, \ldots, K$. Thus, the claim of the lemma follows if we can show the following: If $\mathrm{E}\left[\exp(X_i^2)\right] \leq 2$ for all $i=1, \ldots, K$, then $\mathrm{E}[\exp(\prod_{i=1}^K|X_i|^{2/K} )] \leq 2$. This assertion follows from straightforward calculations:
	\begin{align*}
	\mathrm{E}\left[\psi_{2/K}\left(\prod_{i=1}^KX_i\right)\right] &= \mathrm{E}\left[\exp\left(\prod_{i=1}^K |X_i|^{2/K}\right)\right]\\
	&\leq \mathrm{E}\left[\exp\left(\frac{1}{K}\sum_{i=1}^K |X_i|^2\right)\right]\\
	&=\mathrm{E}\left[\prod_{i=1}^K\exp\left(\frac{1}{K} |X_i|^2\right)\right]\\
	&\leq \frac{1}{K} \left(\sum_{i=1}^K\mathrm{E}\left[\exp\left(|X_i|^2\right)\right] \right)\\
	&\leq 2,
	\end{align*}
	where in the first and second inequalities we have used Young's inequality.
\end{proof}
\begin{proof}[\textbf{Proof of Lemma~\ref{lemma:BoundsCovariance}}]	
	\textbf{Proof of Case (i).}
	We have the following:
	\begin{align*}
	\widehat{\Sigma} - \Sigma &= \frac{1}{n} \sum_{i=1}^n \big(X_iX_i' - \mathrm{E}[X_i X_i']\big) - \frac{1}{n^2}\sum_{1\leq i,j \leq n } X_i X_j'\\
	&\equiv \mathbf{I} - \mathbf{II}.
	\end{align*}	
	For $1 \leq j,k \leq d$, set $\widehat{\Sigma}_{jk} = n^{-1}\sum_{i=1}^n X_{ij}X_{ik}$, $\Sigma_{jk} = n^{-1}\sum_{i=1}^n \mathrm{E}[X_{ij}X_{ik}]$, and $\sigma_k^2 =n^{-1}\sum_{i=1}^n \mathrm{E}[X_{ik}^2]$. By Assumption 1 there exists an absolute constant $K_0 > 1$ such that for all $1 \leq i \leq n$,
	\begin{align*}
	\big\|X_{ij} X_{ik}\big\|_{\psi_1} \leq \|X_{ij}\|_{\psi_2}\|X_{ik}\|_{\psi_2} \leq K_0^2 \mathrm{E}[X_{ij}^2]^{1/2}\mathrm{E}[X_{ik}^2]^{1/2} = K_0^2\sigma_j\sigma_k.
	\end{align*}
	Hence, by union bound and Bernstein's inequality there exists an absolute constant $C > 0$ such that for all $t > 0$,
	\begin{align*}
	\mathrm{P}\Big(\|\mathrm{vec}(\mathbf{I})\|_p > t  K_0^2 \|\sigma\|_p^2\Big) \leq \sum_{1 \leq j,k \leq d} \mathrm{P}\left(\big|\widehat{\Sigma}_{jk} - \Sigma_{jk}\big|> t  K_0^2 \sigma_j\sigma_k\right) \leq 2 d^2 \exp\left(- C \min\left\{t^2, t\right\}n\right).
	\end{align*}
	Above tail bound implies that with probability at least $1- \zeta$,
	\begin{align*}
	\|\mathrm{vec}(\mathbf{I})\|_p &\leq K_0^2 \|\sigma\|_p^2\left( \sqrt{\frac{2}{C}}\sqrt{\frac{\log d + \log (2/\zeta)}{n}} \bigvee \frac{2}{C} \frac{\log d + \log (2/\zeta)}{n} \right).
	\end{align*}
	To bound $\|\mathrm{vec}(\mathbf{II})\|_p$ we directly use the sub-gaussianity of the $X_i$'s. By union bound and Hoeffding's inequality there exists an absolute constant $C > 0$ such that for all $t > 0$,
	\begin{align*}
	\mathrm{P}\big(\|\mathrm{vec}(\mathbf{II})\|_p > t^2K_0^2 \|\sigma\|_p^2\big) &\leq \sum_{1 \leq k \leq d} \mathrm{P}\left(|X_{ik}| > tK_0 \sigma_k n\right) \leq 2d e \exp\left(- Ct^2n\right).
	\end{align*}
	Hence, with probability at least $1 - \zeta$,
	\begin{align*}
	\|\mathrm{vec}(\mathbf{II})\|_p &\leq K_0^2 \|\sigma\|_p^2 \frac{2}{C}\left(\frac{\log d + \log (2/\zeta)}{n}\right).
	\end{align*}
	Conclude that with probability at least $1 - 2 \zeta$,
	\begin{align}\label{eq:lemma:BoundsCovariance-1}
	\begin{split}
	\|\mathrm{vec}(\widehat{\Sigma} - \Sigma)\|_p &\lesssim \|\mathrm{vec}(\mathbf{I})\|_p \vee \|\mathrm{vec}(\mathbf{II})\|_p\\ &\lesssim  K_0^2\|\sigma\|_p^2\left(\sqrt{\frac{\log d + \log (2/\zeta)}{n}} \vee \frac{\log d + \log (2/\zeta)}{n} \right).
	\end{split}
	\end{align}	
	
	\textbf{Proof of Case (ii).} We have
	\begin{align*}
	\|\widehat{\Sigma} - \Sigma\|_{op} &\leq \left\|\frac{1}{n} \sum_{i=1}^n \big(X_iX_i' - \mathrm{E}[X_i X_i']\big)\right\|_{op} + \left\|\frac{1}{n^2}\sum_{1\leq i,j \leq n } X_i X_j'\right\|_{op}\\
	&=  \left\|\frac{1}{n} \sum_{i=1}^n \big(X_iX_i' - \mathrm{E}[X_i X_i']\big)\right\|_{op} + \left\|\frac{1}{n}\sum_{i=1}^n X_i \right\|_2^2\\
	&\equiv \|\mathbf{I}\|_{op} + \|\mathbf{II}\|_2^2.
	\end{align*}
	By Theorem 5.39 and Remarks 5.40 and 5.53 in~\cite{vershynin2011introduction}, with probability at least $1- \zeta$,
	\begin{align*}
	\|\mathbf{I}\|_{op} \lesssim \|\Sigma\|_{op} \left(\sqrt{\frac{\mathrm{r}(\Sigma)\log d + \log (2/\zeta)}{n}} \vee \frac{\mathrm{r}(\Sigma)\log d + \log (2/\zeta)}{n} \right),
	\end{align*}
	Similar arguments as in the second part of Case 1 with $p =2$ and the fact that $\|\sigma\|_2^2 = \mathrm{tr}(\Sigma)$, yield, with probability at least $1- \zeta$,
	\begin{align*}
	\|\mathbf{II}\|_2^2 \lesssim \|\Sigma\|_{op} \left(\frac{\mathrm{r}(\Sigma)\log d + \log (2/\zeta)}{n} \right).
	\end{align*}
	The claim follows from combining both bounds.
	
	\textbf{Proof of Case (iii).} Denote by $\oslash$ the Hadamard division and observe that
	\begin{align*}
	\max_{1 \leq k \leq d} \big|(\widehat{\sigma}_k/\sigma_k)^2 - 1\big| &= \left\|\mathrm{diag}(\widehat{\Sigma})\oslash \mathrm{diag}(\Sigma) - I_d\right\|_{\infty}\\
	&\leq \left\|\mathrm{diag}\left( \left(\frac{1}{n} \sum_{i=1}^n X_iX_i' - \mathrm{E}[X_i X_i']\right)\oslash\mathrm{diag}(\Sigma) - I_d\right)\right\|_\infty\\
	&\quad{} + \left\| \mathrm{diag}\left( \left(\frac{1}{n^2}\sum_{1\leq i,j \leq d } X_i X_j'\right) \oslash \mathrm{diag}(\Sigma) \right)\right\|_\infty.
	\end{align*}
	Moreover, $\|\mathrm{diag}(I_d)\|_{\infty} = 1$. Hence, the claim follows from Case 1 with $p =\infty$.
\end{proof}
\begin{proof}[\textbf{Proof of Lemma~\ref{lemma:BoundsCovariance-FiniteMoments}}]
	\textbf{Proof of Case (i).}
	We have the following:
	\begin{align*}
	\widehat{\Sigma} - \Sigma &= \frac{1}{n} \sum_{i=1}^n \big(X_iX_i' - \mathrm{E}[X_i X_i']\big) - \frac{1}{n^2}\sum_{1\leq i,j \leq d } X_i X_j'\\
	&\equiv \mathbf{I} - \mathbf{II}.
	\end{align*}		
	\textbf{Step 1.} We begin with the analysis of $\mathbf{I}$. Let $X \in \mathbb{R}^d$ satisfy Assumption~\ref{assumption:FiniteMoments}. Then, there exists an absolute constant $C_1 > 0$ such that
	\begin{align}\label{eq:lemma:BoundsCovariance-HeavyTailed-1}
	\mathrm{E}\left[\|\mathrm{vec}(XX') \|_p^2\right] = \mathrm{E}\left[ \left(\sum_{j=1}^d \sum_{k=1}^d |X_jX_k|^p\right)^{2/p}\right]=  \mathrm{E}\left[ \left(\sum_{j=1}^d |X_j|^p\right)^{4/p}\right] \leq C_1 K_{p \vee 4}^4 \|\sigma\|_p^4.
	\end{align}
	Since $\mathrm{E}[X] = 0$, eq.~\eqref{eq:lemma:BoundsCovariance-HeavyTailed-1} implies that
	\begin{align*}
	\sup_{\|u\|_q=1} \mathrm{Var}\left[\mathrm{vec}(XX')'u\right] = \sup_{\|u\|_q=1} \mathrm{E}\left[\big(\mathrm{vec}(XX')'u\big)^2\right] \leq \mathrm{E}\left[\|\mathrm{vec}(XX') \|_p^2\right] \leq C_1 K_{p \vee 4}^4\|\sigma\|_p^4.
	\end{align*}
	(As an aside, we note the following: This bound is obviously loose; a tighter upper bound would involve the operator norm of the population covariance matrix. However, we will see later that this tighter bound would not improve the final rate of convergence.)
	
	Therefore, by Symmetrization Lemma 2.3.7 in~\cite{vandervaart1996weak}, for any $t \geq 8n^{-1/2} C_1^{1/2} K_{p \vee 4}^2 \|\sigma\|_p^2$,
	\begin{align}\label{eq:lemma:BoundsCovariance-HeavyTailed-2}
	\mathrm{P}\left(\left\|\mathrm{vec}\left(\frac{1}{n}\sum_{i=1}^nX_iX_i' - \mathrm{E}[X_iX_i'] \right) \right\|_p > t\right) \leq 4 \mathrm{P}\left(\left\|\frac{1}{n}\sum_{i=1}^n\mathrm{vec}(X_iX_i')\varepsilon_i \right\|_p > t/4\right),
	\end{align}
	where $\varepsilon_1, \ldots, \varepsilon_n$ are i.i.d. Rademacher random variables independent of the $X_i$'s.	Let $\theta > 0$ (to be specified below) and define
	\begin{align*}
	A(\theta) := \left\{ \omega \in \Omega: \left(\mathrm{E}\left[\left\|\frac{1}{n}\sum_{i=1}^n\mathrm{vec}(X_iX_i')\varepsilon_i \right\|_p^2 \mid X_1, \ldots, X_n\right]\right)(\omega) \leq \theta \right\}.
	\end{align*}
	Expand the tail probability on the right hand side in above eq.~\eqref{eq:lemma:BoundsCovariance-HeavyTailed-2},
	\begin{align}\label{eq:lemma:BoundsCovariance-HeavyTailed-3}
	&\mathrm{P}\left(\left\|\frac{1}{n}\sum_{i=1}^n\mathrm{vec}(X_iX_i')\varepsilon_i \right\|_p > t/4, \:\: \mathrm{E}\left[\left\|\frac{1}{n}\sum_{i=1}^n\mathrm{vec}(X_iX_i')\varepsilon_i \right\|_p^2 \mid X_1, \ldots, X_n\right] \leq \theta \right)\nonumber\\
	&\hspace{40pt} + \mathrm{P}\left( \mathrm{E}\left[\left\|\frac{1}{n}\sum_{i=1}^n\mathrm{vec}(X_iX_i')\varepsilon_i \right\|_p^2 \mid X_1, \ldots, X_n\right] > \theta \right)\nonumber\\
	&\leq 2\int_{A(\theta)} \exp\left(- \frac{t^2}{512 \mathrm{E}\left[\left\|\frac{1}{n}\sum_{i=1}^n\mathrm{vec}(X_iX_i')\varepsilon_i \right\|_p^2 \mid X_1, \ldots, X_n\right]}\right)d\mathrm{P}_{X_1, \ldots, X_n}(\omega)\nonumber\\
	&\hspace{40pt} + \mathrm{P}\left( \mathrm{E}\left[\left\|\frac{1}{n}\sum_{i=1}^n\mathrm{vec}(X_iX_i')\varepsilon_i \right\|_p^2 \mid X_1, \ldots, X_n\right] > \theta \right)\nonumber\\
	&\leq 2 \exp\left(-\frac{t^2}{512 \theta }\right) + \theta^{-1} \mathrm{E}\left[\left\|\frac{1}{n}\sum_{i=1}^n\mathrm{vec}(X_iX_i')\varepsilon_i \right\|_p^2\right],
	\end{align}
	where the first inequality follows from the sub-gaussianity of Rademacher random variables~\citep[e.g.][Theorem 4.7 and eq. (4.12) on p. 101]{ledoux1991probability} and the second inequality by Markov's inequality.
	
	We now determine the choice of $\theta > 0$. By Theorem 2.2 in~\cite{duembgen2010nemirovskis} (refinement of Nemirovski's inequality) there exists an absolute constant $C_2 > 0$ such that
	\begin{align}\label{eq:lemma:BoundsCovariance-HeavyTailed-4}
	\begin{split}
	\mathrm{E}\left[ \left\|\sum_{i=1}^n\mathrm{vec}\big(X_iX_i'\big)\varepsilon_i \right\|_p^2 \right] &\leq  C_2\big(p \wedge \log d \big) \sum_{i=1}^n \mathrm{E}\left[\|\mathrm{vec}\big(X_iX_i'\big)\varepsilon_i\|_p^2 \right]\\
	&\leq C_2 \big(p \wedge \log d \big) n \mathrm{E}\left[\|\mathrm{vec}\big(XX'\big)\|_p^2 \right].
	\end{split}
	\end{align}
	Combine eq.~\eqref{eq:lemma:BoundsCovariance-HeavyTailed-1} and eq.~\eqref{eq:lemma:BoundsCovariance-HeavyTailed-4} to conclude that
	\begin{align*}
	\mathrm{E}\left[ \left\|\frac{1}{n}\sum_{i=1}^n\mathrm{vec}\big(X_iX_i'\big)\varepsilon_i\right\|_p^2 \right] \leq C_1 C_2 K_{p \vee 4}^4\|\sigma\|_p^4 \left(\frac{ p \wedge \log d}{n}\right).
	\end{align*}
	Thus, we set $\theta = M (C_1 C_2 K_{p \vee 4}^4 \vee 64) \|\sigma\|_p^4 \left(\frac{ p \wedge \log d}{n}\right)$ and $t = M^{1/2}\theta^{1/2} \geq 8 n^{-1/2}C_1^{1/2} K_{p \vee 4}^2 \|\sigma\|_p^2$, where $M \geq 1$ is a large absolute constant. By choosing $M$ large enough we can make the left hand side of eq.~\eqref{eq:lemma:BoundsCovariance-HeavyTailed-3} arbitrarily small. Hence, we conclude that
	\begin{align}\label{eq:lemma:BoundsCovariance-HeavyTailed-5}
	\left\|\mathrm{vec}\left(\mathbf{I}\right) \right\|_p = O_p\left( K_{p \vee 4}^2 \|\sigma\|_p^2 \sqrt{\frac{ p \wedge \log d}{n}}\right).
	\end{align}
	
	\textbf{Step 2.} We now analyze term $\mathbf{II}$. Let $X \in \mathbb{R}^d$ satisfy Assumption~\ref{assumption:FiniteMoments} and let $\widetilde{X}$ be an independent copy of $X$. Then, there exists an absolute constant $C_1 > 0$ such that
	\begin{align}\label{eq:lemma:BoundsCovariance-HeavyTailed-6}
	\mathrm{E}\left[\|\mathrm{vec}(X\widetilde{X}') \|_p^2\right] = \mathrm{E}\left[ \left(\sum_{j=1}^d \sum_{k=1}^d |X_j\widetilde{X}_k|^p\right)^{2/p}\right]=  \mathrm{E}\left[ \left(\sum_{j=1}^d |X_j|^p\right)^{2/p}\right]^2 \leq C K_p^4 \|\sigma\|_p^4.
	\end{align}
	
	Let $1 < p,q < \infty$ be conjugate exponents such that $1/p + 1/q =1$. By standard decoupling arguments~\citep[e.g.][Theorem 8.11]{foucart2013mathematical} we have
	\begin{align}\label{eq:lemma:BoundsCovariance-HeavyTailed-7}
	\sup_{\|u\|_q=1} \mathrm{E}\left[\left(\sum_{i\neq j} \mathrm{vec}(X_i X_j')'u \right)^2\right] \leq 	16 \sup_{\|u\|_q=1} \mathrm{E}\left[\left(\sum_{1 \leq i,j \leq n} \mathrm{vec}(X_i \widetilde{X}_j')'u \right)^2\right],
	\end{align}
	where $\widetilde{X}_1, \ldots, \widetilde{X}_n$ are mutually independent copies of the corresponding $X_i$'s. Since $\mathrm{E}[X] = \mathrm{E}[\widetilde{X}] = 0$, we can further bound the right hand side of above inequality using eq.~\eqref{eq:lemma:BoundsCovariance-HeavyTailed-6},
	\begin{align}\label{eq:lemma:BoundsCovariance-HeavyTailed-8}
	16 \sup_{\|u\|_q=1} \mathrm{E}\left[\sum_{1 \leq i,j \leq n} \left( \mathrm{vec}(X_i \widetilde{X}_j')'u \right)^2\right] \leq 16 \mathrm{E}\left[\sum_{1 \leq i,j \leq n} \left\|\mathrm{vec}(X_i \widetilde{X}_j')\right\|_p^2\right] \leq 16 n^2 C_4 K_p^4\|\sigma\|_p^4.
	\end{align}
	Therefore, by Symmetrization Lemma 2.3.7 in~\cite{vandervaart1996weak}, for any $t \geq 32 n^{-1}C_4^{1/2} K_p^2 \|\sigma\|_p^2$,
	\begin{align}\label{eq:lemma:BoundsCovariance-HeavyTailed-9}
	\mathrm{P}\left(\left\|\mathrm{vec}\left(\frac{1}{n^2}\sum_{i\neq j} X_i X_j'\right) \right\|_p > t\right) \leq 4 \mathrm{P}\left(\left\|\mathrm{vec}\left(\frac{1}{n^2}\sum_{i\neq j} X_i X_j'\right)\varepsilon_{ij} \right\|_p > t/4\right),
	\end{align}
	where $\varepsilon_{11}, \ldots, \varepsilon_{nn}$ are i.i.d. Rademacher random variables independent of the $X_iX_j'$'s.	Proceeding as in Step 1, we upper bound the tail probability in~\eqref{eq:lemma:BoundsCovariance-HeavyTailed-9} by
	\begin{align}\label{eq:lemma:BoundsCovariance-HeavyTailed-10}
	8 \exp\left(-\frac{t^2}{512 \theta }\right) + 4\theta^{-1} \mathrm{E}\left[\left\|\frac{1}{n^2}\sum_{i\neq j}\mathrm{vec}\left( X_i X_j'\right)\varepsilon_{ij} \right\|_p^2\right],
	\end{align}
	where $\theta > 0$ is arbitrary. Conditional on $X_1, \ldots, X_n$, the summands $\mathrm{vec}\left( X_i X_j'\right)\varepsilon_{ij}$ are independent with mean zero. Thus, by Theorem 2.2 in~\cite{duembgen2010nemirovskis} and eq.~\eqref{eq:lemma:BoundsCovariance-HeavyTailed-6} there exists an absolute constant $C_5 > 0$ such that
	\begin{align*}
	\mathrm{E}\left[\left\|\frac{1}{n^2}\sum_{i\neq j}\mathrm{vec}\left( X_i X_j'\right)\varepsilon_{ij} \right\|_p^2\right] \leq C_5 (p \wedge \log d)\frac{1}{n^4}\sum_{i\neq j}\mathrm{E}\left[\|\mathrm{vec}\big(X_iX_j'\big)\varepsilon_i\|_p^2 \right] \leq C_4 C_5 K_p^4\|\sigma\|_p^4\left( \frac{p \wedge \log d}{n^2} \right).
	\end{align*}
	Hence, we set $\theta = M (C_4 C_5 K_p^4 \vee 64) \|\sigma\|_p^4 \left(\frac{ p \wedge \log d}{n^2}\right)$ and $t = M^{1/2}\theta^{1/2} \geq 32 n^{-1} C_4^{1/2} K_p^2 \|\sigma\|_p^2$, where $M \geq 1$ is a large absolute constant. By choosing $M$ large enough we can make the left hand side of eq.~\eqref{eq:lemma:BoundsCovariance-HeavyTailed-10} arbitrarily small, i.e.
	\begin{align}\label{eq:lemma:BoundsCovariance-HeavyTailed-11}
	\left\|\mathrm{vec}\left(\frac{1}{n^2}\sum_{i\neq j} X_i X_j'\right) \right\|_p = O_p\left( K_p^2 \|\sigma\|_p^2 \sqrt{\frac{ p \wedge \log d}{n^2}}\right).
	\end{align}	
	Lastly, by triangle inequality, eq.~\eqref{eq:lemma:BoundsCovariance-HeavyTailed-5} and eq.~\eqref{eq:lemma:BoundsCovariance-HeavyTailed-1} we have
	\begin{align}\label{eq:lemma:BoundsCovariance-HeavyTailed-12}
	\begin{split}
	\left\|\mathrm{vec}\left(\frac{1}{n^2}\sum_{i=1}^n X_i X_i'\right) \right\|_p &\leq  \left\|\frac{1}{n^2}\sum_{i=1}^n \mathrm{vec}\left(X_i X_i' - \mathrm{E}[X_iX_i']\right) \right\|_p + \frac{1}{n^2}\sum_{i=1}^n\mathrm{E}\left[\| \mathrm{vec}(X_iX_i')\|_p^2\right]^{1/2}\\
	&=O_p\left( K_{p \vee 4}^2 \|\sigma\|_p^2 \sqrt{\frac{ p \wedge \log d}{n^2}}\right) + O\left(\frac{K_{p \vee 4}^2\|\sigma\|_p^2}{n}\right).
	\end{split}
	\end{align}
	Combine eq.~\eqref{eq:lemma:BoundsCovariance-HeavyTailed-11} and eq.~\eqref{eq:lemma:BoundsCovariance-HeavyTailed-12} to conclude that
	\begin{align}\label{eq:lemma:BoundsCovariance-HeavyTailed-13}
	\left\|\mathrm{vec}\left(\mathbf{II}\right) \right\|_p = O_p\left( K_{p \vee 4}^2 \|\sigma\|_p^2 \sqrt{\frac{ p \wedge \log d}{n^2}}\right).
	\end{align}
	Therefore,
	\begin{align}\label{eq:lemma:BoundsCovariance-HeavyTailed-14}
	\|\mathrm{vec}(\widehat{\Sigma} - \Sigma)\|_p \lesssim \|\mathrm{vec}(\mathbf{I})\|_p \vee \|\mathrm{vec}(\mathbf{II})\|_p = O_p\left( K_{p \vee 4}^2 \|\sigma\|_p^2 \sqrt{\frac{ p \wedge \log d}{n}}\right).
	\end{align}
	
	\textbf{Proof of Case (ii).} Since $\frac{1}{n^2}\sum_{1\leq , j\leq n}^n X_iX_j' = \left(\frac{1}{n}\sum_{i=1}^n X_i\right) \left(\frac{1}{n}\sum_{i=1}^n X_i\right)'$ has rank one, we have
	\begin{align*}
	\|\widehat{\Sigma} - \Sigma\|_{op} &\leq \left\|\frac{1}{n} \sum_{i=1}^n \big(X_iX_i' - \mathrm{E}[X_i X_i']\big)\right\|_{op} + \left\|\frac{1}{n^2}\sum_{1\leq i,j \leq n } X_i X_j'\right\|_{op}\\
	&=  \left\|\frac{1}{n} \sum_{i=1}^n \big(X_iX_i' - \mathrm{E}[X_i X_i']\big)\right\|_{op} + \left\|\frac{1}{n^2}\sum_{1\leq , j\leq n}^n \mathrm{vec}(X_iX_j') \right\|_2\\
	&\equiv \|\mathbf{I}\|_{op} + \|\mathrm{vec}(\mathbf{II})\|_2.
	\end{align*}
	By Theorem 5.48 in~\cite{vershynin2011introduction},
	\begin{align}
	\|\mathbf{I}\|_{op} =  O_p\left(\|\Sigma\|_{op}\left(\sqrt{\frac{\mathrm{m}(\Sigma)\log(d \wedge n)}{n}} \vee \frac{\mathrm{m}(\Sigma)\log(d \wedge n)}{n} \right) \right).
	\end{align}
	Since $\|\sigma\|_2^2 = \mathrm{tr}(\Sigma) \leq \mathrm{m}(\Sigma)\|\Sigma\|_{op}$, we have by eq.~\eqref{eq:lemma:BoundsCovariance-HeavyTailed-13} with $p=2$,
	\begin{align*}
	\|\mathrm{vec}(\mathbf{II})\|_2 \lesssim  O_p\left( \|\Sigma\|_{op} \frac{\mathrm{m}(\Sigma)}{n}\right).
	\end{align*}
	The claim follows from combining the last two bounds.
	
	\textbf{Case (iii).} Denote by $\oslash$ the Hadamard division. Suppose that Assumption~\ref{assumption:FiniteMoments} holds with $s \geq 4$. Observe that
	\begin{align*}
	\max_{1 \leq k \leq d} \big|(\widehat{\sigma}_k/\sigma_k)^2 - 1\big| &= \left\|\mathrm{diag}(\widehat{\Sigma})\oslash \mathrm{diag}(\Sigma) - I_d\right\|_{\infty}\\
	&\leq \left\|\mathrm{diag}\left( \left(\frac{1}{n} \sum_{i=1}^n X_iX_i' - \mathrm{E}[X_i X_i']\right)\oslash\mathrm{diag}(\Sigma) - I_d\right)\right\|_\infty\\
	&\quad{} + \left\| \mathrm{diag}\left( \left(\frac{1}{n^2}\sum_{1\leq i,j \leq d } X_i X_j'\right) \oslash \mathrm{diag}(\Sigma) \right)\right\|_\infty\\
	&\leq \left\|\mathrm{diag}\left( \left(\frac{1}{n} \sum_{i=1}^n X_iX_i' - \mathrm{E}[X_i X_i']\right)\oslash\mathrm{diag}(\Sigma) - I_d\right)\right\|_r\\
	&\quad{} + \left\| \mathrm{diag}\left( \left(\frac{1}{n^2}\sum_{1\leq i,j \leq d } X_i X_j'\right) \oslash \mathrm{diag}(\Sigma) \right)\right\|_r\\
	\end{align*}
	Moreover, $\|\mathrm{diag}(I_d)\|_{\infty} = 1$. Hence, from Case 1,
	\begin{align*}
	\max_{1 \leq k \leq d} \big|(\widehat{\sigma}_k/\sigma_k)^2 - 1\big| = O_p\left( K_s^2 d^{1/s} \sqrt{\frac{ s \wedge \log d}{n}}\right).
	\end{align*}
	Suppose Assumption 2 holds with $r\geq 2$. Observe that
	\begin{align*}
	\max_{1 \leq k \leq d} \big|(\widehat{\sigma}_k/\sigma_k)^2 - 1\big| &= \left\|\mathrm{diag}(\widehat{\Sigma})\oslash \mathrm{diag}(\Sigma) - I_d\right\|_{op}.
	\end{align*}
	Note that $\mathrm{tr}\big(\mathrm{diag}(\Sigma) \oslash \mathrm{diag}(\Sigma)\big) \leq \widetilde{\mathrm{m}}\big(\mathrm{diag}(\Sigma)\big)$. Thus, by Case 2,
	\begin{align*}
	\max_{1 \leq k \leq d} \big|(\widehat{\sigma}_k/\sigma_k)^2 - 1\big| &=  O_p\left(\sqrt{\frac{\widetilde{\mathrm{m}}\big(\mathrm{diag}(\Sigma)\big) \log(d \wedge n)}{n}} \vee \frac{\widetilde{\mathrm{m}}\big(\mathrm{diag}(\Sigma)\big) \log(d \wedge n)}{n} \right).
	\end{align*}
\end{proof}
\begin{proof}[\textbf{Proof of Lemma~\ref{lemma:BoundsThresholdEstimators}}]
	\textbf{Proof of Case (i).} By Lemma~\ref{lemma:BoundsCovariance} (i) with probability at least $1- 2\zeta$,
	\begin{align*}
	\|\mathrm{vec}(\widehat{\Sigma} - \Sigma)\|_\infty \lesssim \|\sigma\|_\infty^2\lambda_n.
	\end{align*}
	On this event, it is straight forward to show~\citep[e.g.][p. 181]{wainwright2019high} that
	\begin{align*}
	\big\|\mathrm{vec}\big(T_{\lambda_n}(\widehat{\Sigma}) - \Sigma\big)\big\|_p \lesssim \|\mathrm{vec}(A)\|_p \|\sigma\|_\infty^2 \lambda_n.
	\end{align*}
	and
	\begin{align*}
	\big\|T_{\lambda_n}(\widehat{\Sigma}) - \Sigma\big\|_{op} \lesssim \|\mathrm{vec}(A)\|_{op}\|\sigma\|_\infty^2  \lambda_n.
	\end{align*}
	
	\textbf{Proof of Case (ii).} The claim about the difference in operator norm follows verbatim from the proof of Theorem 6.27 in~\cite{wainwright2019high}. The statement about the difference in vectorized $\ell_p$-norm follows from an easy modification of the proof of Theorem 6.27. For completeness we provide a sketch of the modified argument. Wainright's proofs are easier to generalize to our setup than the original proofs in~\cite{bickel2008regularization}.
	
	Suppose that $\|\mathrm{vec}(\widehat{\Sigma} - \Sigma)\|_\infty \leq \lambda/2$ holds, where $\lambda > 0$ will be specified below. Fix $j \in \{1, \ldots, d\}$ and define
	\begin{align*}
	S_j(\lambda/2) : = \big\{ k \in \{1, \ldots, d\}: |\Sigma_{jk}| > \lambda/2\big\}.
	\end{align*}
	For any $k \in S_j(\lambda/2)$, we have
	\begin{align}\label{eq:lemma:BoundsThresholdEstimators-1}
	|T_\lambda(\widehat{\Sigma}_{jk}) - \Sigma_{jk}| \leq |T_\lambda(\widehat{\Sigma}_{jk}) - \widehat{\Sigma}_{jk}|  + |\widehat{\Sigma}_{jk}-\Sigma_{jk}| \leq \frac{3}{2}\lambda,
	\end{align}
	where the second inequality follows from $\|\mathrm{vec}(\widehat{\Sigma} - \Sigma)\|_\infty \leq \lambda/2$ and property (iii) of the thresholding operator $T_\lambda$.	For any $k \notin S_j(\lambda/2)$, it follows from $\|\mathrm{vec}(\widehat{\Sigma} - \Sigma)\|_\infty \leq \lambda/2$ and property (ii) of the thresholding operator $T_\lambda$ that $T_\lambda(\widehat{\Sigma}_{jk}) = 0$. Hence,
	\begin{align}\label{eq:lemma:BoundsThresholdEstimators-2}
	|T_\lambda(\widehat{\Sigma}_{jk}) - \Sigma_{jk}| = |\Sigma_{jk}|.
	\end{align}
	Combine eq.~\eqref{eq:lemma:BoundsThresholdEstimators-1} and~\eqref{eq:lemma:BoundsThresholdEstimators-2} to conclude that
	\begin{align}\label{eq:lemma:BoundsThresholdEstimators-3}
	\left(\sum_{j=1}^d|T_\lambda(\widehat{\Sigma}_{jk}) - \Sigma_{jk}|^p\right)^{1/p} &\leq \left(\sum_{j \in S_k(\lambda/2)}|T_\lambda(\widehat{\Sigma}_{jk}) - \Sigma_{jk}|^p\right)^{1/p} + 	\left(\sum_{j \notin S_k(\lambda/2)}|T_\lambda(\widehat{\Sigma}_{jk}) - \Sigma_{jk}|^p\right)^{1/p} \nonumber\\
	&\leq |S_k(\lambda/2)| \frac{3}{2}\lambda + \left(\sum_{j \notin S_k(\lambda/2)}|\Sigma_{jk}|^p\right)^{1/p}.
	\end{align}
	To bound the first term on the far right hand side in above display note that
	\begin{align*}
	R_{p,\gamma} \geq \left(\sum_{j=1}^d |\Sigma_{jk}|^{p\gamma}\right)^{1/p} \geq |S_k(\lambda/2)| (\lambda/2)^\gamma.
	\end{align*}
	Re-arranging this inequality yields $|S_k(\lambda/2)| \leq 2^\gamma R_{p,\gamma}\lambda^{-\gamma}$. To bound the second term on the far right hand side in eq.~\eqref{eq:lemma:BoundsThresholdEstimators-3} observe that
	\begin{align*}
	\left(\sum_{j \notin S_k(\lambda/2)}|\Sigma_{jk}|^p\right)^{1/p} = \frac{\lambda}{2} \left(\sum_{j \notin S_k(\lambda/2)}\left|\frac{\Sigma_{jk}}{\lambda/2}\right|^p\right)^{1/p} \leq \frac{\lambda}{2} \left(\sum_{j \notin S_k(\lambda/2)}\left|\frac{\Sigma_{jk}}{\lambda/2}\right|^{p\gamma}\right)^{1/p} \leq \lambda^{1-\gamma}R_{p,\gamma}.
	\end{align*}
	Combine the preceding two inequalities with eq.~\eqref{eq:lemma:BoundsThresholdEstimators-3} and conclude that
	\begin{align}\label{eq:lemma:BoundsThresholdEstimators-4}
	\left(\sum_{j=1}^d|T_\lambda(\widehat{\Sigma}_{jk}) - \Sigma_{jk}|^p\right)^{1/p} \leq 2^\gamma R_{p,\gamma} \lambda^{1-\gamma} \frac{3}{2} + R_{p,\gamma}\lambda^{1-\gamma} \leq 4 R_{p,\gamma} \lambda^{1-\gamma}.
	\end{align}
	We now determine the choice of $\lambda > 0$. Consider the following,
	\begin{align*}
	&\mathrm{P}\left(\|\mathrm{vec}(T_\lambda(\widehat{\Sigma}) - \Sigma)\|_p > 4d^{1/p} R_{p,\gamma} \lambda^{1-\gamma} \right)\\
	&= \mathrm{P}\left(\|\mathrm{vec}(T_\lambda(\widehat{\Sigma}) - \Sigma)\|_p > 4d^{1/p} R_{p,\gamma} \lambda^{1-\gamma}, \:\: \|\mathrm{vec}(\widehat{\Sigma} - \Sigma)\|_\infty \leq \lambda/2 \right)\\
	&\quad{}  + \mathrm{P}\left(\|\mathrm{vec}(\widehat{\Sigma} - \Sigma)\|_\infty \leq \lambda/2 \right) \\
	&\leq \sum_{k =1}^d \mathrm{P}\left(\left(\sum_{j=1}^d|T_\lambda(\widehat{\Sigma}_{jk}) - \Sigma_{jk}|^p\right)^{1/p} > 4R_{p,\gamma} \lambda^{1-\gamma}, \:\: \|\mathrm{vec}(\widehat{\Sigma} - \Sigma)\|_\infty \leq \lambda/2 \right)\\
	&\quad{}+ \mathrm{P}\left(\|\mathrm{vec}(\widehat{\Sigma} - \Sigma)\|_\infty \leq \lambda/2 \right)\\
	&=\mathrm{P}\left(\|\mathrm{vec}(\widehat{\Sigma} - \Sigma)\|_\infty \leq \lambda/2 \right),
	\end{align*}
	where the last line follows from eq.~\eqref{eq:lemma:BoundsThresholdEstimators-4}. Now, set $\lambda = 2 \|\sigma\|_\infty^2\lambda_n$ and conclude by Lemma~\ref{lemma:BoundsCovariance} (i) (applied with $p = \infty$) that with probability at least $1- 2\zeta$,
	\begin{align*}
	\|\mathrm{vec}(T_\lambda(\widehat{\Sigma}) - \Sigma)\|_p \lesssim d^{1/p} R_{p,\gamma} \|\sigma\|_\infty^{2(1-\gamma)} \lambda_n^{1-\gamma}.
	\end{align*}
	
	\textbf{Proof of Case (iii).} Note that for all $s < \infty$ and $\lambda > 0$,  $\|\mathrm{vec}(\widehat{\Sigma} - \Sigma)\|_s \leq \lambda$ implies $\|\mathrm{vec}(\widehat{\Sigma} - \Sigma)\|_\infty \leq \lambda$. Thus, by Lemma~\ref{lemma:BoundsCovariance-FiniteMoments} (i) we have for all $s \geq ( p \wedge \log p) \vee 4$,
	\begin{align*}
	&\big\|\mathrm{vec}(T_{\lambda_n}(\widehat{\Sigma}) - \Sigma)\big\|_p = O_p\left(\|\mathrm{vec}(A)\|_p  K_s^2  \|\sigma\|_s^2 \sqrt{\frac{ s \wedge \log d}{n}}\right),\\
	&\big\|T_{\lambda_n}(\widehat{\Sigma}) - \Sigma\big\|_{op} =  O_p\left(\|A\|_{op} K_s^2 \|\sigma\|_s^2 \sqrt{\frac{ s\wedge \log d}{n}}\right).
	\end{align*}
	
	\textbf{Proof of Case (iv).} The claim follows as Case 2 but using Lemma~\ref{lemma:BoundsCovariance-FiniteMoments} (i) with $s \geq ( p \wedge \log p) \vee 4$,, instead of Lemma~\ref{lemma:BoundsCovariance} (i) with $p = \infty$. See also Case 3.
\end{proof}
\begin{proof}[\textbf{Proof of Lemma~\ref{lemma:BoundsBandedEstimators}}]
	\textbf{Proof of Case (i).} The claim about the difference in operator norm follows from the proof of Theorem 1 in~\cite{bickel2008covariance}. The statement about the difference in vectorized $\ell_p$-norm follows from an easy modification of the proof of Theorem 1. For completeness we give the modified argument below.
	
	Fix $ \ell \in \{0, \ldots, d - 1\}$ and compute
	\begin{align}\label{eq:lemma:BoundsBandedEstimators-1}
	&\|\mathrm{vec}(B_{\ell}(\widehat{\Sigma}) - \Sigma)\|_p\nonumber\\
	&\quad{}\leq \big\|\mathrm{vec}\big(B_{\ell}(\widehat{\Sigma}) - B_{\ell}(\Sigma)\big)\big\|_p + \|\mathrm{vec}(B_{\ell}(\Sigma) - \Sigma)\|_p \nonumber\\
	&\quad{}= \left(\sum_{j=1}^d \sum_{k=1}^d |\widehat{\Sigma}_{jk} - \Sigma_{jk}|^p\mathbf{1}\{|j-k| \leq \ell\}\right)^{1/p} + \left(\sum_{j=1}^d \sum_{k=1}^d |\Sigma_{jk}|^p\mathbf{1}\{|j-k| > \ell\}\right)^{1/p} \nonumber\\
	&\quad{}\leq \big( d + \ell (2d - \ell - 1)\big)^{1/p}\big\|\mathrm{vec}\big(B_{\ell}(\widehat{\Sigma}) - B_{\ell}(\Sigma)\big)\big\|_{\infty} + d^{1/p} B_p \ell^{-\alpha},
	\end{align}
	where the first term on the far right hand side follows since the double sum has $d + \ell(2d - \ell -1)$ nonzero summands and the second term follows from Assumption 4.
	
	Set $\ell = \ell_n \equiv B_p^{p/(1 + p\alpha)} \|\sigma\|_\infty^{-2p/(1 + p\alpha)}\lambda_n^{-p/(1 + p\alpha)}$. By Lemma~\ref{lemma:BoundsCovariance} (i) (applied to $p = \infty$), with probability at least $1- 2\zeta$,
	\begin{align}\label{eq:lemma:BoundsBandedEstimators-2}
	\|\mathrm{vec}(B_{\ell}(\widehat{\Sigma}) - \Sigma)\|_p \lesssim B_p^{1/(1 + p\alpha)} d^{1/p} \|\sigma\|_\infty^{2p\alpha/(1+p\alpha)}  \lambda_n^{p\alpha/(1 + p\alpha)}.
	\end{align}
	
	\textbf{Proof of Case (ii).} The claim follows as Case 1 but using Lemma~\ref{lemma:BoundsCovariance-FiniteMoments} (i) with $s \geq (p \wedge \log d) \vee 4$, instead of Lemma~\ref{lemma:BoundsCovariance} (i) with $p = \infty$. See also Case 1.
\end{proof}

\subsubsection{Proofs for Appendix~\ref{subsec:AuxResultsApplication}}
\begin{proof}[\textbf{Proof of Lemma~\ref{lemma:ComparisonQuantiles}}]
	The proof is identical to the one of Lemma 3.2 in~\cite{chernozhukov2013GaussianApproxVec}. We sketch it for completeness. By Theorem~\ref{theorem:Gaussian-Comparison-KolmogorovDistance}, on the event $\{\Pi_p \leq \delta \}$, we have
	$|\mathrm{P}(S_{n,p}^* \leq t \mid X) - \mathrm{P}(\widetilde{S}_p \leq t)| \leq \pi_p(\delta)$ or all $ t \in \mathbb{R}$; in particular, for $t = \tilde{c}_{n,p}\big(\pi_p(\delta) + \alpha\big)$ we have
	\begin{align*}
	\mathrm{P}\Big(S_{n,p}^* \leq \tilde{c}_p\big(\pi_p(\delta) + \alpha\big) \mid X\Big) &\geq 	\mathrm{P}\Big(\widetilde{S}_p \leq \tilde{c}_p\big(\pi_p(\delta) + \alpha\big) \mid X\Big) - \pi_p(\delta) \\
	&\geq \pi_p(\delta) + \alpha - \pi_p(\delta) = \alpha.
	\end{align*}
	This implies the first inequality in the lemma. The second follows similarly.
\end{proof}
\begin{proof}[\textbf{Proof of Lemma~\ref{lemma:UpperBoundQuantiles}}]
	We first establish the upper bound for all $\alpha \in (0,1)$.Note that $\widetilde{S}_p := \|\Omega^{1/2}Z\|_p$ with $Z \sim N(0, I)$ and that the map $f(Z)=\|\Omega^{1/2}Z\|_p$ is Lipschitz continuous (with respect to the Euclidean norm) with Lipschitz constant $\|\Omega^{1/2}\|_{2 \rightarrow p} := \sup_{\|u\|_2 \leq 1} \|\Omega^{1/2}u\|_p$. Thus, by the Gaussian concentration inequality for Lipschitz continuous functions~\citep[e.g.][Lemma A.2.2]{vandervaart1996weak}, for all $t > 0$,
	\begin{align*}
	\mathrm{P}\left(\widetilde{S}_p - \mathrm{E}[\widetilde{S}_p] \geq t \right) \leq \exp\left\{-\frac{t^2}{2\|\Omega^{1/2}\|_{2 \rightarrow p}^2}\right\}.
	\end{align*}
	In particular,
	\begin{align*}
	\mathrm{P}\left(\widetilde{S}_p > \mathrm{E}[\widetilde{S}_p]  + \sqrt{2\log(1/\alpha)} \|\Omega^{1/2}\|_{2 \rightarrow p}\right) \leq \alpha.
	\end{align*}
	Similarly, by Chebyshev's inequality
	\begin{align*}
	\mathrm{P}\left(\widetilde{S}_p - \mathrm{E}[\widetilde{S}_p] \geq t \right) \leq \frac{\mathrm{Var}[\widetilde{S}_p]}{t^2},
	\end{align*}
	and therefore
	\begin{align*}
	\mathrm{P}\left(\widetilde{S}_p > \mathrm{E}[\widetilde{S}_p]  + \sqrt{1/\alpha}\sqrt{\mathrm{Var}[\widetilde{S}_p]} \right) \leq \alpha.
	\end{align*}
	Now, the upper bound follows from the definition of $\tilde{c}_p(1-\alpha)$.
	
	To establish the lower bound for $\alpha \in (0,1/2]$, recall the following inequality:
	\begin{align*}
	\left|\mathrm{E}[\widetilde{S}_p] - \tilde{c}_p(1/2)\right| \leq \sqrt{\mathrm{Var}[\widetilde{S}_p]}.
	\end{align*}
	Whence, for all $\alpha \in (0,1/2]$ it follows that
	\begin{align*}
	\tilde{c}_p(1-\alpha)  \geq \tilde{c}_p(1/2) \geq \mathrm{E}[\widetilde{S}_p] - \sqrt{\mathrm{Var}[\widetilde{S}_p]}.
	\end{align*}
	To conclude, note that by the Gaussian Poincar{\'e} inequality, $\sqrt{\mathrm{Var}[\widetilde{S}_p]} \leq \|\Omega^{1/2}\|_{2 \rightarrow p}$.
\end{proof}

\subsubsection{Proofs for Appendix~\ref{subsec:AuxResultsDerivatives}}
\begin{proof}[\textbf{Proof of Lemma~\ref{lemma:Lp-Norm-Derivatives}}]
	The claim follows from straightforward computations. The most convenient way to carry out those calculations is to notice that $M_p(x) \equiv M_{f, g}(x) = g^{-1}\left(\sum_{j=1}^d f(x_j)\right)$ for $f(x) = g(x) = x^p$, $x \geq 0$. Now, repeated applications of the implicit function theorem and the chain rule yield the claim.
\end{proof}
\begin{proof}[\textbf{Proof of Lemma~\ref{lemma:Lp-Norm-Derivatives-Stability}}]
	Note that $M_p(x) = \|x\|_p$ for any $x \in \big\{z \in \mathbb{R}^d: z_i \geq 0, i=1, \ldots, d \big\} \setminus \{0\}$ and $p > 1$. With slight abuse of notation, we will also use this formulation when the exponent is less than one or negative. First, since conjugate exponents satisfy $(p-1)q = p$,
	\begin{align*}
	\sum_{k=1}^d \left|\frac{\partial M_p(x)}{\partial x_k}\right|^q = \frac{\|x\|_{(p-1)q}^{(p-1)q}}{\|x\|_p^{(p-1)q}} = \frac{\|x\|_p^p}{\|x\|_p^p} = 1.
	\end{align*}
	Second, suppose that $p \geq 2$. Since $(p-2)q = p - q$ and $(2p-1)q = 2p + q$,
	\begin{align*}
	\sum_{k=1}^d \left|\frac{\partial^2 M_p(x)}{\partial x_k^2}\right|^q &\leq \frac{2^{q-1}(p-1)^q\|x\|_{(p-2)q}^{(p-2)q}}{\|x\|_p^{(p-1)q}} + \frac{2^{q-1}(p-1)^q\|x\|_{2(p-1)q}^{2(p-1)q}}{\|x\|_p^{(2p-1)q}}\\
	&= \frac{2^{q-1}(p-1)^q\|x\|_{p-q}^{p-q}}{\|x\|_p^p} + \frac{2^{q-1}(p-1)^q\|x\|_{2p}^{2p}}{\|x\|_p^{2p + q}}\\
	&\leq \frac{2^{q-1}(p-1)^q d^{q/p}\|x\|_p^{p-q}}{\|x\|_p^p} + \frac{2^{q-1}(p-1)^q\|x\|_{2p}^{2p}}{\|x\|_p^{2p + q}}\\
	&= \frac{2^{q-1}(p-1)^q d^{q/p}}{\|x\|_p^q} + \frac{2^{q-1}(p-1)^q}{\|x\|_p^q},
	\end{align*}
	where the second inequality follows from the power mean inequality. Third, since $(2p-1)q = 2p + q$,
	\begin{align*}
	\sum_{k, \ell} \left|\frac{\partial^2 M_p(x)}{\partial x_k \partial x_\ell}\right|^q \leq \frac{(p-1)^q \|x\|_{(p-1)q}^{2(p-1)q}}{\|x\|_p^{(2p-1)q}} = \frac{(p-1)^q \|x\|_p^{2p}}{\|x\|_p^{2p + q}} = \frac{(p-1)^q}{\|x\|_p^q}.
	\end{align*}
	Fourth, since $(3p-1)q = 3p + 2q$,
	\begin{align*}
	\sum_{k, \ell, m} \left|\frac{\partial^3 M_p(x)}{\partial x_k \partial x_\ell \partial x_m}\right|^q &\leq \frac{(2p-1)^q(p-1)^q\|x\|_{(p-1)q}^{3(p-1)q}}{\|x\|_p^{(3p-1)q}} = \frac{(2p-1)^q(p-1)^q\|x\|_p^{3p}}{\|x\|_p^{3p + 2q}} = \frac{(2p-1)^q(p-1)^q}{\|x\|_p^{2q}}.
	\end{align*}
	Fifth, suppose that $p \geq 2$, and compute
	\begin{align*}
	\sum_{k, \ell} \left|\frac{\partial^3 M_p(x)}{\partial x_k^2 \partial x_\ell}\right|^q &\leq \frac{2^{q-1}(p-1)^{2q} \|x\|_{(p-2)q}^{(p-2)q} \|x\|_{(p-1)q}^{(p-1)q}}{\|x\|_p^{(2p-1)q}} + \frac{2^{q-1}(2p-1)^q(p-1)^q\|x\|_{2(p-1)q}^{2(p-1)q}\|x\|_{(p-1)q}^{(p-1)q}}{\|x\|_p^{(3p-1)q}}\\
	&= \frac{2^{q-1}(p-1)^{2q} \|x\|_{p-q}^{p-q} \|x\|_p^p}{\|x\|_p^{2p+q}} + \frac{2^{q-1}(2p-1)^q(p-1)^q\|x\|_{2p}^{2p}\|x\|_p^p}{\|x\|_p^{3p+ 2q}}\\
	&\leq \frac{2^{q-1}(p-1)^{2q}d^{q/p} \|x\|_p^{p-q}}{\|x\|_p^{p+q}} + \frac{2^{q-1}(2p-1)^q(p-1)^q\|x\|_{2p}^{2p}}{\|x\|_p^{2p+ 2q}}\\
	&\leq \frac{2^{q-1}(p-1)^{2q}d^{q/p}}{\|x\|_p^{2q}} + \frac{2^{q-1}(2p-1)^q(p-1)^q}{\|x\|_p^{2q}},
	\end{align*}
	where the second inequality follows from the power mean inequality. Sixth, suppose that $p \geq 3$. Since $(p-3)q = p- 2q$ and $(2p- 3)q = 2p - q = p + (p-2)q \geq p$,
	\begin{align*}
	\sum_{k=1}^d \left|\frac{\partial^3 M_p(x)}{\partial x_k^3}\right|^q &\leq \frac{2^{2q-1}(p-1)^q(p-2)^q\|x\|_{(p-3)q}^{(p-3)q}}{\|x\|_p^{(p-1)q}} + \frac{2^{2q-1}3^q(p-1)^{2q}\|x\|_{(2p-3)q}^{(2p-3)q}}{\|x\|_p^{(2p-1)q}}\\
	&\quad{} + \frac{2^{2q-1}(2p-1)^q(p-1)^q\|x\|_{3(p-1)q}^{3(p-1)q}}{\|x\|_p^{(3p-1)q}}\\
	&=\frac{2^{2q-1}(p-1)^q(p-2)^q\|x\|_{p -2q}^{p-2q}}{\|x\|_p^p} + \frac{2^{2q-1}3^q(p-1)^{2q}\|x\|_{2p-q}^{2p-q}}{\|x\|_p^{2p + q}}\\
	&\quad{} + \frac{2^{2q-1}(2p-1)^q(p-1)^q\|x\|_{3p}^{3p}}{\|x\|_p^{3p + 2q}}\\
	&\leq \frac{2^{2q-1}(p-1)^q(p-2)^qd^{2q/p}\|x\|_p^{p-2q}}{\|x\|_p^p} + \frac{2^{2q-1}3^q(p-1)^{2q}\|x\|_p^{2p-q}}{\|x\|_p^{2p + q}}\\
	&\quad{} + \frac{2^{2q-1}(2p-1)^q(p-1)^q}{\|x\|_p^{2q}}\\
	&= \frac{2^{2q-1}(p-1)^q(p-2)^qd^{2q/p}}{\|x\|_p^{2q}} + \frac{2^{2q-1}3^q(p-1)^{2q}}{\|x\|_p^{2q}} + \frac{2^{2q-1}(2p-1)^q(p-1)^q}{\|x\|_p^{2q}},
	\end{align*}
	where the second inequality follows from the power mean inequality. If $p = 2$, then the first term vanishes, and we have
	\begin{align*}
	\sum_{k=1}^d \left|\frac{\partial^3 M_2(x)}{\partial x_k^3}\right|^2 &\leq \frac{2^{2q}3^q}{\|x\|_p^{2q}}.
	\end{align*}
\end{proof}
\begin{proof}[\textbf{Proof of Lemma~\ref{lemma:Lp-Norm-Derivatives-Stability-3}}]
	The claim follows from Lemma~\ref{lemma:Lp-Norm-Derivatives-Stability} and the power mean inequality.
\end{proof}

\end{appendices}

\bibliography{GBA_51_ref}
\bibliographystyle{apalike}
\end{document}